\documentclass[3p,final]{elsarticle}

\usepackage{pdf14}

\usepackage{mathpazo}
\usepackage[T1]{fontenc}
\newlength{\alphabet}
\usepackage[reqno,fleqn]{amsmath}
\usepackage{stmaryrd,xifthen,mathrsfs}
\usepackage{xspace}
\usepackage{xifthen}
\usepackage{verbatim}
\usepackage{ifdraft}

\usepackage{graphicx}
\graphicspath{{figs/}}
\usepackage{tabularx,booktabs,soul,multirow}
\usepackage[margin=10pt,font=small,textfont=it,labelfont=bf,justification=justified,position=bottom]{caption}
\usepackage[rightcaption]{sidecap}
\sidecaptionvpos{figure}{m}

\newcommand{\mcaption}[3]{
  \ifthenelse{\isempty{#2}}
             {\caption[#3]{#3 \label{#1}}}
             {\caption[#2]{{\sc #2.} #3 \label{#1}}}
}

\usepackage[format=hang,singlelinecheck=false,font={scriptsize},labelfont=bf,justification=centering]{subfig}

\newcommand\subfigure[2][]{\subfloat[#1]{#2}}

\newif\ifUseTikz
\ifUseTikz
\usepackage{tikz,pgfplots,import}
\pgfplotsset{compat=newest}
\pgfplotsset{plot coordinates/math parser=false}

\usetikzlibrary{arrows.meta}
\usetikzlibrary{backgrounds}
\usetikzlibrary{external}
\usetikzlibrary{pgfplots.groupplots}
\usetikzlibrary{plotmarks}

\pgfkeys{/pgf/images/include external/.code={\includegraphics[]{#1}}}
\newcommand{\includepgf}[2][1]{
  \beginpgfgraphicnamed{#2}%
  \tikzsetnextfilename{external-#2}%
  \scalebox{#1}{\subimport{figs/}{#2.pgf}}%
  \endpgfgraphicnamed%
}

\makeatletter
\tikzset{
    every picture/.style={
        execute at begin picture={
            \let\ref\@refstar
        }
    }
}
\makeatother

\newcommand\figscale{1}
\newcommand\figannotate{false}
\newcommand\pncangle{0}
\newcommand\pnclocx{0}
\newcommand\pnclocy{0}
\else%for ifUseTikz
\newcommand{\includepgf}[2][1]{
\scalebox{#1}{\includegraphics[]{external-#2}}%
}
\fi%for ifUseTikz
\newlength\figureheight
\newlength\figurewidth

\usepackage[
  pdftitle    = {Particulate Flow},
  pdfcreator  = {pdflatex},
  pdfkeywords = {},
  hyperindex  = {true},
  colorlinks  = {true},
  linkcolor   = {blue},
  citecolor   = {blue}
]{hyperref}

\usepackage{titlesec}
\titleformat{\subsubsection}[runin]{\normalfont\normalsize\itshape}{\thesubsubsection.}{.5em}{}[.\hspace{.5em}]
\titleformat{\paragraph}[runin]{\normalfont\normalsize\itshape}{\theparagraph}{}{}[.\hspace{.5em}]

\usepackage[capitalize]{cleveref}
\crefname{equation}{Eq.}{Eqs.}
\Crefname{equation}{Equation}{Equations}
\crefrangelabelformat{equation}{(#3#1#4--#5#2#6)}
\crefmultiformat{equation}{Eqs.~(#2#1#3)}{ and~(#2#1#3)}{, (#2#1#3)}{, and~(#2#1#3)}
\Crefmultiformat{equation}{Equations~(#2#1#3)}{ and~(#2#1#3)}{, (#2#1#3)}{, and~(#2#1#3)}
\newcommand\pr[1]{\cref{#1}}
\renewcommand\Pr[1]{\Cref{#1}}

\usepackage{wasysym} %biological symbols
\usepackage[normal]{engord} %remove normal option to have raised suffix
\usepackage{bm, paralist}
\usepackage{ctable}

\let\d\undefined
\let\O\undefined

\DeclareMathOperator{\Div}{\nabla\cdot}
\DeclareMathOperator{\Curl}{\nabla\times}
\newcommand\defeq{\mathrel{\mathop :}=}
\newcommand\inner[2]{#1\cdot #2}
\newcommand\cross[2]{#1\times #2}
\newcommand\d{\ensuremath{\,\mathrm{d}}}
\newcommand\parderiv[2]{\frac{\partial #1}{\partial #2}}
\newcommand\O[1]{\ensuremath{\mathcal{O}\left(#1\right)}}
\newcommand\ordinal[1]{\ifthenelse{\isin{#1}{abcdefghijklmnopqrstuvwxyz}}{\ensuremath{#1^\mathrm{th}}}{\engordnumber{#1}}}
\newcommand\R{\mathbb{R}}

\let\vector\undefined
\DeclareMathAlphabet{\mathbfsf}{\encodingdefault}{\sfdefault}{bx}{n}
\newcommand\discrete[1]{\mathsf{#1}}
\newcommand\fourier[1]{\widehat{#1}}
\newcommand\scalar[1]{#1}

\newcommand\vector[1]{\bm{#1}}
\newcommand\vectord[1]{\mathbfsf{#1}}%-- also change \discrete
\newcommand\linop[1]{\mathbf{#1}}
\newcommand\linopd[1]{\mathbfsf{#1}}
\newcommand\conv[1]{\mathcal{#1}}

\newcommand\contract[3][{\cdot}]{#2\mathop{#1}#3}

\newcommand\altr[3][-.15em]{#2/{\kern #1}#3}%alternative using a slash

\newcommand\expl[2][{}]{#2^{#1}}%explicit
\newcommand\impl[2][+]{#2^{\impltag{#1}}}%implicit
\newcommand\impltag[1]{#1}
\newcommand\lbl[1]{{#1}}
\newcommand\lh[1]{\textbf{#1}}%line heading
\newcommand\para[1]{\paragraph*{#1}}
\newcommand\uc[1][{}]{\ensuremath{\u^{\Gamma_0#1}}}
\newcommand\ucbar[1][{}]{\ensuremath{\bar{\u}^{\Gamma_0#1}}}
\newcommand\uf[1][{}]{\ensuremath{\u^{F#1}}}
\newcommand\ufbar[1][{}]{\ensuremath{\bar{\u}^{F#1}}}

\newcommand\up[1][{}]{\ensuremath{\u^{P#1}}}
\newcommand\upbar[1][{}]{\ensuremath{\bar{\u}^{P#1}}}
\newcommand\wlbl[1]{\mathrm{#1}}%word label

\def\LB/{Lattice Boltzmann} %not a combined name so no dash/if insisting on dash use -
\def\Ldot{\skew{6}{\dot}{L}}
\def\NS/{Navier--Stokes}
\def\U{\vector{U}}
\def\Xd{\vectord{X}}
\def\X{\vector{X}}

\def\bc{boundary condition}
\def\celegans{C. elegans}
\def\deltat{\delta t}
\def\emdash/{\kern 0.2em---\kern 0.2em}
\def\etal{\textit{et al.}}
\def\fidx{m}
\def\f{\vector{f}}
\def\gl{Gauss--Legendre}
\def\ns{nearly singular}
\def\nystrom{Nystr\"om}
\def\n{\vector{n}}
\def\pidx{n}
\def\rr{\vector{r}}
\def\u{\vector{u}}
\def\x{\vector{x}}
\def\y{\vector{y}}
\def\zhat{\hat{\vector{z}}}

\def\meter{\text{m}}
\def\microm{\text{$\mu$m}}
\def\minute{\text{min}}

\def\newton{\text{N}}
\def\pascal{\text{P{\kern -.09em}a}}
\def\pico{\text{p}}
\def\reynolds{\ensuremath{Re}}
\def\sec{\text{s}}
\def\usep{\kern .1em}

\begin{document}

%%- body
\title{A fast platform for simulating flexible fiber suspensions\\
  applied to cell mechanics
}
\author[nyu]{Ehssan Nazockdast} \ead{ehssan@cims.nyu.edu}
\author[nyu]{Abtin Rahimian} \ead{arahimian@acm.org}
\author[nyu]{Denis Zorin} \ead{dzorin@cs.nyu.edu}
\author[nyu]{Michael Shelley}\ead{shelley@cims.nyu.edu}
\address[nyu]{Courant Institute of Mathematical Sciences, New York
  University, New York, NY 10012}
\begin{abstract}
  We present a novel platform for the large-scale simulation of
  fibrous structures immersed in a Stokesian fluid and evolving under
  confinement or in free-space.
  One of the main motivations for this work is to study the dynamics
  of fiber assemblies within biological cells.
  For this, we also incorporate the key biophysical elements that
  determine the dynamics of these assemblies, which include the
  polymerization and depolymerization kinetics of fibers, their
  interactions with molecular motors and other objects, their
  flexibility, and hydrodynamic coupling.

  This work, to our knowledge, is the first technique to include
  many-body hydrodynamic interactions (HIs), and the resulting fluid
  flows, in cellular fiber assemblies.
  We use the non-local slender body theory to compute the
  fluid-structure interactions of the fibers and a second-kind
  boundary integral formulation for other rigid bodies and the
  confining boundary.
  A kernel-independent implementation of the fast multiple method is
  utilized for efficient evaluation of HIs.
  The deformation of the fibers is described by the nonlinear
  Euler--Bernoulli beam theory and their polymerization is modeled by
  the reparametrization of the dynamic equations in the appropriate
  non-Lagrangian frame.
  We use a pseudo-spectral representation of fiber positions and
  implicit HIs in the time-stepping to resolve large fiber
  deformations, and to allow time-steps not constrained by temporal
  stiffness or fiber-fiber interactions.
  The entire computational scheme is parallelized, which enables
  simulating assemblies of thousands of fibers.
  We use our method to investigate two important questions in the
  mechanics of cell division: (i) the effect of confinement on the
  hydrodynamic mobility of microtubule asters; and (ii) the dynamics
  of the positioning of mitotic spindle in complex cell geometries.
  Finally to demonstrate the general applicability of the method, we
  simulate the sedimentation of a cloud of fibers.
\end{abstract}
\maketitle

\section{Introduction\label{sec:intro}}
%%
%% -- motivation
%%

Semi-flexible biopolymers constitute a principal mechanical component
of intracellular structures \citep{BM2014, Mofrad2009, BK2006}.
Together with molecular motors, fiber networks, consisting of such
polymers,  form the cytoskeleton that is the cell's
mechanical machinery for executing several key tasks including cell
motility, material transport, and cell division \citep{Howard-Book}.
From these semi-flexible filaments (microtubules, actin filaments and
intermediate filaments) and a variety of molecular motors the
cytoskeleton is able to reorganize to supramolecular architectures
that are distinctly designed to perform a particular task
\citep{Fletcher2010, Howard-Book}.

Due to their central role in intracellular structures
the rheology and collective dynamics of semi-flexible fiber
suspensions and networks has attracted increasing interest in
engineering and biology \citep{BM2014}.

The basic physical difference between semi-flexible fibers and their
better-understood counterpart, polymer chains, is the significantly
larger bending rigidity of fibers that then yields larger end-to-end
distances. This results in many interesting differences between the
rheology of semi-flexible polymers and that of the two extreme limits
of flexibility, polymer melts and rigid fiber suspensions.  For
example, the stiffness of semi-flexible networks can increase or
decrease under compression and their suspensions show negative normal
stress differences \citep{Mofrad2009, BM2014, BS2001}.

With advancements in microscopy and data acquisition at small
time-~and~length-scales, we now know a great deal about the
interactions of the individual microscopic filaments and their
associated motor proteins. On the other hand, we know very little about how
they interact collectively and how these interactions determine the
ensemble behavior of cellular matter and structures.

Microrheological measurements provide a strong basis for understanding
the mechanical behavior of cytoskeletal structures and matter
\citep{Wirtz2009,Mofrad2009,Squires2010}, but do not directly inform
us of the relationship between microscopic interactions and
macroscopic behaviors. Moreover, living systems typically operate far
from equilibrium\emdash/due to internally generated forces being much
larger than thermal forces\emdash/and constitutive relationship is
required to extract rheological behavior based on, for example, the
trajectories of probe particles. Finding microscopic constitutive
relations %for exploiting active microrheology
is difficult even for the simplest of out-of-equilibrium complex fluids, a hard-sphere
colloidal suspension, due to the complex and nonlinear relations
between rheological properties and microstructural dynamics
\citep{NM2015}.

Dynamic simulation is a powerful tool to gain insight into the
underlying physical principles that govern the formation and
reorganization of cytoskeletal structures and ultimately obtain
relevant constitutive relationships.  With the continuous
advancements in \emph{in vitro} reconstitution of cellular matter,
comparing the experiments with \emph{in silico} reconstitution (i.e.,
detailed, large-scale, dynamic simulation of cellular structures) is
within reach \citep{BK2006}.  To this end, this paper presents a
computational platform for dynamic simulation of semi-flexible fiber
suspensions in Stokes flow.  Our method explicitly accounts for fiber
flexibility, their polymerization and depolymerization kinetics, their
interactions with molecular motors, and hydrodynamic interactions
(HIs). From a physical point of view, what distinguishes our method is
the inclusion of HIs, which has been almost entirely ignored in the
previous theoretical and numerical studies of cellular structures
\citep{BM2014}.

%%
%% -- method outline
%%
We consider suspensions of hydrodynamically interacting rigid
bodies and flexible fibers immersed in a Stokesian fluid, either under
confinement or in free-space. Our approach is based upon boundary
integral formulations of solutions to the Stokes equations. The flows
associated with the motion of rigid bodies and confining surfaces are
represented through a well-conditioned second-kind boundary integral
formulation \citep{Power1987}. The fluid flows associated with the
dynamics of fibers are accounted for using \emph{non-local slender
  body theory}~\citep{Keller1976, Johnson1980, gotz2000,
  Tornberg2004}.

%%
%% -- related work
%%
\paragraph{Related work}
Modeling approaches to suspensions and networks of fibers can be
roughly categorized into volume- and particle-based methods.  In
volume-based methods, the Stokes (or \NS/) equation is solved by
discretizing the entire computational domain.  Within this class,
immersed boundary methods have been applied to study the dynamics of
single \citep{Stockie1998,Lim2004} or several \citep{WS2015} flexible
fibers. The fibers are typically represented by a discrete set of
points (forming a one-dimensional curve \citep{Stockie1998, WS2015},
or a three-dimensional cylinder \citep{Lim2004}) whose interactions
capture stretching stiffness and internal elastic stresses.

These points on the fibers are Lagrangian and so are moved with the
background fluid flow. The consequent stretching or bending of the
discretized fiber creates elastic forces represented at the Lagrangian
points. These forces are distributed to the background grid, which
then provides forcing terms solving anew the Stokes or \NS/ equations
for the updated background flow. This cycle is then repeated.  A
similar update strategy has been adopted using the \LB/ method
\citep{Wu2010} to study the rheology of flexible fiber suspensions.
Typically, to properly account for fluid-structure interactions in
volume-based methods, the size of the volume grid is taken to be
several times smaller than the smallest dimension of the immersed
bodies.  As a result, these methods become computationally expensive
for simulating slender bodies such as fibers and disks. Moreover,
these methods typically  use explicit time-stepping
to evolve fiber shapes and elastic forces that substantially limits the region of
time-step stability \citep{LS2014} due to temporal stiffness.

Versions of particle-based methods include bead-spring models, regularized
Stokes, variations of dissipative particle dynamics
\citep{Gompper2011, Groot2012}, and slender-body theory (which we use
in this work).  In bead-spring models, the fiber is represented as a chain of
rigid spheres \citep{Yamamoto1995, Yamamoto1996, Fu2015} or ellipsoids
\citep{Ross1997} linked by inextensible connectors with finite bending
rigidity \citep{Joung2001} or by springs with given tensile and
bending stiffnesses \citep{Ross1997, Nedelec2007}.  The system is
evolved by imposing the balance between hydrodynamic and elastic
forces and torques on all beads.  Some implementations of this method
include long-range hydrodynamics as well as short-range lubrication
interactions \citep{Yamamoto1995, Yamamoto1996, Joung2001}, while
others only include local drag on the beads \citep{Ross1997}.  The
advantage of these techniques is their relatively simple
implementation.  Due to the discreteness of the fiber in this
construction, the polymerization of fibers is captured by adding and
removing beads and springs discretely in time \citep{Nedelec2007}.

Alternatively, the Regularized Stokes Method (RSM)
\cite{Cortez2001,Cortez2005} has been used to model the dynamics of
elastic fibers in Stokes flow \cite{Flores2005,Smith2010,Olson2013}.
The RSM represents the fluid velocity as the superposition of smoothed
fundamental solutions to the Stokes equation, with spread $\delta$,
distributed on immersed surfaces \cite{Cortez2005}.  To be convergent,
this method needs area elements to be scaled proportionally to
$\delta^2$ \cite{Cortez2005}.  Smith~\cite{Smith2010} removed this
constraint by formulating the integrals in the context of a boundary
element method.  Flores~\etal~\cite{Flores2005} used this framework to
study the role of hydrodynamic interactions on the dynamics of
flagella.  Olson~\etal~\cite{Olson2013} combined the elastic rod model
developed by \cite{Lim2004} with the RSM to study the dynamics of rods
with intrinsic curvature and twist.  In the context of fibers, the RSM
is similar to slender body theory but to achieve the same order of
accuracy in velocity as in slender body theory, the regularization
factor $\delta$ needs to be in the order of the aspect ratio of the
fibers. This becomes computationally demanding for slender objects.

In Boundary Integral (BI) methods, the Stokes equation is recast in
the form of integrals of distributions of point forces, torques, and
stresses on the boundaries of the fluid domain. As a result, the
computational domain is reduced to the two-dimensions of the immersed
surfaces. In \emph{Slender Body Theory} (SBT), the slenderness of the
fiber is used to asymptotically reduce the BI formulation to
one-dimensional integrals along the fibers' centerlines.  This results
in non-local slender body theory
(NLSBT)~\cite{Keller1976,Johnson1980,gotz2000, Tornberg2004}, which is
asymptotically accurate to $\O{\epsilon^{2}\log \epsilon}$, where
$\epsilon \ll 1$ is the aspect ratio of the
fiber. Shelley~and~Ueda~\cite{Shelley2000} designed numerical methods
for simulating closed flexible fibers based on NLSBT which
\cite{Tornberg2004} extended substantially to fibers with free ends
and devised algorithms to make the problem tractable and numerically
stable.

\paragraph{Contributions} We extend \citep{Tornberg2004} to enable simulation of
$\mathcal{O}(1000)$ actively driven semi-flexible filaments.  To
achieve this, we introduce several enhancements:
\begin{itemize}
\item \lh{Spectral spatial discretization:} A Chebyshev basis is used
  for the representation of fiber position enabling us to compute
  high-order derivatives with uniform accuracy along the fiber.
%  Accurate evaluation of high-order derivatives near the two fiber ends
%  is a challenging task in finite difference schemes with uniform
%  discretization due to degradation in polynomial interpolation accuracy.

\item \lh{Fiber (de)polymerization:} In biological settings the
  polymerization/depolymerization of biopolymers is a crucial part of
  their dynamics. Here we account for this dynamics by
  reparametrization of the dynamic equations into a non-Lagrangian
  frame. We found this formulation to be considerably more stable
  (at least in our setting)
  compared to introducing segments at discrete moments in time,
  as done in \citep{Nedelec2007}.

\item \lh{Removing numerical stiffness:} In the \cite{Tornberg2004}
  framework, to remove temporal stiffness, the bending forces were
  treated implicitly, while the tension equation\emdash/which imposed
  inextensibility of the fibers\emdash/was treated explicitly.
  We found that this formulation imposes severe
  limitations on the time-step magnitude for large numbers
  of fibers. In our method,  both bending and tensile
  forces are treated implicitly. As a result, the time-step in our scheme shows no
  dependency on the number of points per fiber and is only
  weakly dependent on the number of fibers.

\item \lh{Computational cost of $\mathcal{O}(N)$ and full
  parallelization:} Due to implicit treatment of the HIs and use of a
  BI method, at each time-step a dense system of equations must be
  solved.  We solve this system using GMRES with a Jacobi or block
  Jacobi preconditioner. A kernel independent fast multipole method
  (FMM) \citep{ying2004kernel,Malhotra2015} is utilized for fast computation of
  nonlocal hydrodynamic interactions, and fast matrix-vector
  products. The combination of GMRES, efficient preconditioning, and
  FMM results in $\mathcal{O}(N)$ computational cost per time-step,
  where $N$ is the number of unknowns, approximately proportional to
  the   number of fibers.  Due to implicit treatment of bending and
  tensile forces as well as the HIs, the stable time-step is three
  orders of magnitude larger than the stable time-step for the explicit
  method.  The entire computational scheme, including the FMM
  routine, the matrix-vector operations in GMRES, are parallelized and
  scalable to many computational cores.
\end{itemize}

To test our scheme and to demonstrate the variety of problems that can
be studied within this framework, we considered three representative
problems. We first investigate the mobility and viscoelastic behavior
of a spherical particle surrounded by a fibrous shell and compare our
results with analytical results obtained by using the porous medium
Brinkman model for the shell \cite{Brinkman1947, Masliyah1987}.
Through this example we clearly demonstrate that HIs play a critical
role in setting the dynamics of fibrous networks.

As an application of the current framework to cellular mechanics, we
study the effect of cell geometry on the positioning of the
\emph{pronuclear complex} (PNC) in the \emph{prophase} stage of cell
division in \celegans\ \citep{Cowan2004}.  For this purpose, we consider
three cell geometries and two different proposed force transduction
mechanisms for moving the PNC.  We demonstrate that changing the
geometry of the cell and the forcing mechanism both result in
substantial changes in the PNC positioning dynamics.

To demonstrate the utility of our method beyond biological settings,
we then look at the sedimentation of a cloud of flexible fibers, a
classical suspension mechanics problem. We find that the many-body
hydrodynamic interactions rearrange the fibers, resulting in evolution
of the cloud into a torus-like structure.

Our simulation results in all three problems are consistent with
available theoretical predictions and previous experimental and
simulation results. More importantly, our study of each problem
revealed several other interesting directions of research which can be
pursued within our framework.

\section{Formulation\label{sec:formulation}} Consider a suspension of
$N_F$ elastic fibers and $N_P$ rigid particles immersed in a Newtonian
fluid which is either confined by an outer boundary $\Gamma_0$ or
which fills free space. The effect of thermal fluctuations on the
fibers and other immersed bodies are neglected throughout this work.
In the context of intracellular assemblies, the confining boundary
represents a cell wall, and the Newtonian fluid filling the cell is
the cytoplasm that envelopes the assembly of the microscopic
filaments, nuclear complexes, and other organelles. A schematic is
shown in \pr{fig:schematic}.
\begin{figure}[!tb]
  \centering
  \includepgf{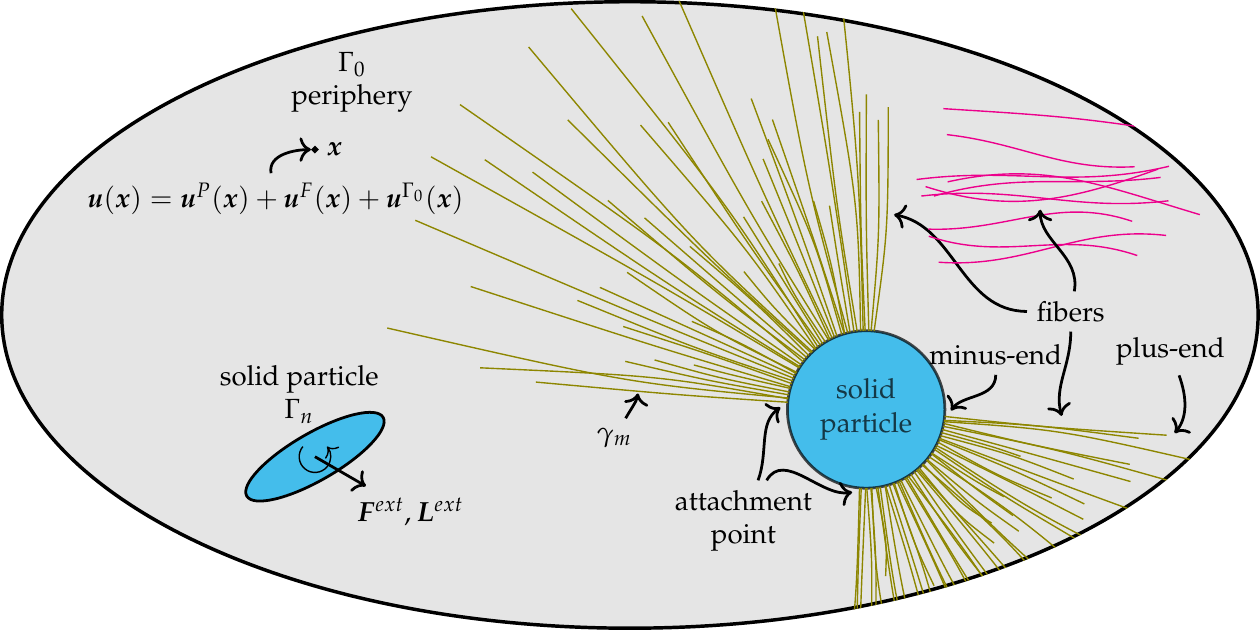}
  \mcaption{fig:schematic}{Schematic}{The gray-shaded region denotes
    the fluid domain $\Omega$.  The outer boundary is denoted by
    $\Gamma_0$.  There are rigid particles and flexible fibers
    suspended within the domain. The fibers and particles interact
    hydrodynamically with each other. Fibers may be attached to
    particles by different means that are specified through the \bc{s}
    on their ends (e.g. hinged or clamped).  When attached, fibers
    exert force $\vector{F}^\lbl{ext}$ and torque
    $\vector{L}^\lbl{ext}$ upon rigid particles.  Fibers can also move
    freely or sediment within the fluid domain.  The velocity at a
    point $\x\in\Omega$ results from the motion of particles $\up$,
    fibers $\uf$, and the backflow due to the confinement $\uc$.}
\end{figure}
The ratio of inertial to viscous forces is the Reynolds number,
$\reynolds=\rho L u/\mu$, where $\rho$ and $\mu$ are the cytoplasmic
density and viscosity, respectively, $L$ is the characteristic length
of the cell, and $u$ is an average velocity magnitude. For cellular
flows, due to the high viscosity of cytoplasm \citep{Howard-Book,
  Wirtz2009} and the small length-scale of the cell, $\reynolds \ll 1$
and so inertial effects can be safely neglected.  Hence, the flow of
cytoplasm is described by the incompressible Stokes equation
\begin{equation}
  \label{eqn:stokes}
    -\mu \Delta \u(\x) + \nabla p(\x) = \vector{0} \quad\text{and}\quad
    \Div\u(\x) = \scalar{0} \quad \text{for} \quad \x \in \Omega,
\end{equation}
where $\Omega$ denotes the domain occupied by the fluid.  Letting
$\Gamma_0$ denote the boundary of the domain and $\Gamma$ denote the
surface of the filaments and other particles immersed in the fluid,
the governing equations are augmented with the no-slip boundary
condition on the surface of these bodies
\begin{align}
  &\u(\x) = \vector{0} \quad\text{for}\quad\x \in \Gamma_0,\label{eqn:bc-cortex} \\
  &\u(\X) = \X_t       \quad\text{for}\quad\X \in \Gamma.
\label{eqn:noslip}
\end{align}
Throughout this paper, we use lowercase letters to denote Eulerian
variables, e.g. $\x$, and uppercase letters to denote Lagrangian
variables, e.g. $\X$.  Partial differentiation with respect to a
variable is denoted by a subscript, e.g. $\X_t\defeq\parderiv{\X}{t}$.
Thus,  $\X_t$ above  is the material surface velocity.
We denote the fiber centerline
positions by $\gamma_\fidx$ with $m=1, \dotsc, N_F$, and particle
surfaces by $\Gamma_\pidx$ with $\pidx=1, \dotsc, N_P$. $P$ 
is the union of all the immersed particle surfaces:
$P=\bigcup_{\pidx=1}^{N_P}\Gamma_\pidx$.

Using the fundamental solutions of the free-space Stokes equation,
\pr{eqn:stokes}, the fluid flow can be directly related to the
dynamics of immersed and bounding surfaces through a boundary integral
formulation. In particular, the solution of the Stokes boundary value
problem can be reformulated as solving a system of singular
integro-differential equations on all immersed and bounding surfaces
\cite{pozrikidis1992, Pozrikidis2001c}, which reduces the
computational domain from three dimensions to two.

The Oseen (Stokeslet) tensor $\linop{G}$, the Stresslet tensor
$\linop{T}$, and the Rotlet (Couplet) tensor $\linop{R}$, are
fundamental solutions of the Stokes equation, and are given by
\begin{align}
    \linop{G}(\rr) &= \frac{1}{8\pi \mu}
    \frac{\linop{I}+\hat{\rr}\hat{\rr}}{|\rr|}, \quad \hat{\rr} \defeq
    \frac{\rr}{|\rr|}, \label{eqn:oseen}\\ \linop{T}(\rr) & = -
    \frac{3}{4\pi \mu} \frac{\hat{\rr} \hat{\rr} \hat{\rr}}{
      |\rr|^2},\label{eqn:stresslet}\\ \contract{\linop{R}(\rr)}{\vector{L}}
    &= \frac{1}{2} \Curl \left( \contract{\linop{G}}{ \vector{L}}
    \right) = \frac{1}{8\pi\mu}\frac{\cross{\vector{L} }{ \hat{\rr}
    }}{|\rr|^2}\label{eqn:rotlet}.
\end{align}
In particular, the solution of
\pr{eqn:stokes,eqn:bc-cortex,eqn:noslip} is expressed as a convolution
of a vector density with the Stokeslet \altr{and}{or} Stresslet
tensors \cite{Pozrikidis2001c}. Therefore, we require convolutions of
these tensors along fibers and surfaces.  In particular, for a fiber
centerline $\gamma$ and a surface $S$ (periphery or an immersed
particle), we define
\begin{align}
  \conv{G}_\gamma [\f](\x) &\defeq \int_{\gamma}
  \d s_y~\contract{\linop{G}(\x-\y)}{\f(\y)},
  \label{eqn:single-layer}\\
  \conv{T}_S [\vector{q}](\x) &\defeq \int_S \d S_y~
  \contract{\n(\y)}{\contract{\linop{T}(\x-\y)}{\vector{q}(\y)}},
  \label{eqn:double-layer}
\end{align}
where $\n$ in \pr{eqn:double-layer} denotes the outward normal to the
surface $S$, and $\f$ and $\vector{q}$ are appropriately defined
vector densities. For $\x \in S$, the integral in
\pr{eqn:double-layer} is interpreted in the principal value
sense. Note that the integral contribution in \pr{eqn:single-layer} is
divergent if $\x\in\gamma$.

It is convenient to represent the fluid velocity as the superposition
of velocities arising from integral contributions from each surface
(periphery and immersed particles) and those from fiber centerlines:
\begin{align}
  \label{eqn:stokes-superpos}
  \vector{u}(\vector{x}) = \up(\x) + \uf(\x) + \uc(\x)
  \text{\quad with\quad}\up(\x)=\sum_{n=1}^{N_P}\up_\pidx (\x),
  \quad \uf(\x)=\sum_{m=1}^{N_F}\uf_{m} (\x),
\end{align}
where we have further decomposed the $N_P$ ($N_F$) immersed particles'
(fibers') velocity contribution, $\up$ ($\uf$), into those from each
individual immersed particle (fiber) with contribution $\up_\pidx$
($\uf_\fidx$).  It is also useful to define a ``complementary''
velocity field: for example, by $\ucbar$ we mean all the velocity
contributions other than those from the periphery surface $\Gamma_0$,
that is, $\ucbar=\u-\uc$ and similarly for $\upbar_{\pidx}$ and
$\ufbar_\fidx$.

We first outline the boundary integral formulation and flow
contributions arising from the bounding surface and immersed rigid
particles, based on the approach developed by \cite{Power1987}. We
then present a slender body formulation for the contributions of the
fibers \citep{Keller1976, Johnson1980, Tornberg2004}.  We proceed by
outlining the mechanics of elastic fibers, microtubule
(de)polymerization kinetics, and boundary conditions for fibers and
particles. We conclude this section with a summary of the formulation.

\subsection{The contribution from the periphery $\Gamma_0$
\label{ssc:cortex-flow}}
The fluid flow in the interior of the periphery can be written as a
double-layer boundary integral over $\Gamma_0$ with an unknown vector
density $\vector{q}_0$ \cite{Power1987}:
\begin{align}
  \uc(\x) &= \conv{T}_{\Gamma_0}[\vector{q}_0](\x) +
  \conv{N}_{\Gamma_0}[\vector{q}_0](\x), \text{\quad where\quad}
  \label{eqn:main-cortex-vel}\\ \conv{N}_{\Gamma_0} [\vector{q}_0](\x)
  &= \int_{\Gamma_0} \d S_y~\contract{[\n(\x)
      \n(\vector{y})]}{\vector{q}_0(\vector{y})},
\end{align}
and $\conv{T}_{\Gamma_0}$ is defined in \pr{eqn:double-layer}.
%The term $\conv{N}_{\Gamma_0}$ has been included for technical completeness
%(see below) but in this context of a fixed outer boundary can be set to zero.
Taking the limit of \pr{eqn:main-cortex-vel} as
$\x\rightarrow\Gamma_0$ and using the boundary condition,
\pr{eqn:bc-cortex}, generate a Fredholm integral equation of the
second kind for the unknown density $\vector{q}_0$
\begin{equation}
  -\frac{1}{2}\vector{q}_0 + \conv{T}_{\Gamma_0}[\vector{q}_0](\x) +
  \conv{N}_{\Gamma_0} [\vector{q}_0](\x) = - {\ucbar(\x)},
  \label{eqn:main-cortex-density}
\end{equation}
The operator $\conv{N}_{\Gamma_0}$ is used to complete the rank of the
operator $\conv{T}_{\Gamma_0}$ \cite{Karrila1989,Ying2006} so that
\pr{eqn:main-cortex-density} is invertible.
\Pr{eqn:main-cortex-density} defines a density for the backflow from
the periphery that offsets the complementary flow on $\Gamma_0$ due to
immersed objects.

\subsection{The contribution from rigid immersed particles
\label{ssc:nsb-flow}} Next, we consider the flow generated by the
motion of rigid immersed bodies, each moving under an externally
imposed force $\vector{F}_\pidx^\lbl{ext}$ and torque
$\vector{L}^\lbl{ext}_\pidx$ and background flow $\upbar_\pidx$.  Let
$\U^P_\pidx$ and $\vector{\Omega}^P_\pidx$ be the particle's unknown
translational and angular velocities, respectively. The external
forces and torques are generally determined by force balance amongst
the fibers and particles but here we assume they are known. The
\emph{mobility problem} for the \ordinal{{\pidx}} particle can be
written as
\begin{align}
  \up_\pidx(\x) = \conv{T}_{\Gamma_\pidx}[\vector{q}_\pidx](\x) +
  \contract{\linop{G}(\x-\X^P_\pidx)}{\vector{F}_\pidx^\lbl{ext}} +
  \contract{\linop{R}(\x-\X^P_\pidx)}{\vector{L}^\lbl{ext}_\pidx},
  \label{eqn:main-nsb-vel}
\end{align}
where $\conv{T}_{\Gamma_\pidx}$ is defined in \pr{eqn:double-layer},
$\X^P_\pidx$ is the \ordinal{{\pidx}} particle center-of-mass, and
$\linop{G}$ and $\linop{R}$ are the Stokeslet and Rotlet tensors
defined in \pr{eqn:oseen,eqn:rotlet} respectively.  Once again, taking
the limit $\x\rightarrow\x'\in\Gamma_\pidx$, and using that the
particle is moving as a rigid body, \pr{eqn:main-nsb-vel} can be
written as
\begin{align}
  \U^{P}_\pidx + \vector{\Omega}^{P}_\pidx \times (\x' - \X^P_\pidx) -
    {\upbar}_\pidx(\x') = &-\frac{1}{2} \vector{q}_\pidx(\x') +
    \conv{T}_{\Gamma_\pidx}[\vector{q}_\pidx](\x') \notag
    \\ &+ \contract{\linop{G}(\x'-\X^P_\pidx)}{
      \vector{F}_\pidx^\lbl{ext}} +
    \contract{\linop{R}(\x'-\X^P_\pidx)}{ \vector{L}_\pidx^\lbl{ext}}.
  \label{eqn:main-nsb-second-kind}
\end{align}
We impose the further constraints
\begin{align}
  &\frac{1}{|\Gamma_\pidx|} \int_{\Gamma_\pidx} \d S_y~
  \vector{q}_\pidx(\y) = \U^{P}_\pidx,\label{eqn:main-nsb-rgbt}
  \\ &\frac{1}{|\Gamma_\pidx|} \int_{\Gamma_\pidx} \d S_y~
  (\y-\X^P_\pidx) \times \vector{q}_\pidx(\y) =
  \vector{\Omega}^{P}_\pidx.\label{eqn:main-nsb-rgbr}
\end{align}
where $|\Gamma_\pidx|$ denotes the surface area of
$\Gamma_\pidx$. There are $N_P$  sets of such integral equations, one
for each immersed body.  The Stokeslet and Rotlet terms in
\pr{eqn:main-nsb-vel} are added to remove a rank deficiency of the
double-layer integral formulation and to account for the net force and
torque on the immersed particle, while the constraints,
\pr{eqn:main-nsb-rgbt,eqn:main-nsb-rgbr}, are added to make the system
fully determined \cite{Power1987,Karrila1989}.

\subsection{Fiber contributions to the flow \label{ssc:slender}}
Consider a single fiber whose centerline is given by $\X(s,t)$ where
$s \in [0,L]$ ($L$ is fiber length), moving in a background
(complementary) velocity field $\ufbar$.  In the biological setting
where the fiber is a microtubule, the length $L$ is a function of time
due to its \altr{polymerization}{depolymerization} kinetics. In that
setting, (de)polymerization typically takes place at $s=L(t)$ (at the
plus-end of the fiber); hence $s=0$ labels the minus-end, which is
stable and has no \altr{polymerization}{depolymerization}
reaction. For clarity in this section, we consider a fixed length $L$.
We assume a circular fiber cross-section with radius $a(s)$ and that
$\epsilon = a(L/2)/{L}\ll 1$.  Then slender body theory uses a matched
asymptotic procedure to relate the force per unit length, $\f(s,t)$,
that fiber exerts upon the fluid to the fiber velocity,
$\vector{V}(s,t)$, through a distribution of Stokeslets along the
fiber centerline \citep{Keller1976,Johnson1980,gotz2000}.

G\"otz~\cite{gotz2000} in particular showed that the velocity induced
by this slender fiber at a distal location $\x$ is given by
\begin{align}
  \uf(\x) &= \conv{G}[\f](\x) + \frac{\epsilon^2}{2}
  \conv{W}[\f](\x),\quad\text{where} \label{eqn:main-sb-vel}\\ \conv{W}[\f](\x)
  &= \frac{1}{8\pi\mu}\int_0^L \d s^\prime~\contract{\frac{\linop{I}
      -3 \hat{\rr} \hat{\rr}}{|\rr|^3}}{\f(s^\prime)},
\end{align}
The singular integrand in $\conv{W}$ is known as the Stokes doublet.
In our biological simulations, $\epsilon \sim \O{10^{-3}{-}10^{-2}}$,
and hence the second term in \pr{eqn:main-sb-vel} is nearly always
negligible except when in the very close proximity to other
structures. We have mostly avoided such situations in the applications
presented here, and we omit this term in the evaluation of velocity
for other particles or fibers.

To leading orders in $\epsilon$, the self-induced motion of the fiber
itself is given by
\begin{align}
  \vector{V}(s)-\ufbar(s) &= (\contract{\linop{M}}{\f})(s) +
  \conv{K}[\f](s),
  \label{eqn:main-sb-local}\\
  \contract{\linop{M}}{\f} &= \frac{1}{8\pi \mu} \contract{\left[ -\ln
      (\epsilon^2 e) \left(\linop{I}+ \X_s \X_s \right)
      +2(\linop{I}-\X_s\X_s)\right]}{\f},\\ \conv{K}[\f](s) &=
  \int_0^L \d s^\prime~\left[ \contract{\linop{G}
      \left(\X(s)-\X(s^\prime)\right)} {\f(s^\prime)} - \frac{1}{8 \pi
      \mu |s-s^\prime|}\contract{[\linop{I}+\X_s(s) \X_s(s)]}{ \f(s)}
    \right], \label{eqn:sb-stokes}
\end{align}
where $\X_s$ is the unit tangent vector to the fiber. The operator
$\conv{K}$ is a so-called finite part integral arising from the
matching procedure that makes an \O{1} contribution to the fiber
velocity.  One typical approximation is to neglect $\conv{K}$ in
comparison with $\linop{M}$, using its dominant contribution which is
proportional to $\ln(\epsilon^2)$. This is termed the \emph{local
  slender body} formulation. For completeness (and asymptotic
consistency) we keep the \emph{non-local self interaction} term
$\conv{K}$. However, we do note that in our particular studies the
$\contract{\linop{M}}{\f}$ and $\ufbar$ terms are dominant and the
nonlocal term has a negligible effect on the dynamics.

As discussed in \cite{Tornberg2004}, \pr{eqn:sb-stokes} is not
well-suited for numerical computation and requires regularization to
achieve stability and to maintain solvability. In the same manner as
\cite{Tornberg2004} we introduce a regularizing parameter $\delta$ to
$\conv{K}$:
\begin{equation}
  \conv{K}_\delta [\vector{f}](s) = \int_0^L \d s^\prime~\left[
    \frac{| \rr |}{\sqrt{|\rr|^2+\delta^2}}
    \contract{\linop{G}\left(\X(s)-\X(s^\prime)\right)}{ \f(s^\prime)}
    - \contract{\frac{\linop{I}+\vector{X}_s(s)\vector{X}_s(s)}{8 \pi
        \mu \sqrt{| s-s^\prime |^2+\delta^2}}}{ \vector{f}(s)}
    \right].
  \label{eqn:sb-stokes-reg}
\end{equation}
In the formulation of \cite{Tornberg2004} the regularization parameter
is a function of $s$ \cite[Equation 16]{Tornberg2004}, resulting in
the asymptotic accuracy of $\O{\delta^2\ln\delta}$. In our formulation
we take $\delta$ as a constant.  In \pr{ssc:numerical-tests}, we
investigate the effect different choices of $\delta$ have on the
overall accuracy and demonstrate that this formulation gives the
asymptotic accuracy of $\O{\delta}$.

We reiterate our finding that the nonlocal term has negligible effect
in our simulations composed of large number of fibers. Also, we note
that the asymptotic accuracy of $\O{\delta^2\ln\delta}$ that was
obtained in \cite{Tornberg2004} is only applicable to fibers with a
quadratic variation of radius with respect to length (\cite{Tornberg2004}
assumes  $r=4s(1-s)r(L/2)/L^2$). In cytoskeletal fibers
such as microtubules, however, the geometry is approximated
more closely by a curved cylinder with
fixed radius along its length. Thus, choosing $\delta$ in a similar
fashion to \cite{Tornberg2004} will not result in the same asymptotic
accuracy of $\O{\delta^2\ln \delta}$.
%% In other words, assuming constant $\delta$ does not reduce the
%% level of accuracy compared with \cite{Tornberg2004}.

\Pr{eqn:main-sb-local} relates the fiber velocity to the fiber forces
acting upon the fluid. Since inertial effects are negligible in the
Stokes regime, the sum of all forces at any point along the fiber is
identically zero. Thus, the hydrodynamic force applied from the fiber
to the fluid, $\f$, balances internally generated forces $\f^I$,
arising for instance from elastic deformations of the fibers, and
external forces applied to the fiber, $\f^E$, say by molecular motors
carrying payloads, or by gravitational body forces. That is,
\begin{align}
  \label{force-balance}
  \f=\f^I+\f^E.
\end{align}
The internal elastic forces are related to fiber configurations
through appropriate constitutive relations, and here we choose to use
the \emph{Euler--Bernoulli} beam theory for elastic rods.  The form of
$\f^E$ is highly dependent on the particular phenomena being modeled,
and one example concerning the positioning of the mitotic spindle
during cell division is discussed in \pr{ssc:position}.

\subsection{Mechanics of elastic fibers \label{ssc:fiber-dynamics}}
For high-aspect-ratio fibers, it is appropriate to use a generalized
form of Euler--Bernoulli beam theory for elastic rods, where the
bending moment $\vector{L}$ and bending force $\vector{F}^B$ are given
by
\begin{align}
  \vector{L} &= -E \cross{\X_{ss}}{\X_s}, \label{eqn:EB-moment}\\
  \vector{F}^{B} &= -E\X_{sss} \label{eqn:EB-force},
\end{align}
where $E$ is the flexural modulus of the fiber.  Twist elasticity is
neglected here \cite{Goldstein1998a,Lim2004}.  The local
inextensibility constraint is satisfied by the determination of a
tensile force, $T\X_{s}$, that acts along the tangent direction of the
fiber. Its magnitude $T$ is computed as a Lagrange multiplier \cite{Tornberg2004}. Consequently, the
total elastic force and elastic force per unit length applied to the
fluid are
\begin{align}
  \label{eqn:main-elastic-force}
  \vector{F}^I &= \vector{F}^{B} + T\X_s = -E\X_{sss} + T\X_s,\\
  \f^I &= \vector{F}^I_s = -E\X_{ssss} + (T\X_s)_s.
\end{align}
Imposing the local inextensibility constraint, i.e., $s$ is a material
parameter and independent of $t$, implies that $\X_t(s,t)=\vector{V}$.
Differentiating the identity $\X_s\cdot\X_s=1$ with respect to time
generates the auxiliary constraint:
\begin{equation}
  \label{eqn:main-inextensibility}
  \inner{\X_{ts}}{\X_{s}}=\inner{\vector{V}_s}{\X_{s}}=0.
\end{equation}

\subsection{Microtubule (de)polymerization kinetics
\label{ssc:polymerization}}
One major factor that allows the cytoskeleton to reprogram itself for
different functionalities is that both actin filaments and
microtubules are highly dynamic structures that continuously nucleate,
polymerize and depolymerize. It is essential to include these effects
in our biophysical simulations.  The time-scales of (de)polymerization
reactions are generally much shorter than those of cytoskeletal
rearrangements. As an example, the lifetime of microtubules in a
mitotic spindle is in the order of a minute, while the entire mitotic
spindle can be maintained stably for hours \citep{WHD2001}. One
approach to simulating (de)polymerization processes is to discretely
(remove) add segments of the filaments in time \cite{Nedelec2007}.
For the problems we aim to solve, this approach  results in severe
limitation of the maximum time-step needed for stability, and difficulty in enforcing boundary conditions.

To overcome these difficulties, we take an alternative route and
reparametrize \pr{eqn:main-sb-local} in terms of a dimensionless
parameter $\alpha(s,t)=2s/L(t)-1$ and write $\X(\alpha,t)\equiv
\X(s(\alpha),t)$.  This gives $\alpha_s=2/L$ and
$(\cdot)_s=2(\cdot)_\alpha/L$.  The chain-rule then gives
\begin{align}
  \parderiv{\X(\alpha,t)}{t} &= \vector{V} + \alpha_t\X_\alpha =
  \vector{V} - \frac{(\alpha+1) \Ldot}{L} \X_\alpha ,
  \label{eqn:chain-rule}
\end{align}
where $\Ldot$ is the rate of
polymerization. \Pr{eqn:main-inextensibility} for tension can easily
be rewritten with respect to $\alpha$ by using the chain rule
\begin{align}
  \label{eqn:inextensibility-mod}
  %% \inner{\left( \totderiv{\x_\alpha}{t }\right)}{\x_\alpha} &=
  %% \inner{\left(\x_{\alpha t} + \x_{\alpha \alpha} \alpha_t
  %%   \right)}{\x_\alpha} = \inner{\x_{\alpha t}}{\x_\alpha} = L\beta,
  %% \\ \inner{\x_{t \alpha}}{\x_\alpha} - L\beta &=
  %% \inner{\left(\totderiv{\x}{t}+\frac{\alpha \beta}{L} \x_\alpha
  %%   \right)_\alpha}{\x_\alpha} - L\beta =
  \inner{\X_{t\alpha}}{\X_\alpha} = 0,
\end{align}
which, because of linear scaling between the parameter $\alpha$ and
$s$, retains its original form.

Note that this particular choice of linear mapping between $s$ and
$\alpha$ gives
\begin{equation}
  \left.\frac{(\alpha+1) \Ldot}{L}\X_\alpha \right|_{s=0} =\vector{0}
  \quad \text{and} \quad \left.\frac{(\alpha+1)
    \Ldot}{L}\X_\alpha\right|_{s=L}
  =\frac{2\Ldot}{L}\X_\alpha=\Ldot\X_s.
\end{equation}
In other words, it is assumed that the fiber does not grow from the
minus-end at $s=0$ ($\alpha=-1$), and only grows by continuous
addition of monomers to the plus-end at $s=L$ ($\alpha=1$). The
underlying reason for this choice is the fact that microtubule are
polar filaments that only grow from their plus-end ($s=L$), while
their minus-end ($s=0$) is stable.  Nonetheless, other forms of growth
along the fiber, say from both ends, can easily be implemented by
modifying the linear relationship between $s$ and $\alpha$, such that
it reflects the known kinetics at the end-points of the fiber.

Note also that in the presence of polymerization, $\alpha$ is not a
material parameter.  Intuitively, incorporating \pr{eqn:chain-rule}
into \pr{eqn:main-sb-local} ensures that only moving a material point
with respect to the background fluid flow would result in an induced
flow and the act of adding and subtracting material elements to the
ends of the fiber through the polymerization reaction does not result
in any flow.  Note that the physics would be different if the
polymerization occurred by opening space and adding monomers in
between the two ends, which requires force and does produce a net
flow.  This condition was considered by \cite{Shelley2000} in their
simulations of closed growing filaments.

\subsection{Boundary conditions\label{ssc:boundary-cond}}
The evolution equations, \pr{eqn:main-sb-local}, are fourth-order in
$s$ for $\X$, while \pr{eqn:main-inextensibility} is second-order in
$s$ for $T$ \cite{Tornberg2004}. Generally the \bc{s} for tension are
obtained by imposing the inextensibility constraint upon the \bc{s}
for $\X$.  Here we discuss two types of \bc{s} that commonly occur in
our modeling of cellular assemblies.
\begin{enumerate}[(i)]
\item \lh{Prescribed external force $\vector{F}^\lbl{ext}$ and torque
$\vector{L}^\lbl{ext}$} on one, or both ends of the fiber. Taking
$s=L$ as an example, then using Euler--Bernoulli theory we have the
boundary conditions:
\begin{align}
  \vector{F}^\lbl{ext} &= \vector{F}^I(s=L) =
  \left. (-E\vector{X}_{sss}+T\vector{X}_{s})\right|_{s=L},
  \label{eqn:force-torque-F} \\
  \vector{L}^\lbl{ext} &= \vector{L}(s=L)
  = \left. E (\cross{\vector{X}_{ss}}{\vector{X}_s})\right|_{s=L}.
  \label{eqn:force-torque-L}
\end{align}
Similar expressions would hold at $s=0$ by changing the sign of the
right-hand-side. This provides two vector \bc{s} at $s=L$.  The
\bc\ for tension is obtained by taking the inner product of
\pr{eqn:force-torque-F} with $\X_{s}$
  \begin{equation}
    \label{eqn:BC1-T}
    T|_{s=L} =  \inner{\vector{F}^\lbl{ext}}{\X_s}|_{s=L}
    + \frac{1}{E}
    \left|\vector{L}^\lbl{ext}\right|^2,
  \end{equation}
which uses \pr{eqn:force-torque-L} after noting that
$\inner{\X_s}{\X_{sss}}=-|\X_{ss}|^2=-|\cross{\X_s}{\X_{ss}}|^2$.

When a fiber end is ``free'', that is no force or torque is applied to
it, we then have
\begin{equation}
  \label{eqn:bc-free}
  \left.\X_{ss}\right|_{s=L} =\vector{0},
  \quad \left.\X_{sss}\right|_{s=L} =\vector{0},
  \quad \left.T\right|_{s=L}=0.
\end{equation}
\item \lh{Prescribed position and velocity} when a fiber is attached
  to a rigid body. Taking $s=0$ as the attachment point we have:
  \begin{equation}
    \label{eqn:BC2-1}
    \left.\parderiv{\X}{t}\right|_{s=0} = \U^P+\vector{\Omega}^{P}
    \times (\left.\X\right|_{s=0}-\X^{P}),
  \end{equation}
  where $\U^P$ and $\vector{\Omega}^P$ are the translational and
  angular velocities of the body and $\X^{P}$ is its center of
  mass. If the fibers are \emph{clamped} at the attachment point, the
  tangent vector there would rotate with the body giving
    \begin{equation}
      \label{eqn:BC2-2}
      \parderiv{}{t}\Big(\X_s\big|_{s=0}\Big) = \cross{
        \vector{\Omega}^P}{\X_{s}}\big|_{s=0}.
    \end{equation}
    If the fiber is hinged and free to change orientation at its point
    of attachment, the torque free condition
    $\vector{X}_{ss}=\vector{0}$ is enforced at the attached end
    and the exerted torque to the particle from the fiber is set to
    zero ($\vector{L}$ in \pr{eqn:force-torque-L}).

    The \bc\ for tension is obtained by taking the inner product of
    \pr{eqn:BC2-1} with $\X_s|_{s=0}$ and using
    \pr{eqn:main-sb-local}, to find
    \begin{equation}
      \label{eqn:BC2-T}
      \inner{\left(\vector{U}^P+ \cross{
          \vector{\Omega}^P}{(\X-\vector{X}^P)} -\ufbar \right)}{\X_s}
      =-\frac{\ln(\epsilon^{2}e)}{4\pi\mu}
      \left(T_{s}-E\inner{\X_{ssss}}{\X_{s}}+\inner{\f^E}{\X_{s}}\right)
      + \inner{\conv{K}_{\delta}[\f](s)}{\X_{s}}.
    \end{equation}
    %% We then apply $\vector{X}_{s}\cdot \vector{X}_{s}=1$ to obtain an
    %% auxiliary boundary condition for tension given by
    %%     \begin{equation}
    %%   \label{eqn:BC2-T}
    %%   \left(\vector{U}^P+
    %%       \cross{ \vector{\Omega}^P}{(\X-\vector{X}^P)}
    %%       -\ufbar \right) \cdot \X_s
    %%       =-\frac{\ln(\epsilon^{2}e)}{4\pi\mu}
    %%        T_{s}-\left(\X_{ssss}+\f^E\right)\cdot\X_{s}
    %%        +\mathcal{K}_{\delta}[\f](s)\cdot\X_{s}.
    %% \end{equation}
\end{enumerate}

Finally, a fiber attached to a body applies a net force and torque to
it, so that
\begin{align}
  \vector{F}^\lbl{ext}_\wlbl{body} &= -\vector{F}|_{s=0},
  \label{eqn:f-ext}\\ \vector{L}^{\lbl{ext}}_\wlbl{body} &=
  -\Big( \left.\vector{L}\right|_{s=0}
  +\cross{(\left.\vector{X}\right|_{s=0}-\vector{X}^P)}{\left.\vector{F}
    \right|_{s=0}} \Big). \label{eqn:t-ext}
\end{align}
If there is more than one fiber attached to the body, then the
right-hand-side of these equations becomes a sum over fiber-end forces
and torques.

\subsection{Formulation summary\label{sec:form-summary}}
For ease of notation, we summarize the formulation in the context of
the biophysical problems we examine, and consider only one immersed
rigid body (i.e., $N_P=1$ with surface $P$) and many fibers all
attached at their minus-ends ($s=0$) to that body. The primary
unknowns of the system are the double-layer densities $\vector{q}_0$
and $\vector{q}_1$ on the periphery and the rigid immersed body
respectively, the translational and angular velocities of the immersed
body, $\U^{P}$ and $\vector{\Omega}^{P}$ respectively, and the
velocities $\vector{V}_\fidx$ and tensions $T_\fidx$ of the fibers
($\fidx=1, \dotsc, N_F$). Given proper constraints and boundary
conditions, coupled
\pr{eqn:main-cortex-density,eqn:main-nsb-second-kind,eqn:main-sb-local}
(repeated below) can be solved for these unknowns.  For convenience in
discussing our numerical formulation we summarize the principal
equations in the form
\begin{equation}
  \mathcal{U}^e + \mathcal{U}^s + \bar{\mathcal{U}} = \vector{0},
  \label{eqn:master-evol}
\end{equation}
where
\begin{alignat}{3}
  &\mathcal{U}^e &&= -
  \left(\begin{array}{c}
    \vector{0}\\
    \U^P + \vector{\Omega}^P \times (\x' - \X^{P})\\
    \vector{V}_\fidx=\parderiv{\X_\fidx(\alpha,t)}{t}
    +\frac{(\alpha+1)\Ldot}{L}(\X_\fidx)_{\alpha}
  \end{array}\right),
  &&\quad\text{(evolution velocities)}\label{eqn:master-ue}\\
  &
  \mathcal{U}^s &&=
  \left(\begin{array}{c}
    -\frac{1}{2}\vector{q}_0(\x)
    + \conv{T}_{\Gamma_0}[\vector{q}_0](\x)
    + \conv{N}_{\Gamma_0} [\vector{q}_0](\x)\\
    \hspace*{-2pt}-\frac{1}{2} \vector{q}_1(\x')
    + \conv{T}_{P}[\vector{q}_1](\x')
    + \contract{\linop{G}}{ \vector{F}^\lbl{ext}}
    + \contract{\linop{R}}{\vector{L}^\lbl{ext}}\hspace*{-2pt}\\
    (\contract{\linop{M}}{\f_\fidx})(\alpha) + \conv{K}_\delta[\f_\fidx](\alpha)
  \end{array}\right), &&\quad\text{(self-interaction)} \label{eqn:master-us}\\
  &
  \bar{\mathcal{U}} &&=
  \left(\begin{array}{c}
    \ucbar(\x)\\
    \upbar(\x')\\
    \ufbar_\fidx(\alpha)
  \end{array}\right). &&\quad\text{(complementary flows)} \label{eqn:master-uc}
\end{alignat}
Here $\x\in\Gamma_0$, $\x'\in P$, and $\alpha\in[-1,1]$ and the terms
for fibers are repeated for $\fidx=1, \dotsc, N_F$.  The complementary
velocities are given by
\begin{align}
  \ucbar(\x)   &= \up(\x) + \sum_{\fidx=1}^{N_F} \uf_\fidx(\x), \label{eqn:complementary-uc}\\
  \upbar(\x')  &= \uc(\x') + \sum_{\fidx=1}^{N_F} \uf_\fidx(\x'), \label{eqn:complementary-up}\\
  \ufbar_\fidx(\alpha)  &= \uc(\X_\fidx(\alpha)) + \up(\X_\fidx(\alpha))
  + \sum_{\substack{l=1\\l\ne \fidx}}^{N_F} \uf_l(\X_\fidx(\alpha)), \label{eqn:complementary-uf}
\end{align}
where
\begin{align}
  \uc(\x)   & = \conv{T}_{\Gamma_0}[\vector{q}_0](\x),
  \label{eqn:cortex-vel}\\
  \up(\x) & = \conv{T}_P[\vector{q}_1](\x)
  + \contract{\linop{G}(\x-\X^{P})}{\vector{F}^\lbl{ext}}
  + \contract{\linop{R}(\x-\X^{P})}{\vector{L}^\lbl{ext}},
  \label{eqn:nsb-vel} \\
  \uf_\fidx(\x) & = \conv{G}_\fidx[\f_\fidx](\x)
  %%+ \frac{\epsilon^2}{2} \conv{W}_l[\f_l](\x)
  .\label{eqn:sb-vel}
\end{align}

For \pr{eqn:master-evol} to be a closed set of equations, this system
requires $12 N_P$ constraints (here $N_P=1$) as well as a constitutive
law relating the configuration of fibers to their elastic force.  We
chose the Euler--Bernoulli model subject to local inextensibility
constraint as the constitutive model.  This choice in turn requires
$14N_F$ constraints for the fibers (four constraints for vector
position and two constraints for scalar tension).  $6N_P$ constraint
are furnished by \pr{eqn:main-nsb-rgbt,eqn:main-nsb-rgbr} and $6N_P$
are furnished by the force and torque balance on the particles,
namely, \pr{eqn:f-ext,eqn:t-ext}:
\begin{align}
  & \frac{1}{|P|} \int_{P} \d S_y~\vector{q}_1(\y) = \U^P, \label{eqn:nsb-rgbt} \\
  & \frac{1}{|P|} \int_{P} \d S_y~(\y-\X^P) \times \vector{q}_1(\y) = \vector{\Omega}^P, \label{eqn:nsb-rgbr}\\
  & \f_\fidx^I = -E(\X_\fidx)_{ssss} + (T_\fidx(\X_\fidx)_s)_s, \label{eqn:elastic-force}\\
  & \inner{(\X_\fidx)_{t\alpha}}{(\X_\fidx)_{\alpha}} =\vector{0},  \label{eqn:inextensibility}\\
  & \vector{F}^\lbl{ext}_\wlbl{body} = -\sum_{\fidx=1}^{N} \vector{F}^\lbl{ext}_\fidx,\label{eqn:fext}\\
  & \vector{L}^\lbl{ext}_\wlbl{body} = -\sum_{\fidx=1}^{N} \left(\vector{L}^\lbl{ext}_\fidx + \cross{(\left.\X_\fidx\right|_{s=0}-\X^P)}{ \vector{F}^\lbl{ext}_\fidx }  \right). \label{eqn:text}
\end{align}
The final $14N_F$ constraints depend on the choice of boundary
condition for the fibers and are chosen from items in
\pr{ssc:boundary-cond}.

\section{Numerical methods\label{sec:numeric}}
Below we outline the numerical evaluation of the dynamic equations
presented in \pr{sec:formulation} for evolving the conformation and
configuration of systems composed of flexible filaments, immersed
rigid-bodies, and the outer boundary.

We first discuss the spatial discretization of surfaces of the
rigid-bodies and the outer boundary and evaluating the related
integral equations on these surfaces. Then we will present the spatial
discretization of the centerline of the fibers in $\alpha$ coordinate
to evaluate the required high-order derivatives and integrals with
respect to $\alpha$.  Afterwards, we will discuss our time
discretization scheme that circumvents the numerical stiffness that
arises as a result of high-order differentiation along fibers as well
as the numerical stiffness induced by many-body hydrodynamic
interactions in crowded suspensions.  The resulting linear system of
equations for an update is solved using a preconditioned Krylov
subspace method in a fast multipole framework.

\subsection{Spatial discretization\label{ssc:spatial}}
We solve the boundary integral equations numerically using the
\nystrom\ method \cite{Kress1989}. On bodies, the integrals are
approximated on piece-wise \gl\ quadrangular surface patches.  We
represent fiber centerline positions and tensions by Chebyshev
expansions in $\alpha$, and use this expansion to evaluate integrals
or derivatives over centerlines.

\subsubsection{Non-singular integrals over surfaces\label{ssc:smooth-integral}}
Let $S$ denote a representative surface for $\x \in S$ let
$\x(\theta,\phi):$ $[0,\pi] \times [0,2\pi) \to \R^3$ be a generalized
  spherical coordinate representation of this surface.  We use a
  uniform trapezoidal grid in the ``polar'' and ``azimuthal'' angles
  to quadrangulate the surface. For more complex geometries,
  high-quality and robust algorithms \cite{Bommes2009a} can be used to
  generate the quadrangular mesh. On each such quadrangle on $S$, we
  use a tensor-product \gl\ grid to approximate surface integrals:
\begin{align}
\int_{S} \d S_{\x}~f(\x) = \sum_i \sum_j \nu_i \nu_j f(u_i, v_j) J(u_i,v_j),
  \label{eqn:patch-gl-quad}
\end{align}
where $(u_i,v_j)$ are the tensor-product \gl\ nodes in the unit
square, with the $\nu_i$'s as the corresponding weights, and $J$ is
the Jacobian of the map from the unit square to each quadrangle.  We
use a $4 \times 4$ \gl\ grid in each patch.  The Jacobian $J$ and
other geometric properties, such as the normal vector, are computed
separately in a preprocessing step and are inputs to our code.

\subsubsection{Singular integrals over surfaces\label{ssc:sing-integral}}
When the evaluation point sits on the surface $S$, the double-layer
operator given in \pr{eqn:double-layer} is singular and the integral
is defined in the principal-value sense.  Numerical evaluation of such
integrals can be done through singularity subtraction
\cite{Zinchenko2000,Keaveny2011}. For this, we use the identity
$\conv{T}_S [\vector{q}(\x)](\x) = \frac{1}{2}\vector{q}(\x)$ for $\x
\in S$ \cite{pozrikidis1992}.  Thus singular integrals in
\pr{eqn:main-cortex-density,eqn:master-us,eqn:main-nsb-second-kind}
can be rewritten as
\begin{align}
  - \frac{1}{2}\vector{q}(\x) + \conv{T}_S[\vector{q}](\x) &=
  \conv{T}_S[\vector{q}-\vector{q}(\x)](\x) = \int_S \d
  S_y~\contract{\n(\y)}{ \contract{\linop{T}(\x-\y)}{\left(
      \vector{q}(\y) - \vector{q}(\x) \right)}}
\end{align}
which is then evaluated using \pr{eqn:patch-gl-quad}, resulting in a
second-order accurate scheme \cite{Keaveny2011}.  Note that while
singularity subtraction makes the integrand bounded, derivatives of
the integrand stay unbounded and so putatively high-order quadrature
such as \pr{eqn:patch-gl-quad} may not exhibit high-order accuracy.
For more complex shapes, high order methods such as partitions of
unity \cite{Bruno2001,Ying2006} can be used.  In our numerical tests,
singularity subtraction showed satisfactory results.

\subsubsection{Fiber representation\label{ssc:numeric-fiber}}
We use a pseudo-spectral method to represent the fibers' centerlines
and to compute derivatives and integrals along them.  We denote the
centerline of a fiber by $\X(\alpha)$ where $\alpha \in [-1,1]$ and
represent $\X$ in the Chebyshev basis
\begin{equation}
  \label{eqn:cheb-basis}
  \X(\alpha) = \sum_{k=0}^p \fourier{\Xd}_k \, T_k(\alpha), \quad
  \alpha\in[-1,1],
\end{equation}
where $T_k$ denotes the \ordinal{k}-order Chebyshev polynomial of the
first kind \cite[Section~A.2]{boyd2001}, $\alpha_k = \cos(k\pi/p),
k=0,\dotsc,p$, denote the collocation points in $\alpha$ and
$\Xd=[\X(\alpha_0), \dotsc, \X(\alpha_p)]$ is the vector of
coordinates at collocation points.  The coefficients $\fourier{\Xd}$
as well as the derivative $\X_\alpha$ at the collocation points can be
computed with spectral accuracy using the FFT \cite{trefethen2000}.
The arclength of the centerline is given by
$s(\alpha)=\int_0^\alpha\left|\X_{\alpha^\prime}\right|\d\alpha^\prime$
and arclength derivative by $\X_s=\X_\alpha/s_\alpha$.  At the
beginning of a simulation, we parameterize each centerline so that the
parameter $\alpha$ coincides with the definition given in
\pr{ssc:polymerization}, i.e., $\alpha=2s/L-1$.

Differentiation of the Chebyshev expansion is performed exactly using
recursive relations for the coefficients of derivatives
\cite[Eq. A.15]{boyd2001}.  We define $\linopd{D}_\alpha$ as the
differentiation operator using Chebyshev series such that $\Xd_\alpha
= \linopd{D}_\alpha \Xd$. Similarly, we define $\linopd{D}_s$.

\subsubsection{Integration over fiber centerlines\label{ssc:integral-fiber}}
~Since the collocation points are extrema of Chebyshev polynomials,
smooth integrals are computed using Clenshaw--Curtis quadrature
weights \cite{trefethen2000}, which gives spectral accuracy.  For
evaluation of integrals in \pr{eqn:sb-vel} we use this
smooth quadrature scheme.  The integral for $\conv{K}_\delta$ given in
\pr{eqn:sb-stokes-reg} is also smooth (see \pr{ssc:numerical-tests})
and is evaluated using this quadrature method.

\subsubsection{Evaluation of nearly singular integrals\label{ssc:near}}
~When the evaluation point is close to the periphery, to a rigid
particle, or to a fiber, the integrals become \ns\ and care must be
taken for their accurate evaluation. There are robust algorithms for
high-order evaluation of \ns\ integrals in two dimensions
\cite{Helsing2008,Ojala2014,Barnett2014}, and the most interesting
recent development is \cite{Klockner2012} that has possible extension
to three dimensions.  In three dimensions, the interpolation
algorithms used by \cite{Ying2006} for boundary integrals and
similarly by \cite{Tornberg2004} for fibers are best suited for our
setting. Therefore, when evaluating near a surface, we apply the
algorithm outlined in \cite[Section 4]{Ying2006}, with the
modification that we use singularity subtraction at the nearest
boundary point to evaluate the \mbox{on-surface} integral. When
evaluating near a fiber, we use the algorithm given in \cite[Section
  3.3.3]{Tornberg2004}. The essence of both of these algorithms is the
high-order interpolation of velocity at targets points that reside in
the near region using nodes in the far region for surfaces or
fibers. The far region is defined based on the spatial grid sizes on
surfaces and fibers.

%\subsubsection{Fiber growth and catastrophe}
%The polymerization and depolymerization rates and the transition rates
%between the growing and shrinking states are given as input and are
%gleaned from the biophysical literature.  The following time-stepping
%scheme is used to update the centerline length and the location of its
%collocation points:
%\begin{align}
%  L(\impl{t}) & = L(\expl{t}) + \Ldot(\expl{t}) \deltat, \\
%  \impl{\left(\parderiv{\X(\alpha,t)}{t}\right)} & = \impl{\vector{V}}(\X) -
%  \frac{(\alpha+1) \expl{\Ldot}}{\impl{L}}
%  \expl{\X_\alpha}.
%  \label{eqn:discrete-mixed-update}
%\end{align}
%The parametrization $\alpha = 2s/L(t) - 1$ is automatically
%maintained during the simulation by this scheme.

\subsection{Time discretization\label{ssc:time}}
Due to the presence of high-order spatial derivatives in the bending
force, an explicit treatment of the evolution equation,
\pr{eqn:main-sb-local}, yields an essentially fourth-order stability
constraint on the time-step size.  To circumvent this, we use a
variation of the implicit-explicit (IMEX) method
\cite{Ascher1995a,Tornberg2004}, where \pr{eqn:main-sb-local} is
linearized and the numerically stiff terms (e.g., bending force) are
treated implicitly. The linearization is done by computing the
geometric properties of surfaces and fibers' centerlines (e.g., tangent
vector, Jacobian, collocation points, etc.) explicitly and treating
the forces and densities defined on them implicitly.
%% can be formally written as
%% \begin{equation}
%%   \X_t=\mathscr{L}(\X,T)+\mathscr{N}(\X),
%% \end{equation}
%% where the term $\mathscr{L}$ is the stiff linear operator that needs
%% to be treated implicitly and $\mathscr{N}$ is the possibly non-linear
%% terms that can be treated explicitly. The operator $\mathscr{L}$ is
%% obtained by linearization of the integral operators

Combined with the spectral spatial discretization, we find that the
implicit treatment of the bending and tensile forces removes the
time-step constraint (see \pr{tbl:bending-stiffness}).  In dilute
suspensions, the hydrodynamic interactions of the particles and the
fibers, i.e. $\bar{\u}^P_\pidx$ and $\bar{\u}^F_\fidx$ terms, can be
treated explicitly without any strict constraint on the step size.
However, as the volume fraction of fibers increase, explicit treatment
of interaction imposes strict limits on the time-step
\cite{Rahimian2015}.  Therefore, we treat the interaction of flow
constituent implicitly as well (see \pr{tbl:gmres-iter}).

We use the ``$\impltag{+}$'' superscript to mark the unknowns
to be determined at the next time-step. The unmarked variables are calculated
at the current time-step. Discretizing the evolution equations,
\pr{eqn:master-evol},
%% the constraints \crefrange{eqn:nsb-rgbt}{eqn:text}, and using the
%% disturbance velocities \pr{eqn:cortex-vel,eqn:nsb-vel,eqn:sb-vel}
using backward Euler method, we have
\begin{subequations}
  \label{set:time-stepping}
  \begin{align}
    &\mathcal{U}^e +\mathcal{U}^s+ \bar{\mathcal{U}} = \mathcal{R},
    \label{eqn:discrete-update} \\
    &\mathcal{U}^e =
    -\left(\begin{array}{c}
      \vector{0}\\
      \U^{P,+} + \vector{\Omega}^{P,+} \times (\x' - \X^P)\\
      \frac{\X^{+}_\fidx-\X_\fidx}{\deltat}
    \end{array}\right),\\
    &\mathcal{U}^s =
    \left(\begin{array}{{>{\displaystyle}c}}
      -\frac{1}{2}\impl{\vector{q}}_0(\x) + \conv{T}_{\Gamma_0}[\impl{\vector{q}}_0](\x) + \conv{N}_{\Gamma_0} [\impl{\vector{q}}_0](\x)\\
      -\frac{1}{2} \impl{\vector{q}}_1(\x^\prime)+ \conv{T}_{P}[\impl{\vector{q}}_1](\x^\prime) + \contract{\linop{G}}{ \vector{F}^{\lbl{ext},\impltag{+}}} + \contract{\linop{R}}{\vector{L}^{\lbl{ext},\impltag{+}}}\\
      \vphantom{\frac{1}{2}}(\contract{\linop{M}}{\impl{\f}_\fidx})(\alpha) + \conv{K}_\delta[\impl{\f}_\fidx](\alpha) \end{array}\right),\\
    &\bar{\mathcal{U}}=
    \left(\begin{array}{c}
      \bar{\u}^{\Gamma_{0},{+}}(\x)\\
      \bar{\u}^{P,+}(\x')\\
      \bar{\u}^{F, +}_\fidx(\alpha)
    \end{array}\right), \\
    &\mathcal{R} = \left(\begin{array}{c}
      \vector{0}\\
      \vector{0}\\
      \frac{(\alpha+1)\Ldot}{L}(\X_\fidx)_{\alpha}
    \end{array}\right) \label{eqn:discrete-rhs}.
  \end{align}
  Note that the complementary flow $\bar{\mathcal{U}}$ is treated
  implicitly. It is evaluated using
  \pr{eqn:complementary-uc,eqn:complementary-up,eqn:complementary-uf}
  given that
  \begin{align}
    &\uc[,\impltag{+}](\x)   = \conv{T}_{\Gamma_0}[\impl{\vector{q}}_0](\x),\label{eqn:cortex-vel-disc}\\
    &\up[,\impltag{+}](\x)  = \conv{T}_{P}[\impl{\vector{q}}_1](\x) + \contract{\linop{G}(\x-\X^{P})}{\vector{F}^{\lbl{ext},\impltag{+}}} + \contract{\linop{R}(\x-\X^P)}{\vector{L}^{\lbl{ext},\impltag{+}}},\label{eqn:nsb-vel-disc}\\
    &\uf[,\impltag{+}]_\fidx(\x)  = \conv{G}_\fidx[\impl{\f}_\fidx](\x).\label{eqn:sb-vel-disc}
  \end{align}
  %The right-hand-side, $\mathcal{R}$, is due to \pr{eqn:discrete-mixed-update}.
  The force along the fiber is computed as
  \begin{equation}
    \impl{\f}_\fidx = -E \linopd{D}^4_s \impl[+]{\X}_\fidx + \linopd{D}_s \left(\impl{T}_\fidx \expl(\X_\fidx)_s\right),
  \end{equation}
  where $\linopd{D}_s$ is the spatial differentiation operator defined
  in \pr{ssc:numeric-fiber}. Discretization of the constraints gives
  us
  \begin{align}
    &\frac{1}{|P|} \int_{P}\d S_y~\impl{\vector{q}}_1(\y) =  {\U}^{P,+}, \\
    &\frac{1}{|P|} \int_{P}\d S_y~({\x'}-{\X}^P) \times \impl{\vector{q}}_1(\y) = {\vector{\Omega}}^{P,+},\\
    %
    %&\linopd{D}_\alpha \left( \X^{+}_\fidx-\X_\fidx \right) \cdot \expl{\X}_{s} =\vector{0},\\
    &\inner{\left( \linopd{D}_\alpha  \X^{+}_\fidx \right)}{(\X_\fidx)_\alpha} = \inner{(\X_\fidx)_\alpha}{(\X_\fidx)_\alpha},
  \end{align}
  where the last equation is the inextensibility constraint,
  \pr{eqn:inextensibility-mod}, where we substituted the time
  derivative with its backward Euler approximation. The force and
  torque balance for the rigid particle are
  \begin{align}
    &\vector{F}^{\lbl{ext},\impltag{+}}_\wlbl{body} = -\sum_{\fidx=1}^{N} \vector{F}^{\lbl{ext},\impltag{+}}_\fidx,\\
    &\vector{L}^{\lbl{ext},\impltag{+}}_\wlbl{body} = -\sum_{\fidx=1}^{N} \left(\vector{L}^{\lbl{ext},\impltag{+}}_\fidx + \cross{({\X}-{\X}^P)}{ \vector{F}^{\lbl{ext},\impltag{+}}_\fidx }  \right) \label{eqn:summary-torque}.
  \end{align}
\end{subequations}
The boundary conditions are treated implicitly as well and linearized
if necessary.  The position of the rigid particle and the length of
each fiber are updated using $\X^{P,+} = \X^P + \deltat \U^p$ and
$L(\impl{t}) = L(\expl{t}) + \Ldot(\expl{t}) \deltat$ respectively.

\subsection{Linear solver and preconditioner\label{ssc:lin-solver}}
The system of \crefrange{eqn:discrete-update}{eqn:summary-torque} is
solved iteratively using a preconditioned GMRES method
\cite{saad2003}.  For the fibers, we use the Jacobi block
preconditioning scheme natural to the problem, where the
self-interaction blocks are formed and inverted directly using
Gaussian elimination.  For immersed particles and the periphery, we
use a Jacobi preconditioner considering only the diagonal elements of
the self-interaction blocks. To be concrete, lets consider a system
with one rigid particle and two fibers, the system of equations is
schematically
\begin{align}
  \linopd{A} =
  \def\stt{\leftarrow}%source to target
  \begin{bmatrix}
    \discrete{A}^{P\stt P}   & \discrete{A}^{P\stt F_1}   & \discrete{A}^{P\stt F_2}   & \dots  \\
    \discrete{A}^{F_1\stt P} & \discrete{A}^{F_1\stt F_1} & \discrete{A}^{F_1\stt F_2} & \dots  \\
    \discrete{A}^{F_2\stt P} & \discrete{A}^{F_2\stt F_1} & \discrete{A}^{F_2\stt F_2} & \dots  \\
    \vdots                   & \vdots                     & \vdots                     & \ddots \\
  \end{bmatrix}
  ,\quad \discrete{u}^{+} =
  \begin{bmatrix}
    \begin{pmatrix}
      \vector{U}^P\\
      \Omega^P\\
      \vector{q}_1\\
    \end{pmatrix}\\
    \noalign{\vspace{2pt}}
    \begin{pmatrix}
      \X_1\\
      T_1
    \end{pmatrix}\\
    \noalign{\vspace{2pt}}
    \begin{pmatrix}
      \X_2\\
      T_2
    \end{pmatrix}\\
    \vdots
  \end{bmatrix}^{+}
  ,\quad \discrete{b} =
  \begin{bmatrix}
    \begin{pmatrix}
      \vector{0}\\
      \vector{0}\\
      \vector{0}\\
    \end{pmatrix}\\
    \noalign{\vspace{2pt}}
    \begin{pmatrix}
      \vector{a}_1\\
      \vector{b}_1
    \end{pmatrix}\\
    \noalign{\vspace{2pt}}
    \begin{pmatrix}
      \vector{a}_2\\
      \vector{b}_2
    \end{pmatrix}\\
    \vdots
  \end{bmatrix},
  \quad \linopd{A}\discrete{u}=\discrete{b},
  \label{eqn:sche-linear-eq}
\end{align}
where, for example, the matrix block $\discrete{A}^{P\leftarrow F_1}$
denotes the discrete operator that maps the (candidate) position and
tension of the first fiber to the disturbance velocity at collocation
points of the particle.  Other matrix blocks are defined similarly and
are constructed using
\crefrange{eqn:discrete-update}{eqn:summary-torque}.
% DZ the matrix below is unclear to me
We use a preconditioner in the form
\begin{align}
  \linopd{P} =
  \def\stt{\leftarrow}%source to target
    \begin{bmatrix}
      \mathrm{diag}(\discrete{A}^{P\stt P}) & 0                          & 0                          & \dots  \\
      0                                     & \discrete{A}^{F_1\stt F_1} & 0                          & \dots  \\
      0                                     & 0                          & \discrete{A}^{F_2\stt F_2} & \dots  \\
      \vdots                                & \vdots                     & \vdots                     & \ddots \\
    \end{bmatrix},
    \label{eqn:prec}
\end{align}
where $\mathrm{diag}(\discrete{A}^{P\leftarrow P})$ denotes the
diagonal matrix constructed in turn by the diagonal entries of
$\discrete{A}^{P\leftarrow P}$. The preconditioned system is then
$\linopd{P}^{-1} \linopd{A}\discrete{u} = \linopd{P}^{-1}
\discrete{b}$.  As is demonstrated in \pr{ssc:numerical-tests}, this
preconditioning scheme combined with the implicit treatment of the
bending force and the implicit handling of hydrodynamic interactions
remove the stiffness due to both high-order derivatives and the
fiber-fiber HIs.

\subsection{Matrix-vector products and FMM \label{ssc:fmm}}
The direct evaluation of the HIs in
\pr{eqn:cortex-vel-disc,eqn:nsb-vel-disc,eqn:sb-vel-disc}, namely
through calculation of $\conv{T}_{\Gamma_\pidx}\; (\pidx=0, \ldots,
N_P)$ and $\conv{G}_\fidx\; (\fidx=0, \ldots, N_F)$, have quadratic
complexity with respect to the total number of collocation points $N$
on all surfaces and fibers. This makes simulations with a large number
of points very expensive.  The Fast Multipole Method (FMM) is used to
reduce this complexity to linear \O{N}.  In particular, we use a
kernel-independent FMM code \cite{Malhotra2014} for evaluation of
complementary velocities,
\pr{eqn:complementary-uc,eqn:complementary-up,eqn:complementary-uf}.
This FMM code uses OpenMP and MPI and is highly optimized for
distributed memory. To further speed up the computation, other pieces
of our code are parallelized. The fibers and particles are distributed
among processors. To compute the matrix-vector multiplication for the
iterative solver, given the candidate unknown vector, the forces are
initially computed locally. Afterwards, the FMM is called to compute
the complementary flow. The local interactions are computed on local
memories and added to the complementary flow and the result is
returned to the iterative solver.

For the applications presented in \pr{sec:results} we typically use
one or two computing nodes with sixteen to twenty processors per
node. A number of problems related to cytoskeletons, however, are
computationally much larger than the problems studied in this paper. For
example the mitotic spindle structure in the cell division can contain
up to one hundred thousand fibers. Simulating such structures demands
a substantial increase in computational power and the number of
processors. This is achievable within our platform and we are
currently pursuing this direction.

\subsection{Numerical tests \label{ssc:numerical-tests}}
In this section, we explore the effectiveness of our numerical
techniques.  We demonstrate the speedup obtained using the preconditioning given by
\pr{eqn:prec} in solving the linear system
\pr{eqn:sche-linear-eq} using GMRES. We study
the effectiveness of the time-stepping scheme in removing temporal
stiffness. Finally, we demonstrate the spectral accuracy in spatial
operators along fibers and explore the effect of regularization factor
in $\conv{K}_\delta$.

As a test case, we consider a sphere of radius $r$ that is enclosed by
a larger sphere of radius $6r$. $N_F$ fibers are clamped to the inner
sphere, as shown in the schematic \pr{fig:schematic-confinement}. A
fixed external force $\vector{F}$ is applied to the inner sphere which
drives the sphere-fiber assembly into motion. We evolve this system
for a short fixed time.

The first row of \pr{tbl:gmres-iter} shows the average number of GMRES
iterations per time-step\emdash/for solving \pr{eqn:sche-linear-eq},
preconditioned with \pr{eqn:prec}\emdash/given the number of fibers,
$N_F$.  Since the volume of the computational domain is fixed by the
dimensions of the outer boundary, $\Gamma_0$, the volume fraction of
the fibers increases proportionally with the number of fibers. In this
case, increasing the number of fibers reduces the average distance
between fibers and causes more pronounced HIs between them.  Because
the matrix in \pr{eqn:sche-linear-eq} becomes less diagonally dominant
(with increasing HIs) and we are not preconditioning the off-diagonal
blocks, the number of iterations shows a mild increase with the number
of fibers. The preconditioner can be improved for this case by
including off-diagonal HIs and using fast direct solvers to invert
such non-sparse matrices \cite{corona2015}.
%% Modifying the preconditioning such that it includes the off-diagonal
%% HIs without compromising the computational efficiency would be a nice
%% modification to our current implementation.

\begin{table}[!b]
  \newcolumntype{M}{>{\small\centering\hspace{0pt}\arraybackslash$}c<{$}}
  \newcolumntype{T}[1]{>{\small\centering\arraybackslash\hspace{0pt}}m{#1}}
  \centering
  \begin{tabular}{T{2.5in} M M M M}\toprule
    Number of fibers $N_F$                          & 32  & 128 & 512 & 2048 \\ \midrule
    Average number of GMRES iterations              & 3   & 5   & 9   & 16   \\
    $10^3 \times \deltat^\wlbl{max}/\tau^{E}$       & 1.0 & 1.3 & 0.5 & 0.25 \\
    Number of GMRES Iterations$/\deltat^\wlbl{max}$ & 3   & 3.8 & 18  & 64   \\ \bottomrule
  \end{tabular}
  \captionsetup{width=.8\linewidth}
    \mcaption{tbl:gmres-iter}{Average number of GMRES iterations and
    stable time-step vs. number of fibers}{The average number of GMRES
    iterations and the largest stable time-step for the mobility of a
    sphere with attached fibers in confined geometry, where
    $\tau^{E}=-{8\pi\eta L^{4}}/{\ln(\epsilon^{2}e)}$ is the
    elastic relaxation time of fibers with $L=2r$.  A sphere with
    radius $r$ is enclosed by another sphere of radius $6r$. The
    fibers are clamped to the inner sphere, where an external force
    $\vector{F}$ is exerted on the structure (see
    \pr{fig:schematic-confinement} for a schematic) and
    the velocity of the inner sphere is computed.  The length of
    the fibers is $2r$ and the force is $1/32$ of the elastic force
    scale $F^\wlbl{elastic} \propto E/L^2$ resulting in small
    fiber deformations.  The interactions between the fibers, the
    sphere, and the outer boundary are treated implicitly.  There are
    $31$ Chebyshev points per fiber.  The GMRES tolerance for relative
    residual was set to $10^{-5}$.  }
\end{table}

\begin{table}[!b]
  \newcolumntype{M}{>{\small\centering\hspace{0pt}\arraybackslash$}c<{$}}
  \newcolumntype{T}[1]{>{\small\centering\arraybackslash\hspace{0pt}}m{#1}}
  \centering
  \begin{tabular}{T{1.5in} M M M M}\toprule
    Chebyshev order                             & 11  & 21   & 31  & 41   \\ \midrule
    $10^{3}\times \delta t^\wlbl{max}/\tau^{E}$ & 0.5 & 0.45 & 0.5 & 0.45 \\ \bottomrule
  \end{tabular}

  \captionsetup{width=.8\linewidth}
  \mcaption{tbl:bending-stiffness}{Stable time-step vs. number of
    points on fibers}{The setup is identical to that of
    \pr{tbl:gmres-iter} but with fixed number of fibers $N_F=512$ and
    changing number of points per fiber.  For an explicit method,
    the bending stiffness imposes a fourth order time-step limit as
    number of points per fiber is increased.  The implicit treatment
    of the bending force removes such constraint.  }
\end{table}

We also observe a mild decrease in stable time-step as the number of
fibers (and volume fraction) is increased (second row).  Nevertheless,
the computational cost per unit of time measured as the number of
global matrix-vector applies per unit time (third row) is very
favorable in the our case compared to the mixed explicit and implicit
treatments of tension and bending forces respectively as done in
\cite{Tornberg2004}.  For example, in the case of $1024$ fibers, the
stable time-step for the scheme presented in \cite{Tornberg2004} was
found to be three orders of magnitude smaller than the implicit
formulation used here. Nevertheless we still do see a roughly linear
increase in the number of global matrix-vector products per unit time
with increasing $N_F$ (third row).

In \pr{tbl:bending-stiffness}, we investigate the effect of the
implicit treatment of high-order spatial derivatives on the stable
time-step.  As is shown in the table, increasing the number of points
on the fiber has no tangible effect on the stable time-step.  By
treating both the bending and tensile forces implicitly, we have
apparently removed the stability constraint due to having derivatives
in the computation of elastic forces.

Since the fibers positions and tensions are represented in the
Chebyshev basis, they are expected to be spectrally accurate with
respect to the number of points on the fibers.  We show this in the
context of a sphere sedimenting in free space with 32 fibers of equal
length $L=5r$ \emph{hinged} to it. \Pr{fig:spectral-sediment} shows
the error in the sphere's velocity after sedimenting one radius, as a
function of the number of points on fibers.  We used a fine grid with
$120$ points as the reference to compute the error. As expected, the
method does show spectral accuracy with respect to the spatial
resolution.
\begin{SCfigure}[.5][!bt]
  \centering
  \setlength\figurewidth{2.6in}
  \setlength\figureheight{2in}
  \includepgf{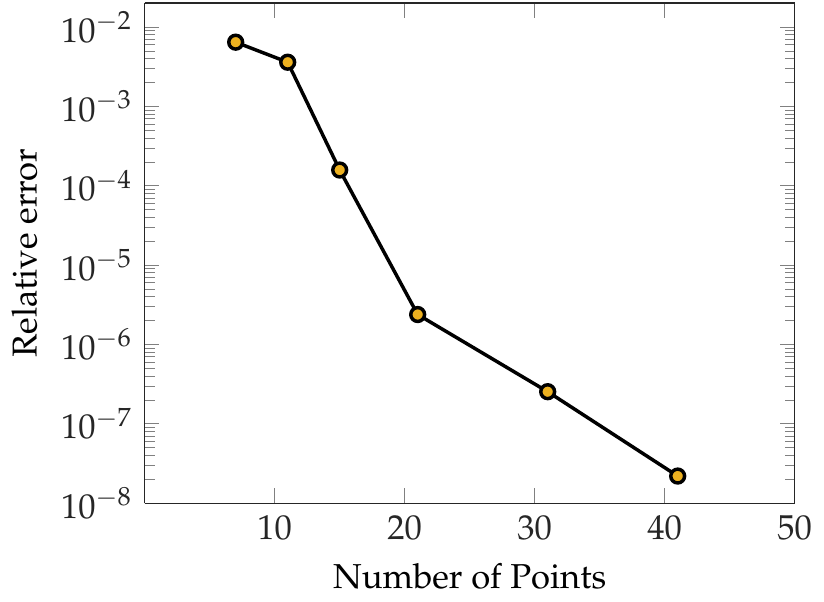}
  \mcaption{fig:spectral-sediment}{Spectral accuracy of spatial
    calculation}{A sedimenting sphere with $32$ hinged fiber is
    considered.  The abscissa is the number of points per fiber and
    the ordinate is the error in the velocity of the sphere after
    moving one radius (about one thousand time steps).  As the
    reference, we used a high-resolution case with $120$ points per
    fiber.}
\end{SCfigure}

\begin{SCfigure}[.7][!bt]
  \centering
  \setlength\figurewidth{2.6in}
  \setlength\figureheight{2in}
  \hspace*{-7pt}\includepgf{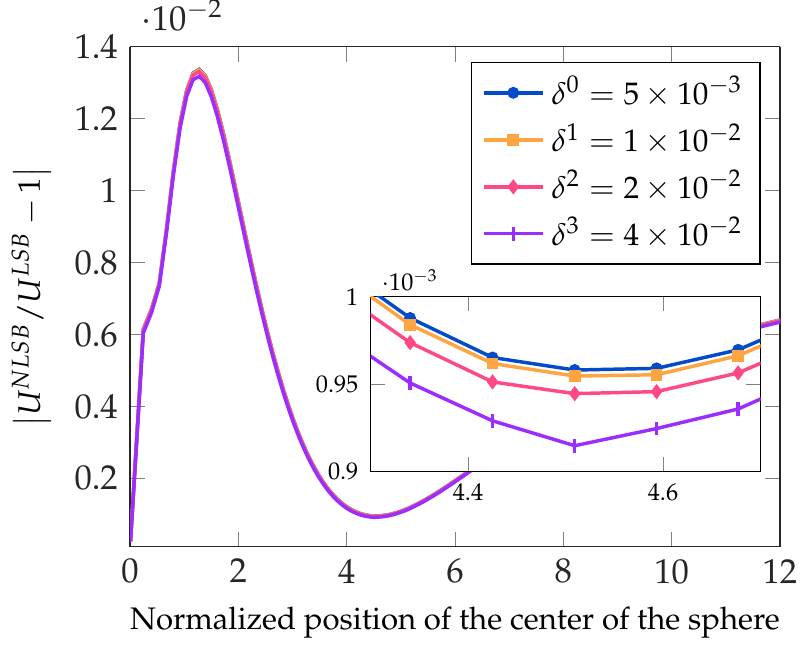}
  \mcaption{fig:regularized}{Effect of the regularization factor in
    $\conv{K}_\delta$}{The relative difference between the computed
    velocity of a sphere with $32$ attached fibers under a given
    force using local and nonlocal slender body theories at different
    values of the regularization parameter $\delta$. The maximum
    relative error associated with local slender body is
    $1.4\%$.
    The inset figure is a closeup of the same figure which shows that
    the difference between the computed velocities at different
    $\delta$  decreases with $\delta$.
    For the two smallest values of $\delta$, this difference is of
    $\O{10^{-5}}$ which is within the same range as the
    residual error for GMRES.
  }
\end{SCfigure}

Finally, we explore the effect on the overall dynamics of the
regularization factor $\delta$ that appears in the formulation of the
non-local self-interaction terms in \pr{eqn:sb-stokes-reg} .  To do
so, we use an identical simulation setup as the one used to study the
spectral accuracy of the method and compute the transient velocity
$U^{NLSB}$ of a sphere with $N_F=32$ fibers attached (hinged) to
it. For different values of $\delta$ we compare non-local slender body
theory (NLSB; Eq.~\ref{eqn:sb-stokes-reg}) against the local slender
body (LSB; dropping the term $\conv{K}_\delta$), while still taking
into account the many-body HIs.  In \pr{fig:regularized} we plot
$\left|U^{NLSB}/U^{LSB}-1\right|$ versus the travel distance of the
sphere.

As it is shown, the relative difference between the NLSB and LSB is at
most $1.4\%$, indicating that including the non-local self interaction
terms has very little effect on the dynamics of the fibers and the
attached bodies.  The inset plot of \pr{fig:regularized} shows a
closeup of the figure around $4.5r$ traveled distance.  As it can be
seen, the difference between the computed velocities at different
values of $\delta$ decreases (almost linearly) with decreasing
$\delta$, which indicates that numerically evaluated
$\mathcal{K}_\delta[\f]$ is a smooth function of $\delta$. For the two
smallest values of $\delta^{0}=5\times10^{-3}$ and
$\delta^{1}=10^{-2}$ the relative difference in the computed values
are of $\mathcal{O}(10^{-5})$. In contrast to \cite{Tornberg2004}
where specialized quadratures were used to evaluate this integral, we
find that our spectral integration is sufficient for computing
$\mathcal{K}_\delta[\f]$ accurately.

\section{Computational experiments \label{sec:results}}
In this section we consider three representative experiments.  First,
we verify the consistency of our numerical framework by studying the
effect of confinement on the mobility of a ``microtubule aster''
\citep{Nedelec2001, Weisenberg1975}. We consider an aster located in
the center of a spherical shell.  In this simple geometry, in certain
parameter ranges, the aster can be modeled as a porous medium and the
hydrodynamic drag coefficients on the body can be computed
analytically using the Brinkman model.  We show that our numerical
results are in excellent agreement with this model.  We further
demonstrate that the porous medium model has considerable shortcomings
as it fails to capture the elastic behavior of the complex when the
timescale of imposed force is shorter than the elastic relaxation time
of the fibers.

As a primary biophysical application of our framework, we study the
effect of confinement on the dynamics of ``pronuclear migration''
\citep{Siller2009}.  The precise and timely positioning of the
pronuclear complex, and the ensuing mitotic spindle, within cells is
necessary for the proper development of eukaryotic organisms
\citep{McNally2013a, Siller2009, Cowan2004}.  To gain further insight
into the mechanics of pronuclear positioning, we consider models of
positioning for cells of varying geometries. We study the time-scale
required for proper positioning and show it depends sensitively on the
choice of model and cell geometry. We also investigate the effect of
varying model parameters on the dynamics of migration.  These results
demonstrate the potential of \emph{in silico} experimentation to study
cases not easily amenable to \emph{in vitro} or \emph{ex vivo}
experiment, and to complement theoretical and experimental
understanding of cellular processes.

To show the more general applicability of our framework, we conclude
by studying a cloud of sedimenting fibers.  Our simulations reveal
that the sedimenting cloud has gross characteristics similar to those
of a sedimenting cloud of particles in either the Stokesian or
inertial regimes, but with strong internal alignment dynamics.

In all three studies, \emph{microtubule} (MT) filaments are our model
system for semi-flexible fibers. The biophysical and mechanical
parameters related to microtubules and their associated molecular
motors that are used in simulating these three conditions are listed
in \pr{tbl:sim-pars}. These values are reproduced from
\cite[Table~1]{Kimura2005}. The references related to these
measurements are also provided in \cite{Kimura2005}.

\begin{table}[b]
  \centering
  \small
  \begin{tabular}{lll}
    \toprule
    Parameter description                                  & Values used in simulations                          \\\midrule
    MT growth velocity ($V_g^0$)                           & 0.12                       $\usep\microm\usep\sec^{-1}$  \\
    MT shrinkage velocity ($V_s$)                          & 0.288                      $\usep\microm\usep\sec^{-1}$  \\
    MT rate of catastrophe ($f_\lbl{cat}^0$)               & 0.014                      $\usep\sec^{-1}$              \\
    MT rate of rescue ($f_\lbl{res}$)                      & 0.014                      $\usep\sec^{-1}$              \\
    MT bending modulus ($E$)                               & 10                         $\usep\pico\newton\meter^{2}$ \\
    MT's stall force for polymerization reaction ($F^P_S$) & 4.4                        $\usep\pico\newton$           \\
    Cytoplasmic dynein's stall force ($F^{\wlbl{stall}}$)  & 1                          $\usep\pico\newton$           \\
    Viscosity of cytoplasm ($\eta$)                        & 1                          $\usep\pascal\usep\sec$       \\
    Longest axis of the cell ( $2\times a_\lbl{cor}$)      & 30                         $\usep\microm$                \\
    Radius of pronuclear complex ($a_\lbl{PNC}$)           & 5                          $\usep\microm$                \\\bottomrule
  \end{tabular}
  \captionsetup{width=.8\linewidth}
  \mcaption{tbl:sim-pars}{The biophysical parameters used in our
    Simulation}{These values are taken from
    \cite[Table~1]{Kimura2005}. The references related to each
    measurement is also given in that article.}
\end{table}

\subsection{The effect of confinement on the hydrodynamic mobility of
microtubule asters\label{ssc:confinement}}

One of the main structural elements of the cytoskeleton is the
\emph{microtubule} (MT) filament. Some MT assemblies are formed by MTs
nucleating from \emph{microtubule organizing centers} (MTOC) and
radially growing out into the cytoplasm.  In some organisms, these
\emph{astral} structures can grow to be as large as the outer boundary
of the cell. The confining geometry of the cell may increase the force
required to move the aster and any attached structures.  To study the
effects of confinement on aster mobility, we construct a very simple
model of it by modeling the MTOC as a solid sphere where MTs of equal
length are radially \emph{clamped} to it.

\Pr{fig:schematic-confinement} shows a schematic of this setup, where
the entire structure is centered inside a spherical cell with $r_{c}=6
r_{m}$ where $r_{m}$ is the radius of the MTOC and $r_{c}$ is the
radius of the cell. We apply an external force to the aster and
compute the resulting velocity. The ratio of the external force to
this velocity is the hydrodynamic drag coefficient at this location.
\begin{SCfigure}[1.4][!bt]
  \centering
  \hspace*{-3pt}\includepgf{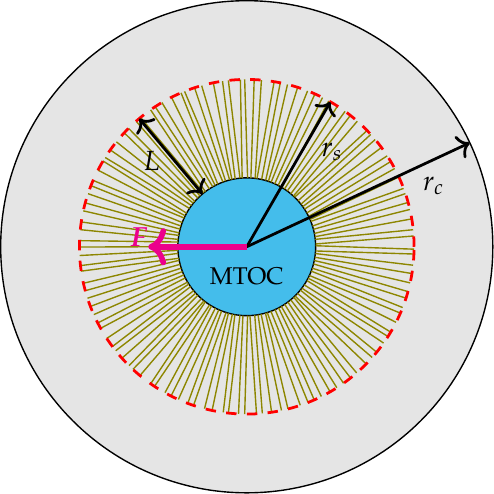}
  %
  %\captionsetup{width=.8\linewidth}
  \mcaption{fig:schematic-confinement}{Schematic of Microtubules Aster
    in Confinement}{A schematic presentation of the simulation setup
    for studying the effect of the cell boundary on the mobility of
    the microtubule asters.  $r_{m}$ denotes the radius of the MTOC
    (the blue sphere), $r_c=6r_m$ denotes the cell radius, $L$ is the
    microtubule length, $r_{s}=r_{m}+L$ is the radius of the astral
    structure, and $\vector{F}$ is an imposed force.}
\end{SCfigure}

Aside from direct numerical simulation of the complex, one could
attempt to model the astral structure as a porous medium and compute
the drag coefficient as a function of its porosity.  To verify the
physical consistency of our numerical framework, we give an analytical
calculation of the drag coefficient by representing the attached
fibers as a porous medium where the flow inside the porous domain is
modeled using Brinkman~equation \cite{Brinkman1947}:
\begin{equation}
  \mu \nabla^{2} \u-\nabla p = \frac{\mu}{\kappa}
  \left(\u-\vector{U}_\wlbl{body}\right) \mbox{~~and~~}
  \nabla\cdot\u=0,
  \label{eqn:brinkman}
\end{equation}
where $\u$ and $\vector{U}_\wlbl{body}$ are the velocities of the fluid
and the porous domain, respectively.  The term
$(\mu/\kappa)(\u-\vector{U}_\wlbl{body})$ is the frictional force
applied by the porous media on the fluid due to their relative motion,
and $\kappa$ is the permeability coefficient that is generally a
function of the orientation, aspect ratio, and volume fraction of the
fibrous region. It is important to note that the underlying assumption
of the Brinkman equation is that the entire porous domain moves as a
solid body which in our case is the velocity of the aster. This
assumption is not valid in many flow regimes, for example, in the
limit of having large enough external force, or velocity fields that
deform the MTs causing the MTs to locally have different velocity than
the MTOC. If, however, this constraint is met, previous comparisons of
simulations based on slender-body theory and boundary integral
calculations with Brinkman theory confirm that this model gives an
excellent representation of fibrous networks over a wide range of
volume fractions \citep{Higdon1996}. That said, given the geometric
structure of an aster array of MTs, a more accurate description would
(at the least) use a spatially dependent permeability
coefficient. Here for simplicity we make the assumption that we can
describe the porous shell through a single constant permeability.

A similar problem in this context has been worked out by
\cite{Masliyah1987} where they solved for the flow and the resulting
drag coefficient of a sphere with a porous shell that is being pulled
in an infinite fluid domain. Here we directly extend their results to
a confined flow with a spherical outer boundary. The velocity field in
the porous shell $r_{m} < r < r_{s}$ is modeled using the Brinkman
equation, \pr{eqn:brinkman}, while for $r_{s} < r < r_{c}$ the Stokes
equation governs the fluid motion. Flow incompressibility is applied
in both regions. The \bc{s} for velocity at $r=r_{m}$ and $r=r_{c}$
are no-slip.  On the interface of the porous and fluid, $r=r_{s}$, the
\bc{s} are continuity of fluid stress and velocity. For convenience we
rewrite the equations in terms of
$\tilde{\u}=\u-\vector{U}_\wlbl{body}$, and represent the velocity in
spherical coordinates, $\u=(u_r,u_\theta,u_\varphi)$. Since the inner
sphere is located at the center, the flow is axisymmetric, i.e.,
$u_{\varphi}=0$ and $\u(r,\theta)$ with $\theta \in [0,\pi]$. The
\bc{s} in this situation simplify to:
\begin{equation}
  \left\{
  \begin{aligned}
    r & = r_{m}:& \tilde{u}^{B}_r(r_{m},\theta)   & = 0,                                 & \tilde{u}^{B}_\theta(r_{m},\theta)   & = 0,                                                                              \\
    r & = r_{s}:& {\tilde{u}}_r^{B}(r_{s},\theta) & = {\tilde{u}}_{r}^{F}(r_{s},\theta), & {\tau}_{r\theta}^{B}(r_{s},\theta)   & =\tau_{r\theta}^{F} (r_{s},\theta), & p^{B}(r_{s},\theta) & =p^{F}(r_{s},\theta), \\
    r & = r_c  :& \tilde{u}^{F}_{r}(r_{c},\theta) & =-{U}_\wlbl{body}\cos(\theta),     & \tilde{u}^{F}_{\theta}(r_{c},\theta) & ={U}_\wlbl{body}\sin(\theta),                                                   \\
   % r & = r_c: & \tilde{u}^{F}_{r}(r_{c},\theta) & =-\inner{\vector{U}_\wlbl{body}}{\hat{\vector{r}}},       & \tilde{u}^{F}_{\theta}(r_{c},\theta) & =-\inner{\vector{U}_\wlbl{body}}{\hat{\theta}},
  \end{aligned}
  \right.
  \label{eqn:BC-brinkman}
\end{equation}
where the $B$ and $F$ superscripts refer to the Brinkman and fluid
domains, respectively, and the porous sphere moves in the $\zhat$
direction, i.e., $\U_\wlbl{body}=U_\wlbl{body}\zhat$. We solve for the
flow using the axisymmetric stream function in spherical coordinates
subject to the conditions given in \pr{eqn:BC-brinkman}. The total
force on the porous sphere can be obtained by integrating the stress
distribution over the sphere $r=r_{s}$ to find
\begin{equation}
  F=2\pi r_{s}^{2} \int_{0}^{\pi} \d\theta~\Big[ \tau_{rr} \cos \theta
    - \tau_{r\theta} \sin \theta \Big]_{r=r_{s}} \sin \theta.
\end{equation}
We can then analytically compute the drag coefficient.

As the number of MTs increases, the penetration length of the fluid
into the porous layer is reduced due to hydrodynamic screening.  In
the limit of an infinite number of MTs ($\kappa \to 0$) we expect the
drag coefficient of the structure to approach the drag coefficient of
a sphere with an effective hydrodynamic radius of $r_{s}$. The drag
coefficient of a sphere with radius $r_{s}$ centered within a sphere
of radius $r_{c}$ is given by \citep{Happel1983}:
\begin{equation}
  \label{eqn:concentric-sphere}
  \gamma_{b}= 6\pi \eta r_{s}\frac{4(1-\lambda^{5})}
        {4-9\lambda+10\lambda^3-9\lambda^5+4\lambda^6}, \quad
        \text{where} \, \lambda=\frac{r_{s}}{r_{c}}.
\end{equation}
In \pr{fig:Brenner-Brinkman} we compare the predictions of the
Brinkman model for a range of permeabilities from the highly
permeable, $\kappa=100$, to nearly impermeable, $\kappa=0.001$, as
$r_{s}$ is increased from $2$ to $5$ (recall that $r_{m}=1$ and
$r_{c}=6$).  The predicted drag coefficients from the Brinkman
equation as $\kappa \to 0$ correctly asymptotes to the limit of having
a completely impermeable sphere with $r=r_{s}$, given by
\pr{eqn:concentric-sphere}.
\begin{figure}[!b]
  \centering
  \setlength\figurewidth{3.5in}
  \setlength\figureheight{1.6in}
  \hspace*{-28pt}\includepgf{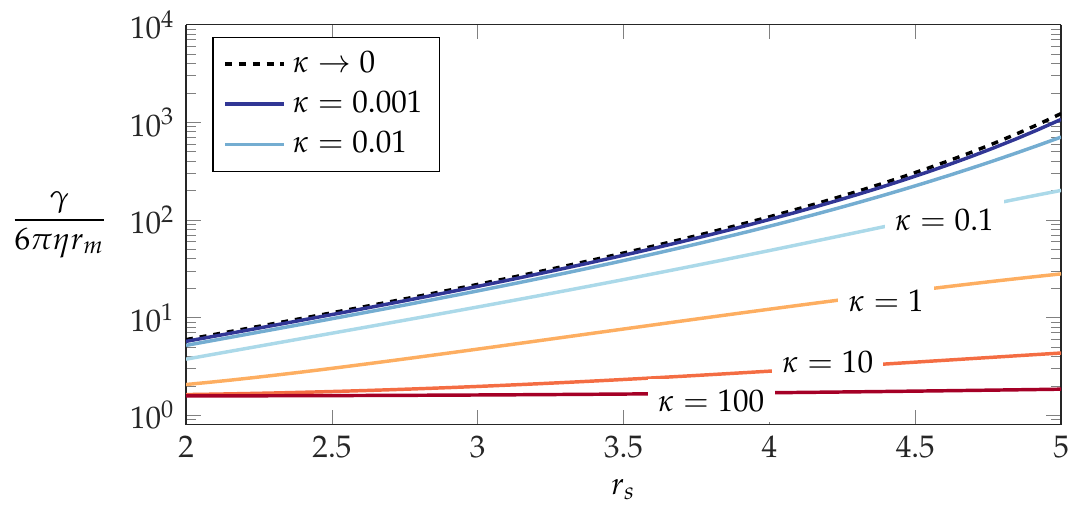}%
  \captionsetup{width=.85\linewidth}
  \mcaption{fig:Brenner-Brinkman}{Normalized drag Coefficients of a
    porous Shell}{The predicted drag coefficient based on modeling the
    porous shell with Brinkman equation compared against the drag
    coefficient of the completely impermeable case $\kappa \to 0$,
    computed using \pr{eqn:concentric-sphere}.  }
\end{figure}

Next we compare the computed drag coefficient from our direct
numerical method against the predictions of the Brinkman equation.  To
ensure that the flow regime meets the requirements of the Brinkman
model and the attached fibers and the aster move as a single body, we
chose the external force small enough, $\vector{F} = (1/32) E/L^2
\zhat$, to guarantee that the MTs remain nearly straight and that
their elastic relaxation time is much shorter than the time required
to change the position of the aster more than $1\%$ of its size.  As a
result, at any given point, the computed drag coefficient represents
the instantaneous drag coefficient (if the fiber were assumed
completely rigid) with good accuracy.  The physical parameters
including flexural modulus of fibers and fluid viscosity are chosen
from those given in Table \ref{tbl:sim-pars}, and all the lengths are
made dimensionless with $a_{PNC}=5\usep\microm$, i.e., $L=1$ is
equivalent to $5\usep\microm$.

\Pr{fig:confinement} shows the computed drag coefficient from our
simulation versus the number of the MTs for $N_F=32 \text{ to } 1568$
for three different lengths of the MTs $L=r_{m}, 2r_{m}, 3r_{m}$.  The
drag coefficient is non-dimensionalized by the value of the drag
coefficient of an impermeable sphere with radius $r=r_{s}$ from
\pr{eqn:concentric-sphere}. As expected, the drag coefficient
increases with increasing number of fibers. For example, for $r_{s}=2$
the ratio of the drag coefficients at $N_F=1568$ to $N_F=32$ is
$2.22$. This ratio is $4.3$ and $7.52$ for $r_{s}=3$ and $4$,
respectively, which shows that the effect of confinement becomes more
pronounced with the increase in fiber length, as they get closer to
the periphery.
\begin{figure}[!tb]
  \centering
  \setlength\figurewidth{3.5in}
  \setlength\figureheight{2.0in}
  \hspace*{-14pt}\includepgf{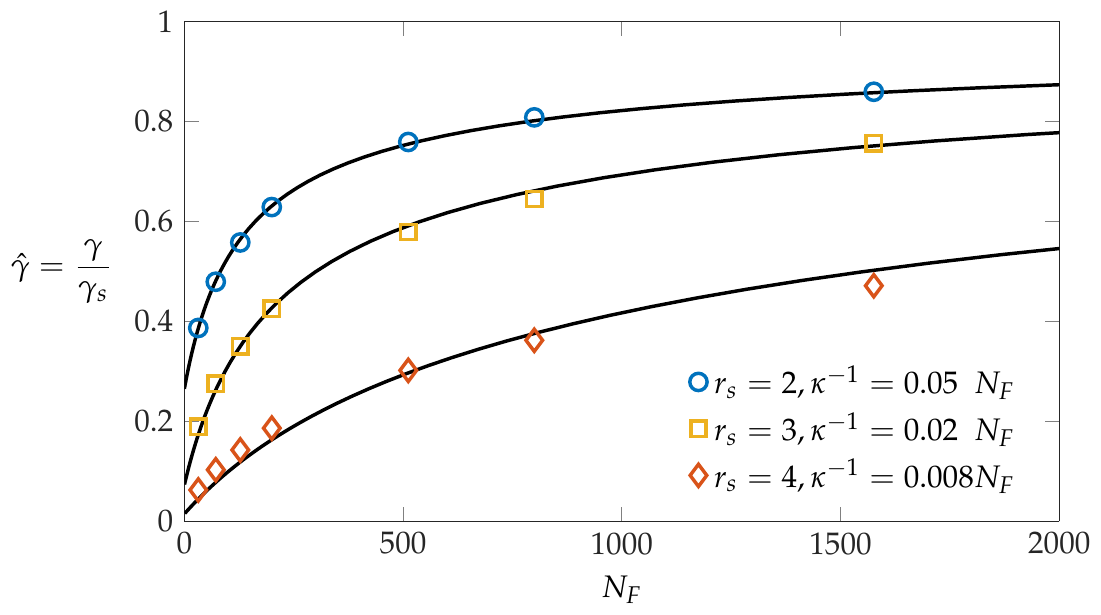}%
  \captionsetup{width=.87\linewidth}
  \mcaption{fig:confinement}{Drag Coefficient vs. Number of MTs}{The
    nondimensional drag coefficient as a function of the number of MTs
    with $r_{m}=1$, and $r_{c}=6$, and $L/r_{m}=1, 2, \text{and } 3$.
    The solid lines present the analytical calculations of the drag
    coefficient based on Brinkman equation for the porous domain.  For
    each line, the permeability coefficient used in the Brinkman
    equation is assumed to be inversely dependent on the number of
    fibers $\kappa^{-1} = c N_F$, where the coefficient of
    proportionality $c$ is computed through regression to the data.  }
\end{figure}

The other important observation is that for a given number of fibers,
the longer the length, the farther the drag is from its maximum value
corresponding to $\kappa \to 0$ in the theory and $N_F \to \infty$ in
simulations\emdash/compare $\hat{\gamma}/\gamma_b=0.86, 0.76,
\text{and }0.47$ for $L/r_{m}=1, 2, \text{and }3$ for $N_F=1568$.
This is due to the fact that the effective volume fraction of the
fibers decreases away from the MTOC as $r^{-2}$ for a given number of
attached MTs and therefore the average permeability of the porous
domain increases with $L$.

The solid curves in \pr{fig:confinement} are the Brinkman predictions
for the nondimensional drag coefficient $\hat{\gamma}$ for different
MT lengths, as $N_F$ is increased. These estimates are in very good
agreement with the numerical results (open symbols). A single
parameter is fit for each curve. Once the domain of the porous shell
is set (i.e., $r_{m,s}$ are chosen), the permeability $\kappa$ must be
specified. To remove dependence upon the number of fibers $N_F$ we
invoke the expected linear scaling of drag forces with volume fraction
and take $\kappa^{-1} = c N_F$. It is the parameter $c$ that is
determined to give the best fit of the Brinkman theory to the
simulation results.  The fit values of $c$ reduce from $0.05$ for
$L=1$, to $0.008$ for $L=3$. This shows that as the length of fibers
are increased, the average permeability in the theory increases which
is in line with our earlier observation based on numerical
simulations.

The Brinkman approach to modeling the interactions of fibrous
assemblies with the fluid has substantial limitations. We noted
earlier that we are using a constant permeability approximation.  Even
more notable is that in the Brinkman equations, the response of the
porous domain to the fluid is immediate, resulting in a purely viscous
response. However, an aster (and many other fibrous structures) are
composed of flexible fibers that deform in response to an external
force. As a result, the response of the entire system can behave as a
viscoelastic material on the time-scales relevant to many cellular
processes.

To demonstrate viscoelastic behavior of our model aster and the
limitations of Brinkman equation, we consider a setup that is
identical to the previous problem except that the external force is
now an oscillatory function of time, $F(t)=F_{0}\cos (\omega t)$. We
take $L=2r_{m}$, $F_{0}=10 E/L^{2}$ and $\omega= a\pi
\tau_{E}^{-1}$. Here $\tau_{E}=-\ln(\epsilon^{2}e) E/(8\pi L^{4}\eta)$
is the elastic relaxation time of a single fiber with length $L$ and
flexural modulus $E$. For different values of $a$ (0.32, 1.60, 8.00,
40.00, 200.00) we study the response of the aster over many periods of
oscillation.

Note that in the range of values of $F_0$ and $\omega$ used in our
simulations, the oscillation amplitude of the position of the aster is
less than $1\%$ of both its radius and the distance between the aster
and the cortex. Thus, we can safely assume that variations of the drag
coefficient with time are not due to the change in the configuration
of the aster with respect to the cell boundary.

In \pr{fig:lissajous} we plot the dimensionless velocity of the aster,
$U(t)/U_{0}$, against the dimensionless applied force, $F(t)/F_0$,
over a temporal period of oscillation (at long times) for different
frequencies $\omega$. The velocity is non-dimensionalized by the
velocity $U_{0}$ of the spherical core in the absence of MTs and under
the same force. The Brinkman predictions are also shown.
\begin{figure}[!bt]
  \centering
 \subfigure[]{\label{fig:lissajous}
   \setlength\figurewidth{1.4in}
   \setlength\figureheight{1.4in}
   \hspace{-.35in}\includepgf{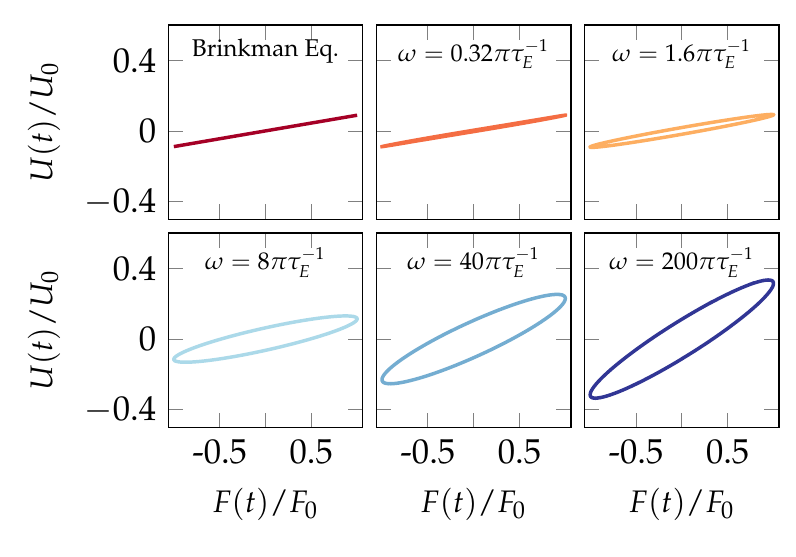}%
 }
 \subfigure[]{\label{fig:viscoelastic}
   \setlength\figurewidth{1.8in}
   \setlength\figureheight{1.62in}
   \hspace*{-6pt}\includepgf{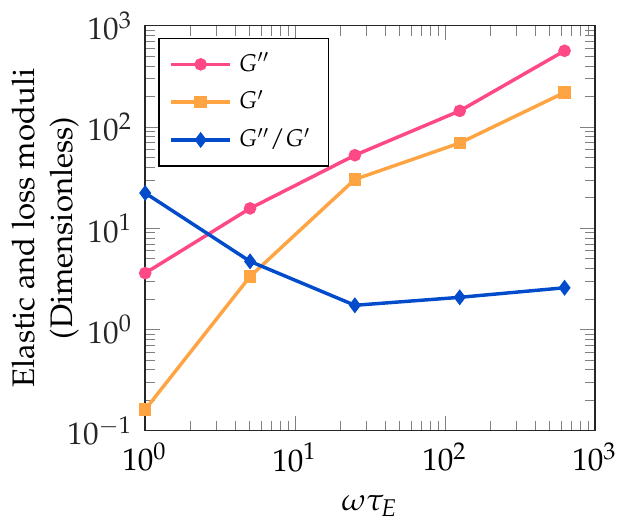}%
 }
 \mcaption{fig:oscillation-confinement}{Viscoelastic behavior of the
   microtubule asters}{\subref{fig:lissajous} The instantaneous
   velocity of the microtubule aster under an oscillatory force with
   various frequencies.  \subref{fig:viscoelastic} The loss and
   elastic moduli and their ratio versus frequency computed from the
   results shown in \subref{fig:lissajous}.  }
\end{figure}
For the Brinkman equation the response is entirely viscous, and so the
velocity is completely in phase with the applied force i.e., in the
linear response regime ${U(t)}/U_0=(\gamma_0/\gamma) \cos(\omega t)$,
where $\gamma_0$ is the drag coefficient of the spherical core in the
absence of fibers.  Thus, $U(t)/U_0$ versus $F(t)/F_0$ gives a line
with slope $\gamma_0/\gamma$. However when fibers are flexible,
their shape, and thus their resistance to the flow, evolves with time
resulting in both a frequency-dependent drag and a delay in the
response of the velocity of the aster with respect to the applied
force, i.e., $U(t)/U_0=(\gamma_0/\gamma)\cos(\omega t+\delta)$. In
this definition, $\delta=0$ and $\delta=\pi/2$ correspond to purely
viscous and elastic behaviors, respectively.  This change in behavior
is visualized by the tilted ellipses in the plots of $U(t)/U_0$ versus
$F(t)/F_0$. The area within each ellipse is proportional to the stored
elastic energy.  If $\delta=0$ the ellipse reduces to a line (zero
stored elastic energy), which is the Brinkman result. Our simulations
at the two lowest frequencies approach this limit. In this limit, the
time-scale of deformation of fibers becomes shorter than the
time-scale of oscillation, $\omega^{-1}$. As a result, the dynamics
are essentially quasi-static and the behavior is predominantly
viscous. As the frequency is increased, these two time-scales become
comparable, and the dynamics become viscoelastic and deviates
significantly from the Brinkman predictions. For example, the velocity
amplitude in simulation at the largest frequency is about $3.6$ times
larger than the value predicted by the Brinkman equation.

The dimensionless elastic and loss moduli, $G'$ and $G''$,
respectively, can be computed as $G'(\omega)=(\gamma/\gamma_0) \omega
\sin(\delta)$ and $G''(\omega)= (\gamma/\gamma_0) \omega
\cos(\delta)$. The computed values of $G'(\omega)$, $G''(\omega)$, and
their ratio are shown \pr{fig:viscoelastic}.  The minimum in $G''/G'$
gives the characteristic frequency $\omega^*$ corresponding to the
slowest relaxation time $\tau^*=1/{\omega^*}$ of the aster. For a
single fiber $\omega^*=\tau_{E}^{-1}$.  For the aster, the critical
frequency is roughly $\omega^*=8\pi \tau_E^{-1}$, that is, the aster
relaxation time is $8\pi\approx 25$ times smaller than that for an
individual fiber. Thus, in addition to changing the drag on the aster,
hydrodynamic interactions within the ensemble also give large changes
in its relaxation dynamics. To summarize, this simple example clearly
shows the viscoelastic nature of a microtubule aster, and demonstrates
the shortcomings of Brinkman equation in describing it.

\subsection{Pronuclear positioning in complex geometries
\label{ssc:position}}

The precise positioning of the pronuclear complex (PNC), and of the
ensuing mitotic spindle, is indispensable for the proper development
of eukaryotic organisms. Positioning of the spindle in the center of
the cell produces equally sized daughter cells while asymmetric
positioning leads to daughter cells of different sizes which is
essential for producing cell diversity \cite{McNally2013a}.
Positioning in organisms with MTOCs is carried out by astral MTs that
nucleate from the MTOC and grow into the cytoplasm. Previous studies
on such organisms have shown that the interaction of MTs with the cell
cortex (the periphery) can be a key factor in defining the position
and orientation of the mitotic spindle. For example, experiments in
early stages of cell division in the \celegans\ embryo have
shown that when the cell cortex is deformed from its native elliptical
shape to spherical, the mitotic spindle can fail to properly align
with the anterior-posterior AP-axis of the cell
\cite{Park2008}. Considering that cells evolve through a variety of
shapes during development, it is important to identify the
relationship between cell shape and spindle positioning.

To study this, \cite{Minc2011} molded individual sea urchin eggs into
microfabricated chambers of different geometries and so took the shape
of the chamber. They analyzed the location and orientation of the
nucleus and mitotic spindle throughout the first cell division and
found that the nucleus moves to the center of the chamber and aligns
along its longest direction. We demonstrate the applicability of our
method in studying the effect of confinement on spindle positioning by
taking a similar approach and numerically studying the positioning of
a model PNC for shapes other than spheres and ellipsoids.

\subsubsection{Stages and mechanisms of pronuclear positioning}
We concentrate on positioning of the pronuclear complex in the
single-cell \emph{\celegans\ embryo}. After fertilization and the
introduction of the male pronucleus, the female pronucleus approaches
and fuses with its male counterpart, which has two MTOCs and
associated astral MT arrays. Together this ensemble forms the PNC and
its motions are associated with its astral MTs. After centering and
rotation of the PNC leading to the alignment of the axis between the
MTOCs with the AP-axis (the two primary aspects of positioning), the
mitotic spindle forms, chromosomes condense, chromatid pairs are
formed, chromosomes are divided and pulled to opposite cell sides, and
cell division proceeds \citep{Cowan2004, Siller2009,
  McNally2013a}. Here we only study the positioning of the pronuclear
complex prior to the formation of spindle. A schematic representation
of the important structural elements involved in PNC positioning, and
its important stages, are shown in \pr{fig:sche}.
\begin{figure}[!tb]
  \begin{minipage}[c]{.45\columnwidth}
    \subfigure[Key structural elements involved in the PNC migration\label{sfg:schematic-pnc}]{
      \def\figscale{.45}
      \def\figannotate{true}
      \def\pnclocx{1.5}
      \def\pnclocy{-.2}
      \def\pncangle{90}
      \hspace*{-7pt}\includepgf{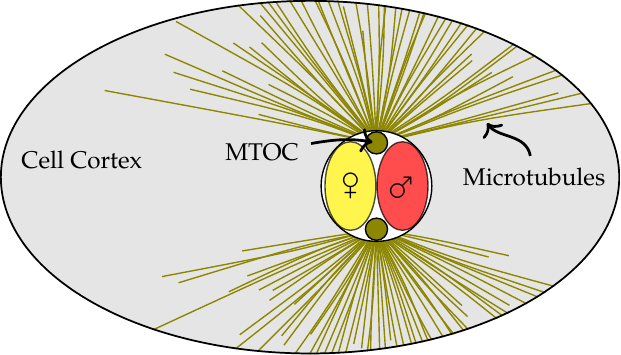}}
  \end{minipage}
  \begin{minipage}[c]{.45\columnwidth}
    \def\figannotate{false}
    \subfigure[$t_{(-1)}$, egg fertilization\label{sfg:schematic-pnc-fert}]{
      \def\figscale{.25}
      \def\pnclocx{3.5}
      \def\pnclocy{0}
      \def\pncangle{90}
      \hspace*{-5pt}
      \includepgf{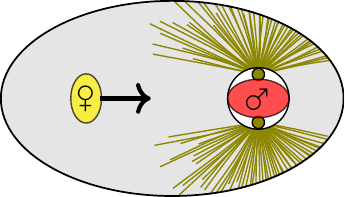}}
    \subfigure[$t_{0}$, PNC formation\label{sfg:schematic-pnc-form}]{
      \def\figscale{.25}
      \def\pnclocx{3}
      \def\pnclocy{-.75}
      \def\pncangle{90}
      \hspace*{-5pt}
      \includepgf{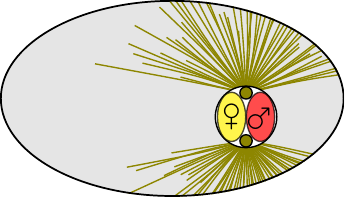}}\\
    \subfigure[$t_{1}$, PNC migration\label{sfg:schematic-pnc-mig}]{
      \def\figscale{.25}
      \def\pnclocx{1.5}
      \def\pnclocy{-.4}
      \def\pncangle{130}
      \hspace*{-5pt}\includepgf{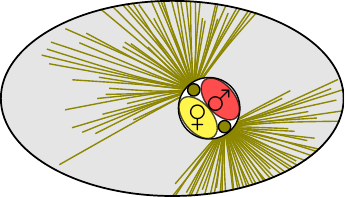}}
    \subfigure[$t_{2}$, PNC centering\label{sfg:schematic-pnc-cen}]{
      \def\figscale{.25}
      \def\pnclocx{0.2}
      \def\pnclocy{-.2}
      \def\pncangle{175}
      \hspace*{-5pt}
      \includepgf{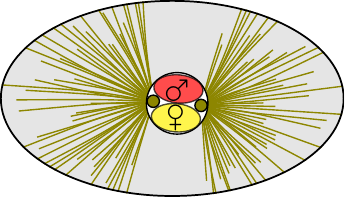}}
  \end{minipage}
  \mcaption{fig:sche}{Schematics of the Pronucleus Complex and stages
    of PNC migration}{
    %\subref{sfg:schematic-pnc} Important structural elements involved in the PNC migration.
    \subref{sfg:schematic-pnc-fert}--\subref{sfg:schematic-pnc-cen}
    The important stages of the migration process in the first cell
    division of \celegans\ embryo.
    \subref{sfg:schematic-pnc-fert} At $t_{(-1)}$ the egg is
    fertilized by the male pronucleus and the female pronucleus
    migrates towards the male pronucleus.
    \subref{sfg:schematic-pnc-form} At $t_0$ the PNC is formed in the
    posterior side of the cell.  This stage is the starting point of
    our simulations.  \subref{sfg:schematic-pnc-mig} Afterwards, from
    $t_{1}$ to $t_{2}$ the PNC migrates to the center of the cell,
    while the PNC rotates and the MTOC axis starts to align with the
    AP-axis.  \subref{sfg:schematic-pnc-cen} Finally, at $t_{3}$ the
    centering and rotation process of the PNC is completed
    corresponding to the final stage of our simulations.}
\end{figure}

The force driving PNC positioning is thought to be generated by one,
or all, of three potential force-transduction mechanisms operating on
astral MTs: \emph{cortical pushing, cytoplasmic pulling, or cortical
  pulling}. In the \emph{cortical pushing} mechanism, pushing forces
are applied by MTs growing against the cell cortex \cite{Reinsch1998}.
In the \emph{cytoplasmic pulling} mechanism, MTs are pulled by
molecular motors located in the cytoplasm \cite{Kimura2005,
  Kimura2011}. In both of these proposed mechanisms, the rotation of
PNC and its alignment with the AP-axis is achieved by asymmetric
ellipsoidal shape of \celegans\ embryo \citep{Shinar2011,
  NRNS2015}. Finally, in the \emph{cortical pulling} mechanism, MTs
are pulled by dyein motors attached to the cortex
\cite{Grill2001}. The activation of these motors is believed to be
regulated by other protein complexes that are distributed
asymmetrically throughout the cell boundary, and that because of this
asymmetry the details of cell shape may not be crucial to achieving
proper positioning \cite{Park2008}.

In a concurrent work \cite{NRNS2015} using the framework presented
here, we have studied the positioning of the PNC in the single-cell
\celegans\ embryo while focusing on the effect of hydrodynamic
interaction and the generated cytoplasmic flows in these three
mechanisms. For more details the reader is referred to
\cite{NRNS2015}. Here, we instead consider PNC positioning under
deformations of the cell shape.

\subsubsection{Biophysical models}
We consider simple instantiations of the \emph{cortical pushing} and
\emph{cytoplasmic pulling} models. Both of these models are generic in
the sense that they rely upon rather nonspecific elements of cellular
physiology. The \emph{cortical pulling} model, on the other hand,
involves several biophysical elements that are specific to
\celegans\ \cite{McNally2013a}. Thus, to keep our study as
general as possible with respect to the choice of the organism, we do
not consider the \emph{cortical pulling} here. The general features of
that model are discussed in \cite{NRNS2015}. Below we briefly outline
\emph{cortical pushing} and \emph{cytoplasmic pulling} models and
their implementation within our numerical method. First we begin with
a discussion of MT dynamics.

\para{Microtubule polymerization kinetics} Microtubules are polar
protein polymers that primarily grow and shrink from their so-called
\emph{plus-end} while the \emph{minus-end} remains stable.  The
process of abrupt stochastic transitions between growth and shrinkage
is termed dynamic instability \cite{Desai1997}. The distribution of MT
lengths is determined by their rates of growth $V_{g}$ and shrinkage
$V_s$, and their frequencies of \emph{catastrophe} $f_\lbl{cat}$
(changing from growing to shrinking), and of \emph{rescue} $f_{r}$
(changing from shrinking to growing).  Previous \emph{in vitro}
measurements and theoretical studies show these rates change under
mechanical load \citep{Peskin1993, Dogterom1997, Janson2003}. We use
an empirical relationship based on the \emph{in vitro} measurements of
\cite{Dogterom1997} to relate the rate of MT growth to an applied
compressional load on its plus-end:
\begin{equation}
V_{g}=V_{g}^0\exp \left(-\frac{7}{3}
  \frac{\inner{\vector{F}(L)}{\X_{s}}}{F_{S}^{P}} \right),
\end{equation}
where $V_{g}^0$ is the growth rate under no compressive load, $F_S^P$
is the stall force for MT's polymerization reaction and
$\vector{F}(L)$ is the end-force of the MT. The values used in our
simulation for these parameters are listed in \pr{tbl:sim-pars}.

The \emph{in vitro} measurements of \cite{Janson2003} suggest that the
turnover time from growth to shrinkage,
$\tau_\lbl{cat}=1/f_\lbl{cat}$, of MTs under compressive
plus-end-loading is proportional to its growth velocity (itself
modulated by the end-force), while other \emph{in vivo} observations
suggest that the turnover time of MTs touching cortex in \celegans\ is
$1$ to $2$~seconds \cite{McNally2013a}. We incorporate these two
observations and model the rate of catastrophe as
\begin{equation}
  f_\lbl{cat}=\max \left( f^{0}_\lbl{cat} \frac{V_{g}^0}{V_{g}},
  \frac{1}{\tau_\lbl{cat}} \right),
  \label{eq:tau_cat}
\end{equation}
where $f^{0}_\lbl{cat}$ is the rate of catastrophe under no
compressive end-load and $\tau_\lbl{cat}$ is the turnover time. In all
simulations presented here, unless specified otherwise, we chose
$\tau_\lbl{cat}=2\usep\sec$.

\para{The cortical pushing model} The movement of the PNC in the
cortical pushing model is achieved by compressive forces being applied
from the cortex to MTs growing against it. See the schematic of this
model in \pr{sfg:sche-CorPush}. If these forces are not strong enough
to stop the growth process altogether, then newly polymerized MT is
pushed out from the wall at the polymerization rate, thus lengthening
the MT. In turn, this lengthening either pushes the attached PNC away
from the cortex, or it deforms the MT.  A rough, but useful, estimate
of the likelihood of the deformation of MTs upon reaching the cell
boundary can be achieved by comparing the stall force for the growth
reaction, $F^{P}_{S}=4.4\usep\pico\newton$ \cite{VanDoorn2000}, to the
characteristic force required for buckling, $F_b=\pi^{2} E/L^{2}$. In
the single cell \celegans\ embryo, the average MT length is roughly
$10\usep\microm$, resulting in $F_b/F^{S}_{P} \approx 0.20$ implying that
MT interactions with the cortex most likely result in buckling (or
bending) rather than the stall of polymerization reaction. As a
result, a larger compressive force is applied to the MTOC with shorter
average length of anchored MTs, since they both buckle less easily and
more directly apply force to the PNC, and are more numerous due to
dynamic instability. This results in the PNC being pushed away from
the side with shorter MTs, and migration towards the side with longer
MTs. The combination of torque and force balance on the PNC eventually
determines its position and alignment of the centrosomes with respect
to the AP-axis \citep{NRNS2015}.
\begin{figure}[!bt]
  \centering
  \def\figscale{.45}
  \subfigure[Cortical Pushing\label{sfg:sche-CorPush}]{
    \fbox{\hspace*{-7pt}\includepgf{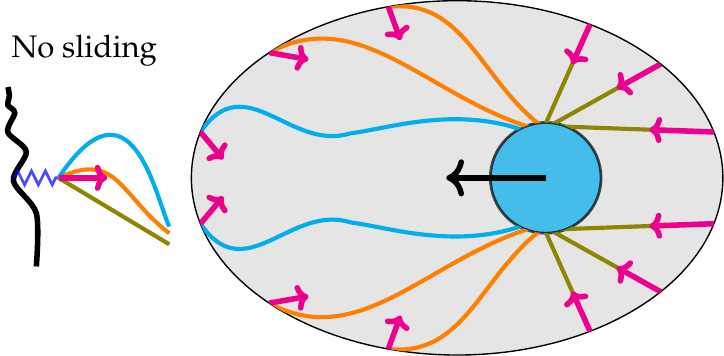}}
  }
  \subfigure[Cytoplasmic Pulling\label{sfg:sche-CytoPull}]{
    \fbox{\hspace*{-7pt}\includepgf{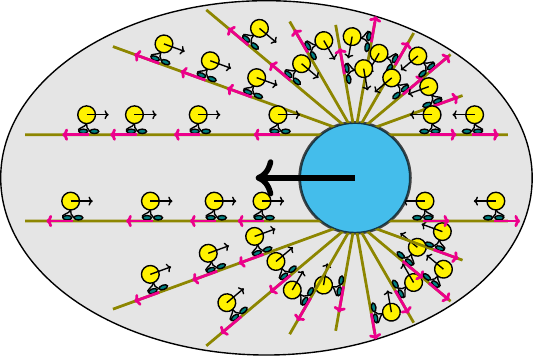}}
  }
  \mcaption{sche-models}{Schematics of the biophysical models of PNC
    migration used in this study}{ \subref{sfg:sche-CorPush} Cortical
    pushing model, where compressive forces are applied from the
    cortex to the MTs growing against it due the polymerization (the
    forces are modeled as spring forces at the attachment point)
    \subref{sfg:sche-CytoPull} cytoplasmic pulling model, where the
    active force for movement is generated by cargo carrying
    cytoplasmic dynein motors.}
\end{figure}

In our simulations, we assume that the positions of the growing MTs
remain fixed on the periphery once they make contact. We capture this
assumption by applying a constraining spring force once an MT plus-end
reaches within the critical distance $\Delta r=0.5\usep\microm$ from the
periphery. The spring force,
\begin{equation}
   \label{eqn:cortical-pushing}
   \vector{F}(L)=-K(\X(L) -\X_\lbl{att}),
\end{equation}
is directed towards the attachment point $\X_\lbl{att}$ and proportional
to the distance from it. We set $K=10\usep\pico\newton(\microm)^{-1}$.  At
$0.5\usep\microm$ of displacement, this choice of $K$ results in
$5\usep\pico\newton$ of force, which is bigger than the stall force for
polymerization reaction ($F^{P}_{S}=4.4\usep\pico\newton$). Thus, the MTs
will stop growing prior to reaching to the cell periphery.

We also assume that the attachment to the cortex is a \emph{hinged}
attachment, i.e. the net torque on the growing plus-end is zero. As a
result, for growing MTs pushing on the cortex the force
$\vector{F}^\lbl{ext}$ is given by \pr{eqn:cortical-pushing} and the
net torque is zero. When the plus-ends are not at the cortex then the
external force and torque on the plus-end are both zero. The minus-end
of all MTs are \emph{clamped} to the MTOC and the \bc{s} are
prescribed by \pr{eqn:BC2-1,eqn:BC2-2,eqn:BC2-T}. The MTs apply a net
force and torque to the PNC through their boundary conditions which is
computed using \pr{eqn:fext,eqn:text}. Note that in this model there
are no active forces from molecular motors on the MTs i.e. $\f^E=0$.

In \cite{NRNS2015} we use another variation of cortical pushing
mechanism where the growing MTs can bend, grow, and slide along the
outer boundary. We found that this variation, unlike the one
considered here, does not properly align the PNC centrosomes along the
AP-axis in physiologically reasonable times, and so do not discuss
that variant here.

\para{Cytoplasmic pulling} Cytoplasmic dynein is a minus-end directed
molecular motor that attaches and walks along MTs to carry cargo.
Through this action, dyneins apply pulling forces on the MTs that are
equal in magnitude and opposite in direction to the force they need to
exert to drag the cargo through the cytoplasm. Assuming that dynein
motors are uniformly distributed in the cytoplasmic volume, the number
of motors on MTs increases linearly with their length. As a result,
the PNC is pulled in the direction of its longest centrosomal MTs. A
schematic of this mechanism is shown in \pr{sfg:sche-CytoPull}. This
pulling mechanism can result in the proper centering and alignment of
the PNC in computational models of \celegans\ embryo \cite{Kimura2005,
  Kimura2011, Shinar2011, NRNS2015} that have various degrees of
biophysical verisimilitude. In a simple instantiation of this
mechanism, we treat the density of the attached dyneins as a continuum
field with constant number of attachments per unit length of MT. Thus,
the force per unit length, $\f^E_i(s)$, and total force,
$\vector{F}^E_i$, applied by cargo-carrying dyneins on the \ordinal{i}
MT is given by
\begin{equation}
  \f^E_i(s)=F^{\lbl{dyn}} n_{\lbl{dyn}} \X_{i,s}(s),
  \quad\text{and}\quad \vector{F}^E_i=\int_0^L \d
  s~f^E_i(s)=F^{\lbl{dyn}} n_{\lbl{dyn}} \left[\X_i(L)-\X_i(0)\right],
\end{equation}
where $F^\lbl{dyn}$ is the magnitude of the force applied from a
single motor to an MT, $n_\lbl{dyn}$ is the number density of the
dyneins per unit length, and $\vector{X}_{i,s}$ is the tangent vector
to \ordinal{i} MT, where $F^\lbl{dyn}$ is related to the walking speed
of the motor by a force-velocity relationship as $F = F ^\wlbl{stall}
(1 - \max(|V |, V_\lbl{max} )/V_\lbl{max} )$, where $F^\wlbl{stall}$
is the maximum force applied by the dynein motor on an MT.  Moreover,
the force and velocity are linearly related through the drag
coefficient of the cargo: $F^\lbl{dyn}=\gamma V$.  If we take the
average radius of the cargos as $r_\wlbl{cargo}=0.1\usep\microm$, and
the cytoplasmic viscosity as $\eta = 1\usep\pascal\usep\sec$ (see
\pr{tbl:sim-pars}), we can compute $F^\lbl{dyn}$ by combining these
two force-velocity relationships which gives
$F^\lbl{dyn}=0.83F^\wlbl{stall}=0.91\usep\pico\newton$
\citep{NRNS2015}. We used $F^\lbl{dyn}=0.91\usep\pico\newton$ in all
the simulations of the cytoplasmic pulling mechanism, presented in
this paper.

The \bc{s} on MTs are clamped \bc{s} for the minus-ends, prescribed by
\pr{eqn:BC2-1,eqn:BC2-T}, and zero force and torque at the plus-end,
given by \pr{eqn:bc-free}. Also when the MTs reach the periphery
(practically, within the small distance $\Delta r$) they
instantaneously go through catastrophe. Hence, no pushing forces are
applied on the cortex by MTs.

Finally, we note that in both the cytoplasmic and the cortical pushing
mechanisms the rotation of PNC and the alignment of the axis of MTOC
with the AP-axis is achieved by a symmetry-breaking torque instability
that arises from the coupling between the asymmetry of the shape of
the cell periphery and the length dependency of the active forces
($\vector{F}^{E}_i=F^\lbl{dyn} \left[\X_i(L)-\X_i(0)\right]$ and
$F_b=\pi^2 E/L^2$) \citep{NRNS2015}. For a spherically shaped
eggshell, both models predict centering of the PNC; however MTOC axis
will not align with the AP-axis and due to the stochastic nature of
the dynamic instability of the microtubules' polymerization dynamics
all the alignment directions are equally sampled over long times.

\subsubsection{Other assumptions and simulation setup}
Below we outline other assumptions made in simulating pronuclear
migration using both cytoplasmic pulling and cortical pushing models.

Previous studies have shown the nuclear envelope encompassing the
nucleus is much stiffer than the plasma membranes of the cell
\cite{Dahl2004}. Based on this and no observations of substantial
deformation of the PNC prior to mitosis, we model the PNC as a rigid
sphere of radius $a_\lbl{PNC}=5\usep\microm$ \cite{Kimura2005}. The
mechanics of the cell cortex, however, are in principle more involved.
Many important cellular processes, including cytokinesis, cell
crawling, and early motion of the female pronucleus towards the male,
are achieved by elaborate spatiotemporal deformations of the cell
cortex \cite{Gilden2010}. Nevertheless, the cortex is not known to
undergo large distortions during the PNC migration stage. For
simplicity then we model the cell cortex as a rigid, fixed surface.

We model the two MTOCs attached to the PNC as two regions set at
opposite poles on the spherical PNC surface ($\theta=0,\,\theta=\pi$)
defined by the polar angle as $\theta_1\le \pi/5$ and $\theta_1\ge
4\pi/5$. We assume that microtubules are mechanically \emph{clamped}
to the PNC and the anchoring sites are uniformly distributed on the
MTOC surface area.  To obtain a smooth velocity field induced by
hydrodynamic interaction of the PNC and the MTs near the anchoring
points, we place the anchoring sites slightly away from the surface of
the PNC at $1.05 a_\lbl{PNC}$. Our numerical experiments show that
changing this radius in the range of $1.02$ to $1.15$ gives changes in
the dynamics that are smaller than the fluctuations in the dynamics
due to the stochastic dynamic instability. The same procedure was used
for the interaction between MTs and the outer boundary, i.e. the MTs
cannot approach closer than $0.95 a_\lbl{cor}$ to the periphery.
Again, changing this distance slightly yielded insignificant changes
in the dynamics.

We initialize the simulation by assuming that all MTs start with the
same length equal to one-half of PNC radius. The PNC is initially
positioned towards the posterior side and, in most simulations, on the
AP-axis (denoted here by $\zhat$). We also performed a number of
simulations where the PNC started $0.25 a_\lbl{PNC}$ away from the
AP-axis, and we did not observe any qualitative change in the dynamics
or its time-scales. Finally, the axis of MTOCs is set initially at a
$90$ degree angle to the AP-axis. This setup approximately replicates
the \emph{in vivo} observations and is consistent with previous
modeling efforts \cite{Kimura2005,Shinar2011}. In all simulations we
consider $N_F=300$ fibers, and discretize each with $N=31$ points
except for one set of simulations where $N=51$ was chosen to resolve
the large fiber deformations arising in cortical pushing simulations
(see \pr{sfg:S1_tauvar_2}).

\subsubsection{Positioning of the {PNC} in various cell geometries}
We represent the shape of the cell periphery in generalized spherical
coordinates for an axisymmetric body:
\begin{align}
  x_\lbl{cor}=R(\theta)\cos(\varphi) \sin (\theta),\quad
  y_\lbl{cor}=R(\theta) \sin(\varphi) \sin (\theta), \text{ and }\quad
  z_\lbl{cor}=R(\theta)\cos (\theta),
\end{align}
where $\varphi\in[0,2\pi)$, and $\theta\in[0,\pi]$.  We consider three
  different cell shapes, $S_{1,2,3}$, defined by
\begin{equation}
  \left\{
  \begin{aligned}
      & S_1:\quad R_\lbl{cor}^{1}(\theta)=a_\lbl{cor}\left[1-0.35\cos^{3}\theta-0.15\sin^{3}\theta\right], \\
      & S_2:\quad R_\lbl{cor}^{2}(\theta)=a_\lbl{cor}\left[1+0.2\cos(2\theta)-0.2\sin(2\theta) \right], \\
      & S_3:\quad R_\lbl{cor}^{3}(\theta)=a_\lbl{cor}\left[1-0.6\sin^{4}\theta\right].
  \end{aligned}
  \right.
  \label{eqn:S-shapes}
\end{equation}
Here $a_\lbl{cor}=25\usep\microm$ is chosen in accordance with the
longest axis of \celegans\ embryo (see \pr{tbl:sim-pars}).  The
corresponding shapes are shown in \pr{fig:S-shapes}. In all cases, the
$\zhat$ corresponds to the longest dimension of the cell and its axis
of symmetry. In shapes $S_1$ and $S_2$ the symmetry is only imposed
about the $\hat{\vector{z}}$ axis, while shape $S_3$ has a mirror
symmetry about the $z=0$ plane that cuts through its middle
(corresponding to the cell division plane). Also the extent of
asymmetry is increased from $S_1$ to $S_2$. Thus we can simultaneously
observe the effect of the type and extent of symmetry on the migration
dynamics of the PNC in shapes $S_1$, $S_2$, and $S_3$.  These choices
of cell confinements in comparison with the typical ellipsoidal shape
of the embryos also allow us to study how local deformation of the
internal geometry of the cell\emdash/by, say, the local contraction of
the actin network of the cortex \citep{Cowan2004}\emdash/can affect
the dynamics of migration.

A major effect of changing the cell shape is changing the length
distribution of MTs within the cell. Since both cortical pushing
forces, $\vector{F}(L)$, and cytoplasmic pulling active forces,
$\vector{F}^E$, depend on the length of the MTs\emdash/$\vector{F}
(L)\propto L^{-2}$ and $\vector{F}^E \propto L$\emdash/and average MT
length is not well-separated from typical cell sizes, we expect that
the variations in the cell shape result in variations in the PNC
migration time. Similarly, since asymmetric shape of the cell is a
necessary factor in PNC rotation and achieving alignment with the AP-axis
in both models, rotation time is also expected to change with the
geometry.

\begin{figure}[!bt]
  \centering
  \subfigure[$S_1$]{\includegraphics[height=1.80in]{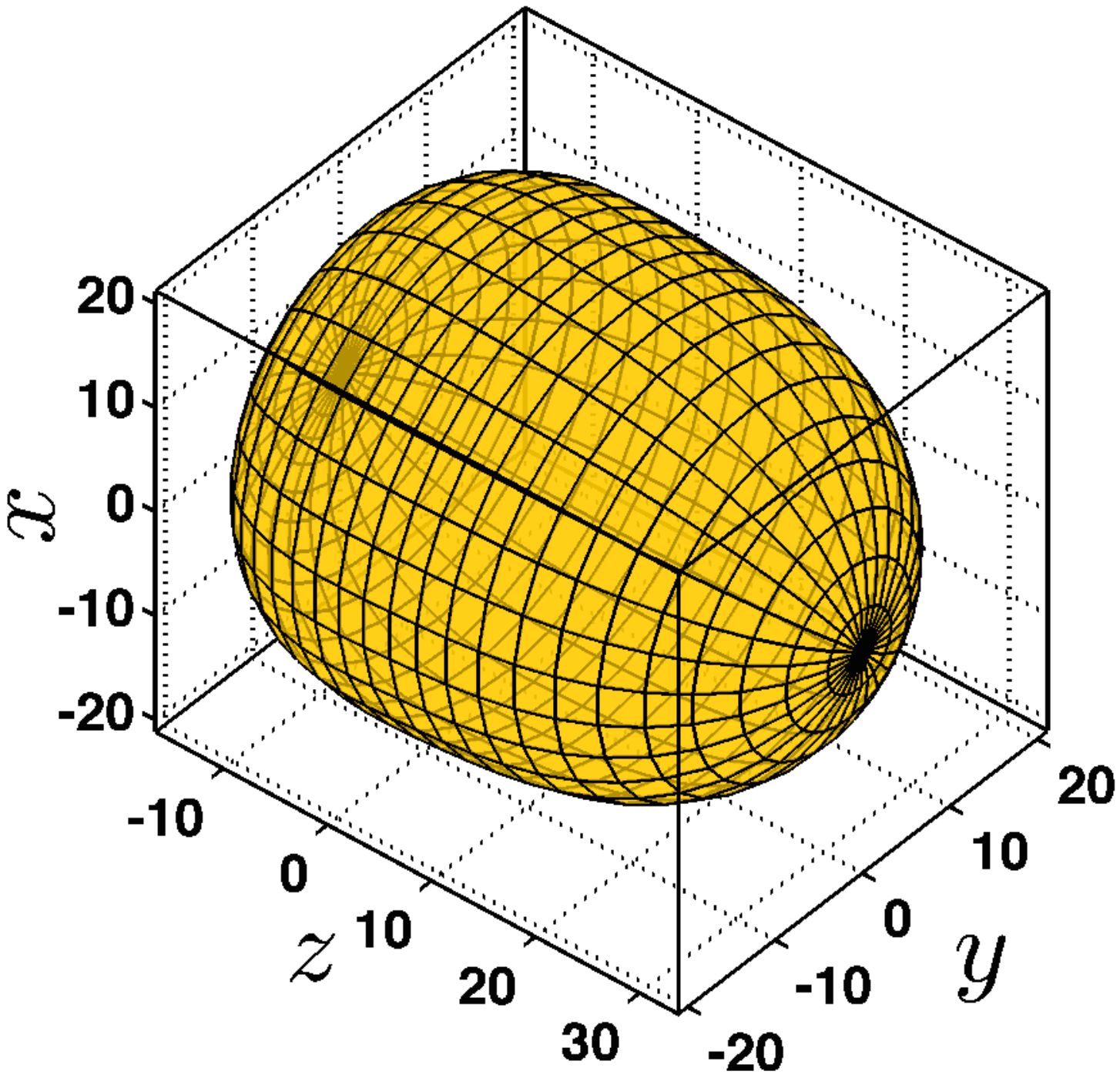}}
  \subfigure[$S_2$]{\includegraphics[height=1.80in]{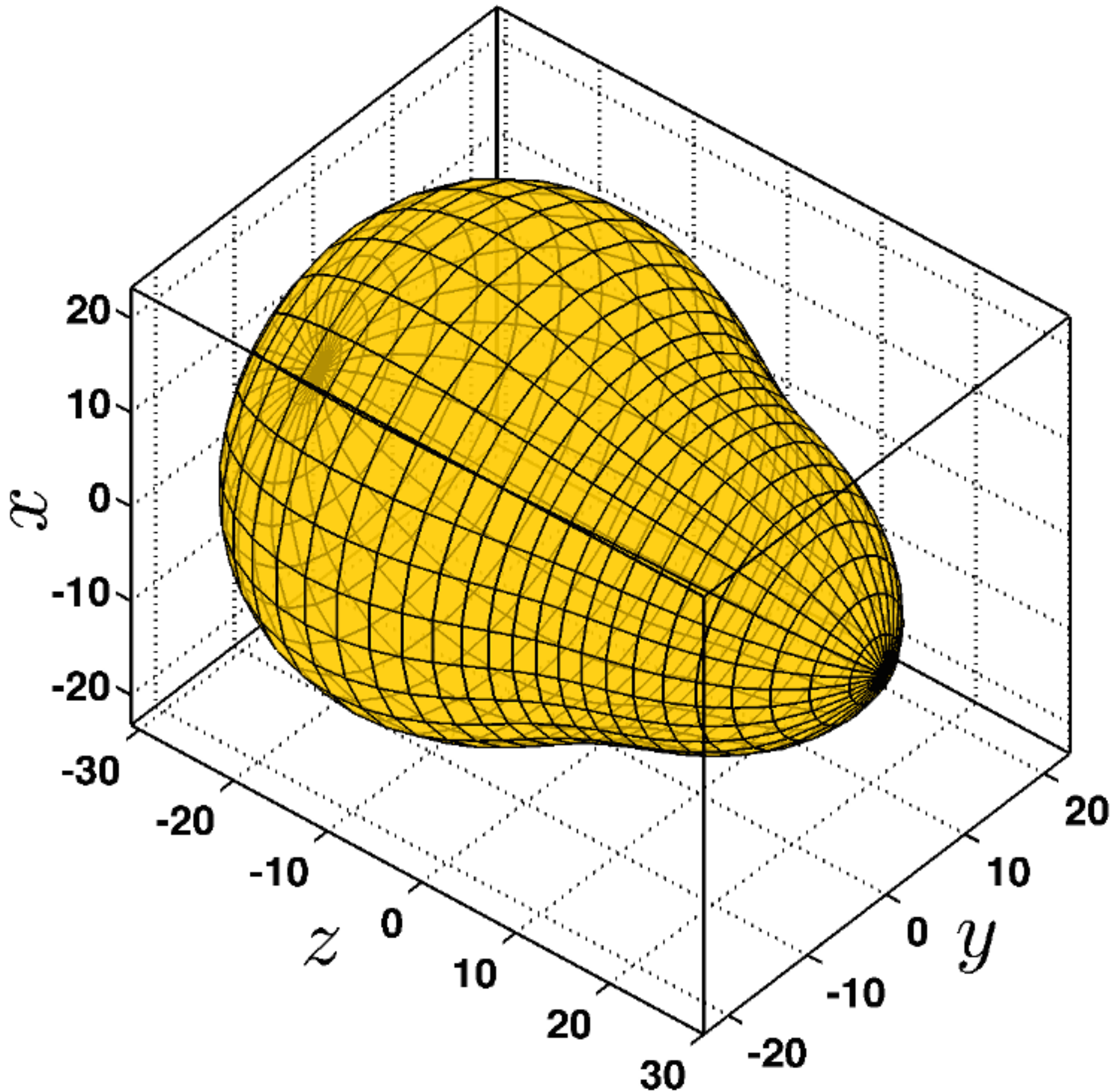}}
  \subfigure[$S_3$]{\includegraphics[height=1.80in]{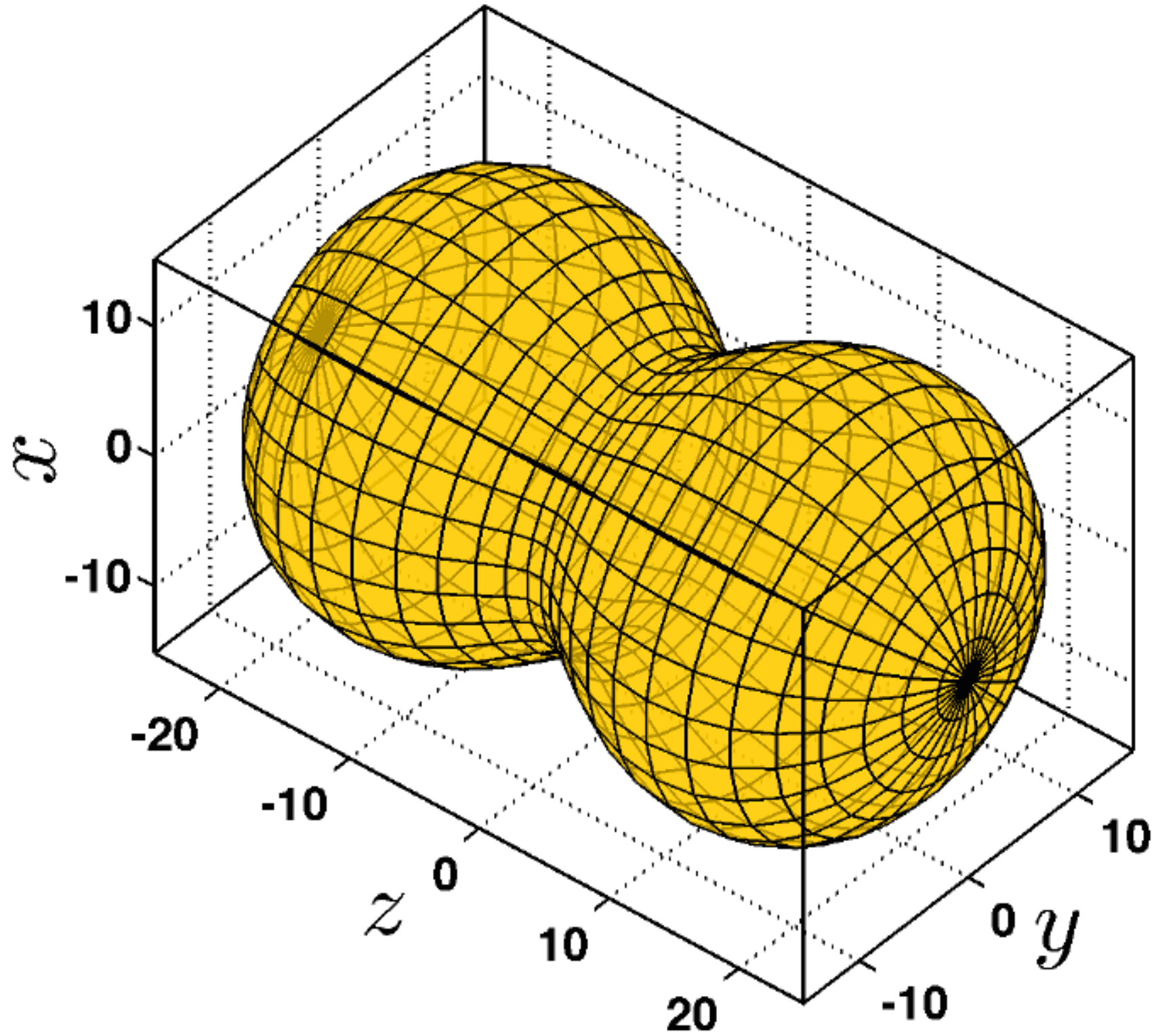}}
  \mcaption{fig:S-shapes}{Shapes of cell periphery}{Three different
    shapes of the cell periphery used in this study corresponding to
    surfaces given by expressions in \pr{eqn:S-shapes}.}
\end{figure}

\Pr{fig:S1,fig:S2,fig:S3} show snapshots of the PNC and its
centrosomal MTs as the structure moves and rotates within the $S_1$,
$S_2$, and $S_3$ cell geometries, respectively, for simulations of
both the cortical pushing and cytoplasmic pulling models. Also plotted
for each case is the temporal evolution of PNC position and the angle
$\alpha^{\star}$ between its centrosomal axis (the vector connecting
the two MTOC) and the $\zhat$ axis. The snapshots correspond to
different stages in the dynamics: soon after migration begins near the
posterior end of cell, after the centering of the PNC, and after the
rotation of PNC to its steady-state alignment with the $\zhat$
axis. We see that the PNC eventually aligns with the $\zhat$ axis in
all three geometries and for both models. The conformation of the MTs
is very different between the two models. In the cortical pushing
simulations several MTs are strongly buckled due to the compressive
forces at the periphery, while in the cytoplasmic pulling simulations
MTs are under extensile loads and so remained nearly undeformed.
\begin{figure}[!tb]
  \centering
  \small
  \subfigure[$t=  2\usep\minute$\label{sfg:S1_push_1}]{
    \includegraphics[height=1.38in]{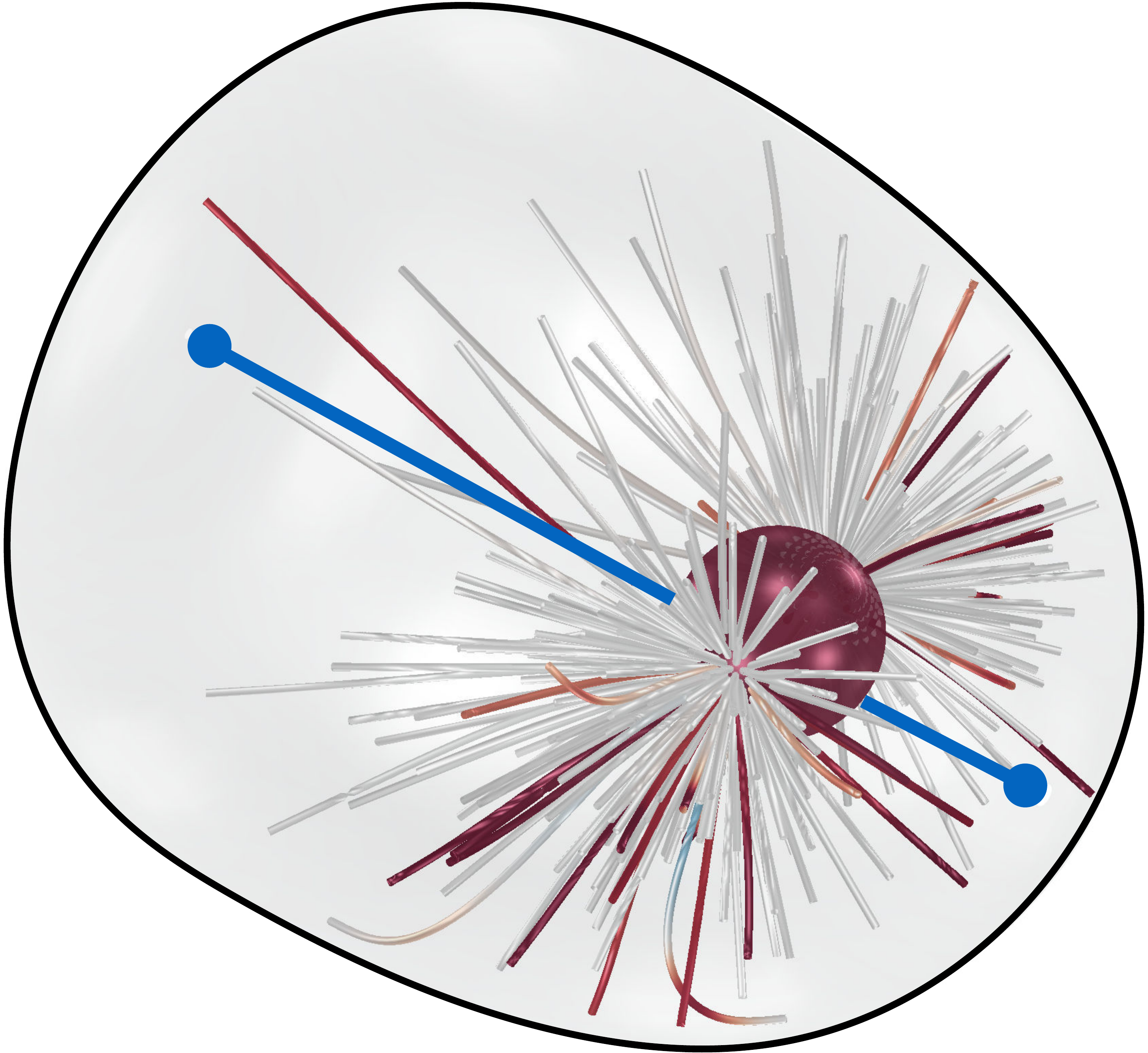}}
  \subfigure[$t= 200\usep\minute$\label{sfg:S1_push_2}]{
    \includegraphics[height=1.38in]{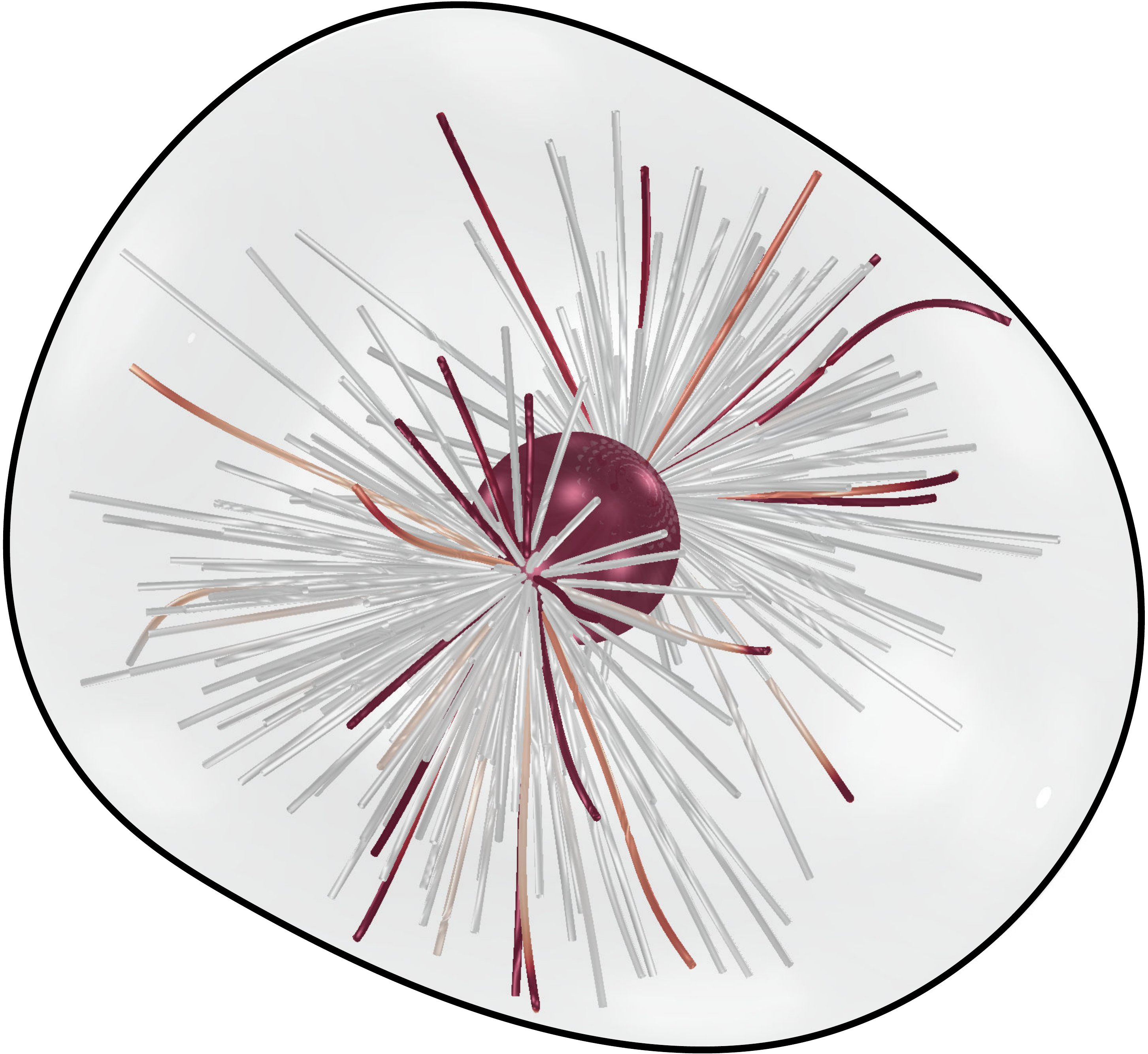}}
  \subfigure[$t=960\usep\minute$\label{sfg:S1_push_3}]{
    \includegraphics[height=1.38in]{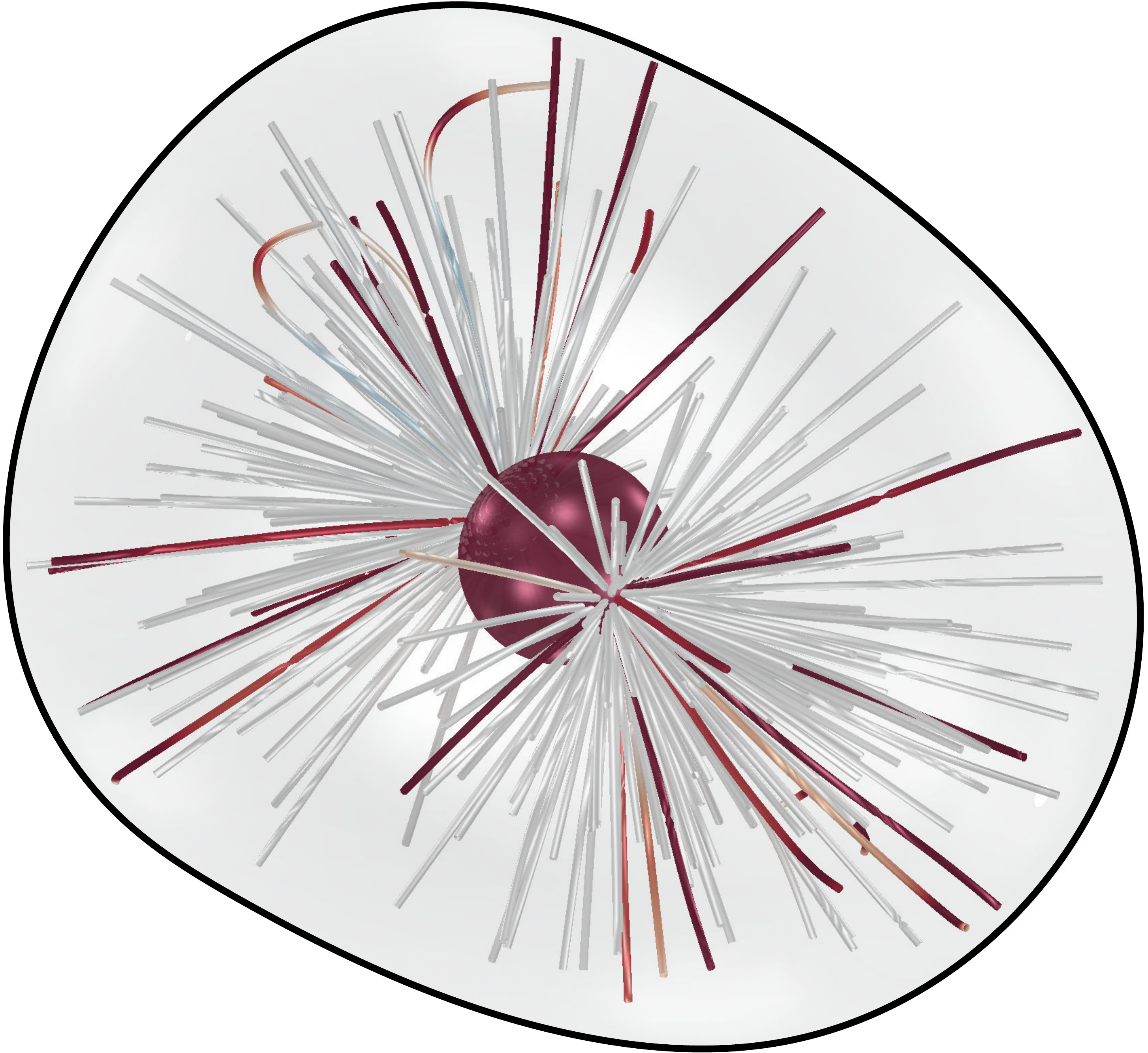}}
  \subfigure[\label{sfg:S1_pos_push}]{
    \setlength\figurewidth{1.0in}%
    \setlength\figureheight{1.05in}%
    \hspace*{-7pt}%
    \includepgf{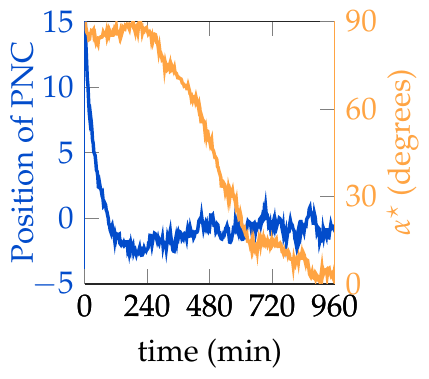}%
    \hspace*{-5pt}%
  }\\
  \subfigure[$t=  2\usep\minute$\label{sfg:S1_pull_1}]{\includegraphics[height=1.38in]{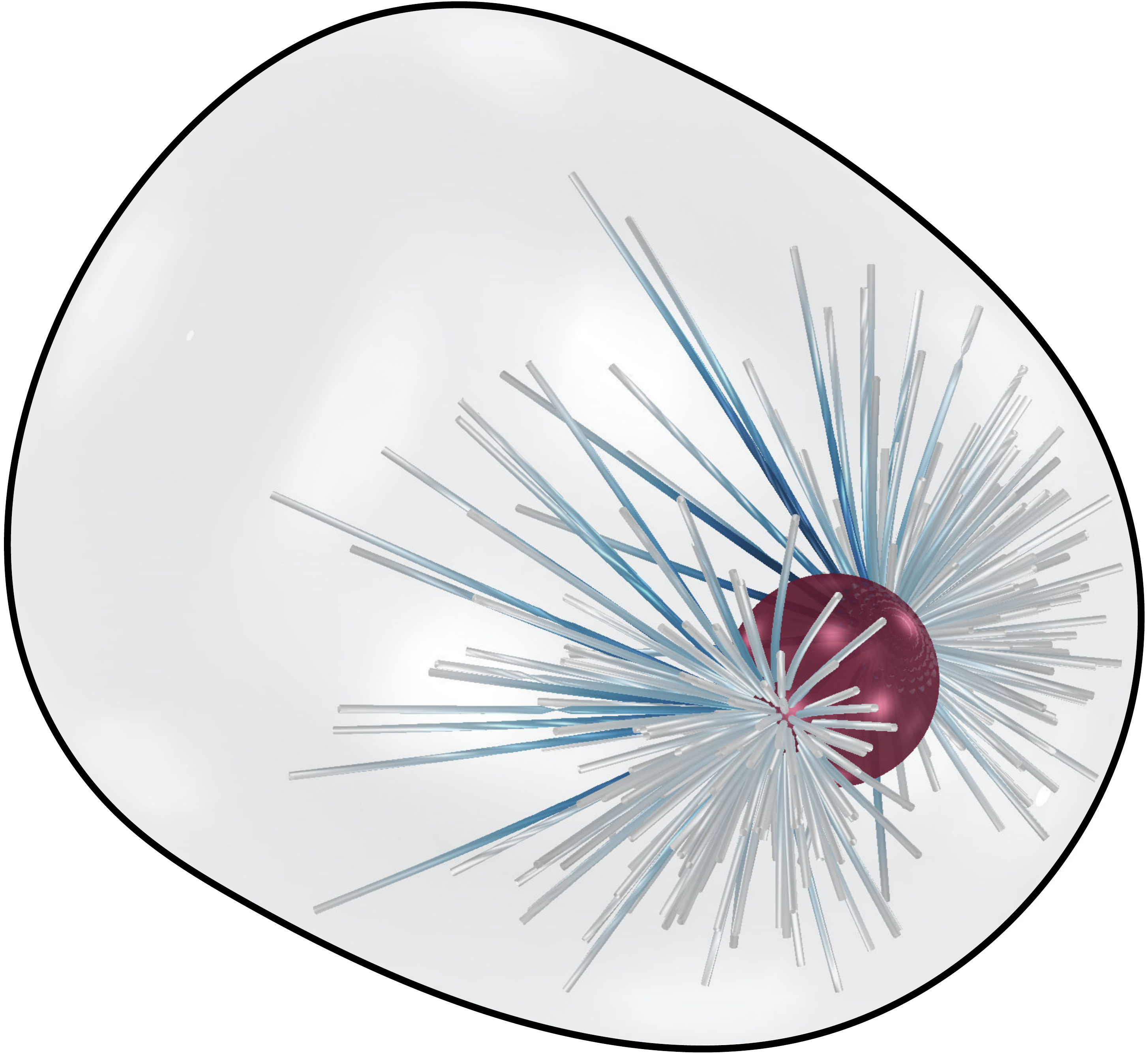}}
  \subfigure[$t= 15\usep\minute$\label{sfg:S1_pull_2}]{\includegraphics[height=1.38in]{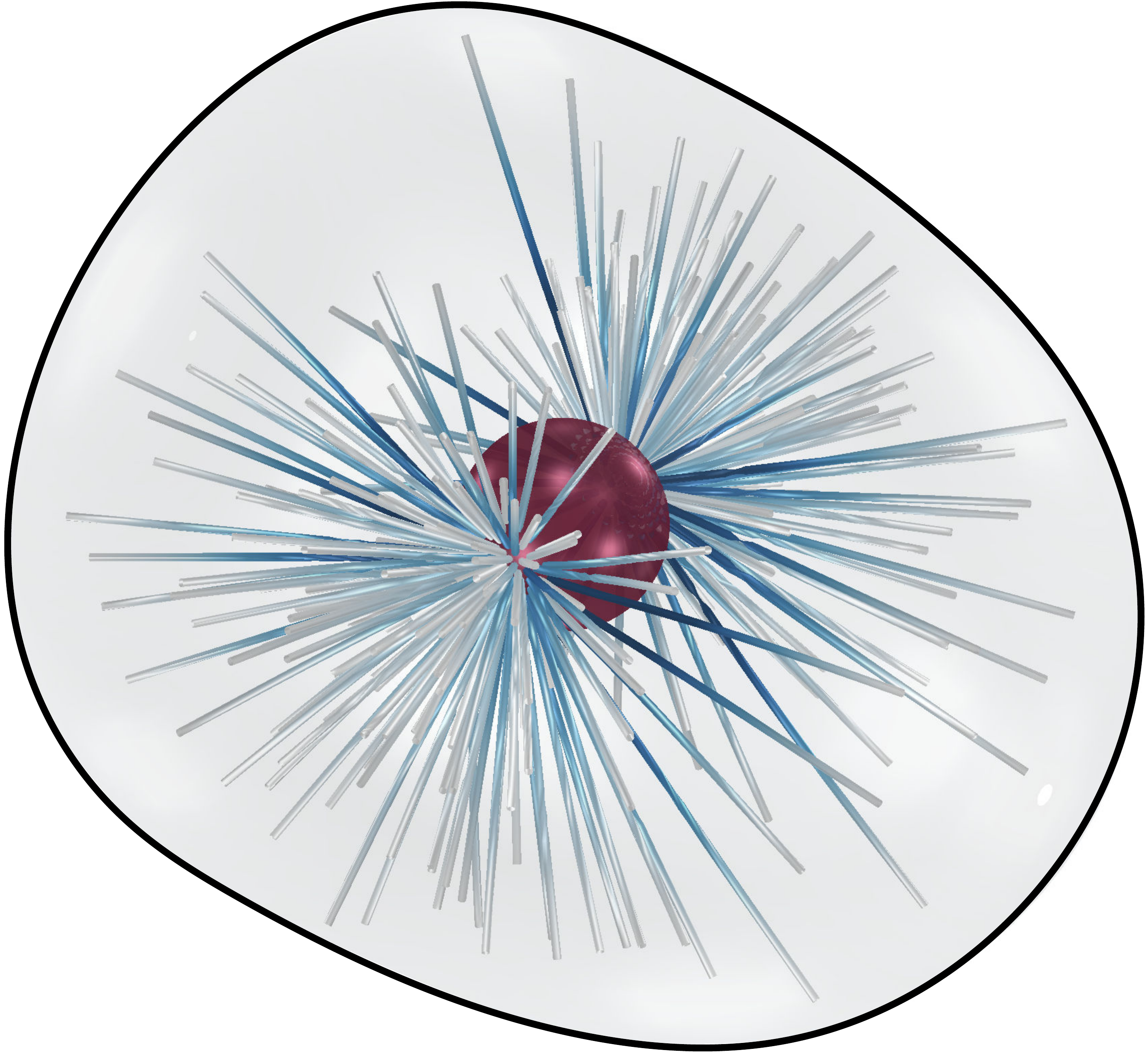}}
  \subfigure[$t=240\usep\minute$\label{sfg:S1_pull_3}]{\includegraphics[height=1.38in]{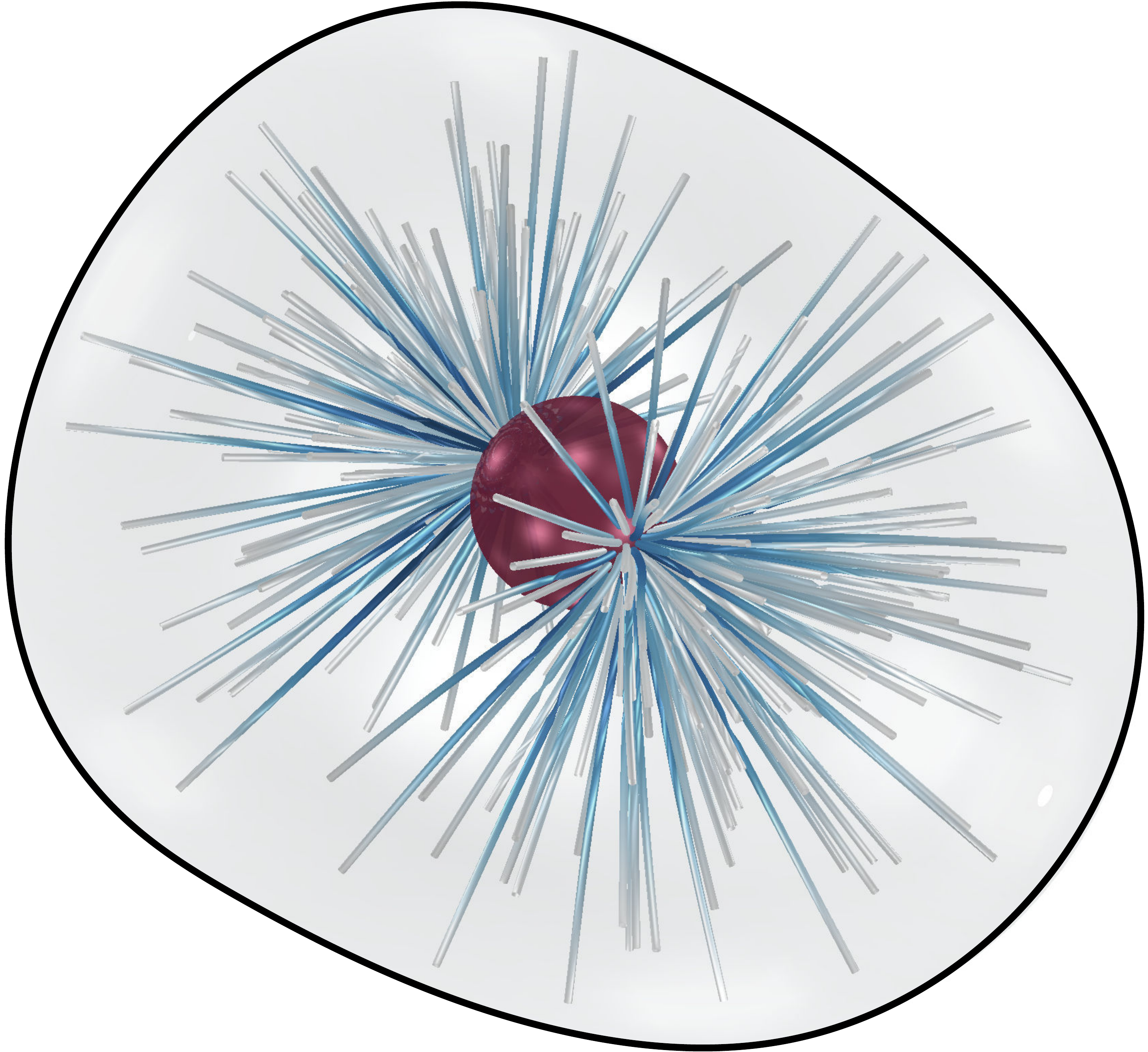}}
  \subfigure[\label{sfg:S1_pos_pull}]{
    \setlength\figurewidth{1.0in}%
    \setlength\figureheight{1.05in}%
    \hspace*{-7pt}%
    \includepgf{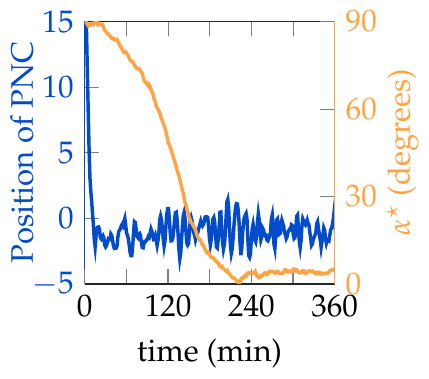}%
    \hspace*{-5pt}%
  }
    \mcaption{fig:S1}{Snapshots of centering process within $S_1$}{(a-c)
    cortical pushing model; (e-g) cytoplasmic pulling model; (d) and
    (h) the $\zhat$-component of the position of the PNC and angle
    between MTOC and the $\zhat$ axis as a function of time, for
    cortical pushing and cytoplasmic pulling models respectively. The
    fibers in all the simulations are colored based on the tension
    along their length. The colors change from dark red (compressive
    tension) to dark blue (extensile tension) and the white color
    denotes no tension. For all the simulations of the PNC migration,
    $N_F=300$. The turnover time in cortical pushing model was taken
    as $\tau_\lbl{cat}=4$ and the number of motors per unit length of
    the MTs in cytoplasmic model was taken as
    $n_{dyn}=0.04\usep(\microm)^{-1}$.}
\end{figure}

\begin{figure}[!tb]
  \centering
  \small
  \subfigure[$t=  2\usep\minute$\label{sfg:S2push_1}]{\includegraphics[height=1.28in]{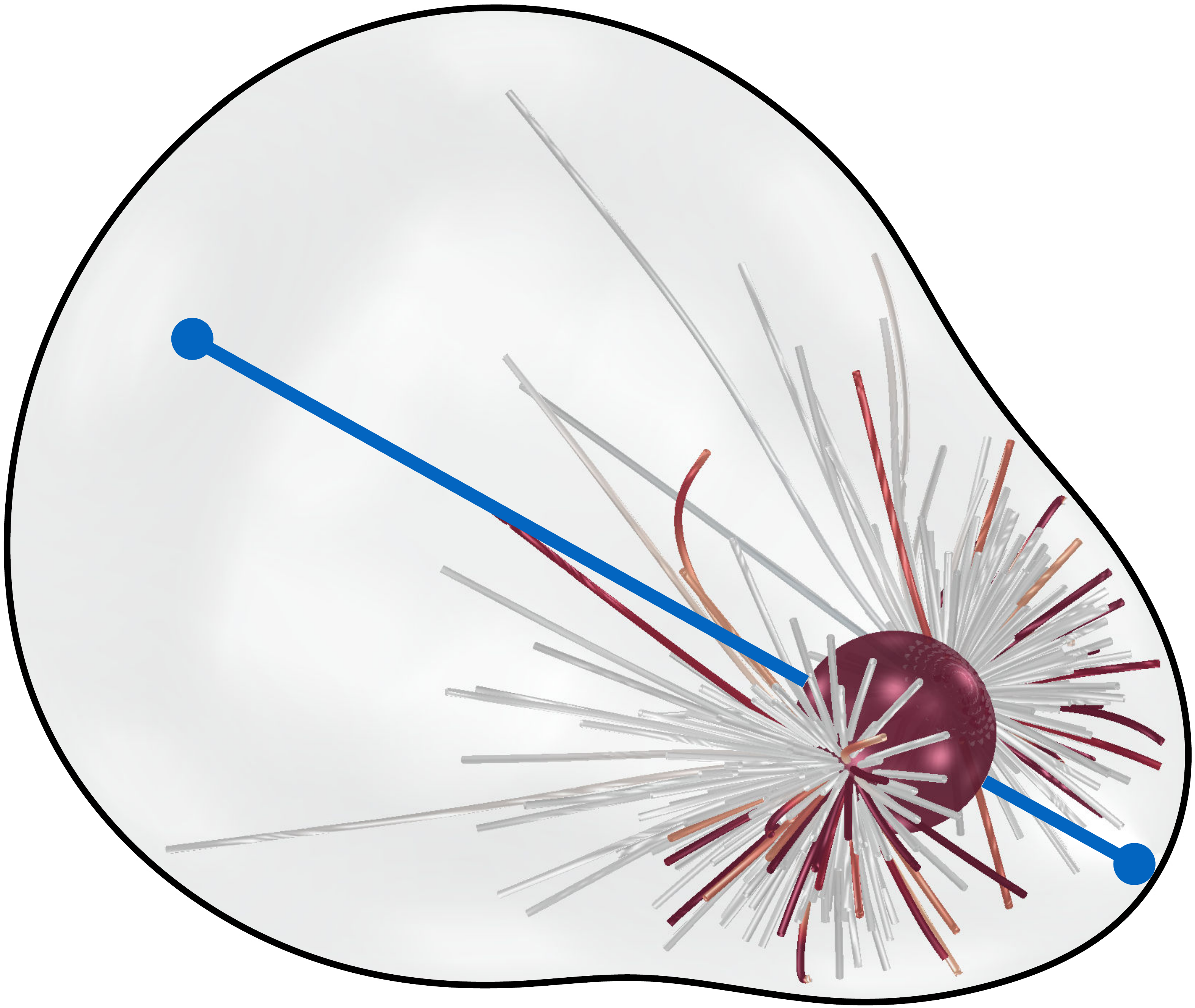}}
  \subfigure[$t= 40\usep\minute$\label{sfg:S2_push_2}]{\includegraphics[height=1.28in]{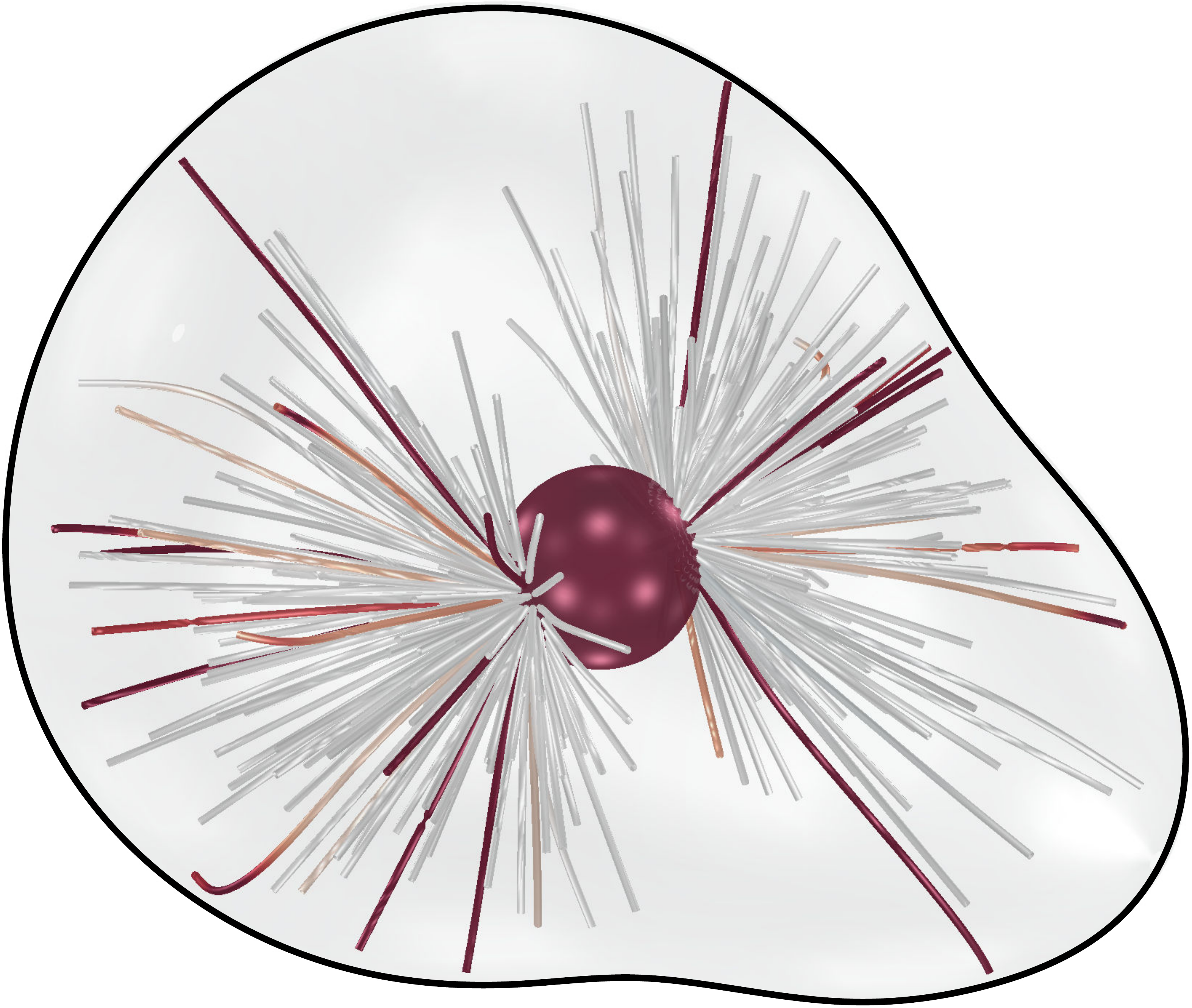}}
  \subfigure[$t=240\usep\minute$\label{sfg:S2_push_3}]{\includegraphics[height=1.28in]{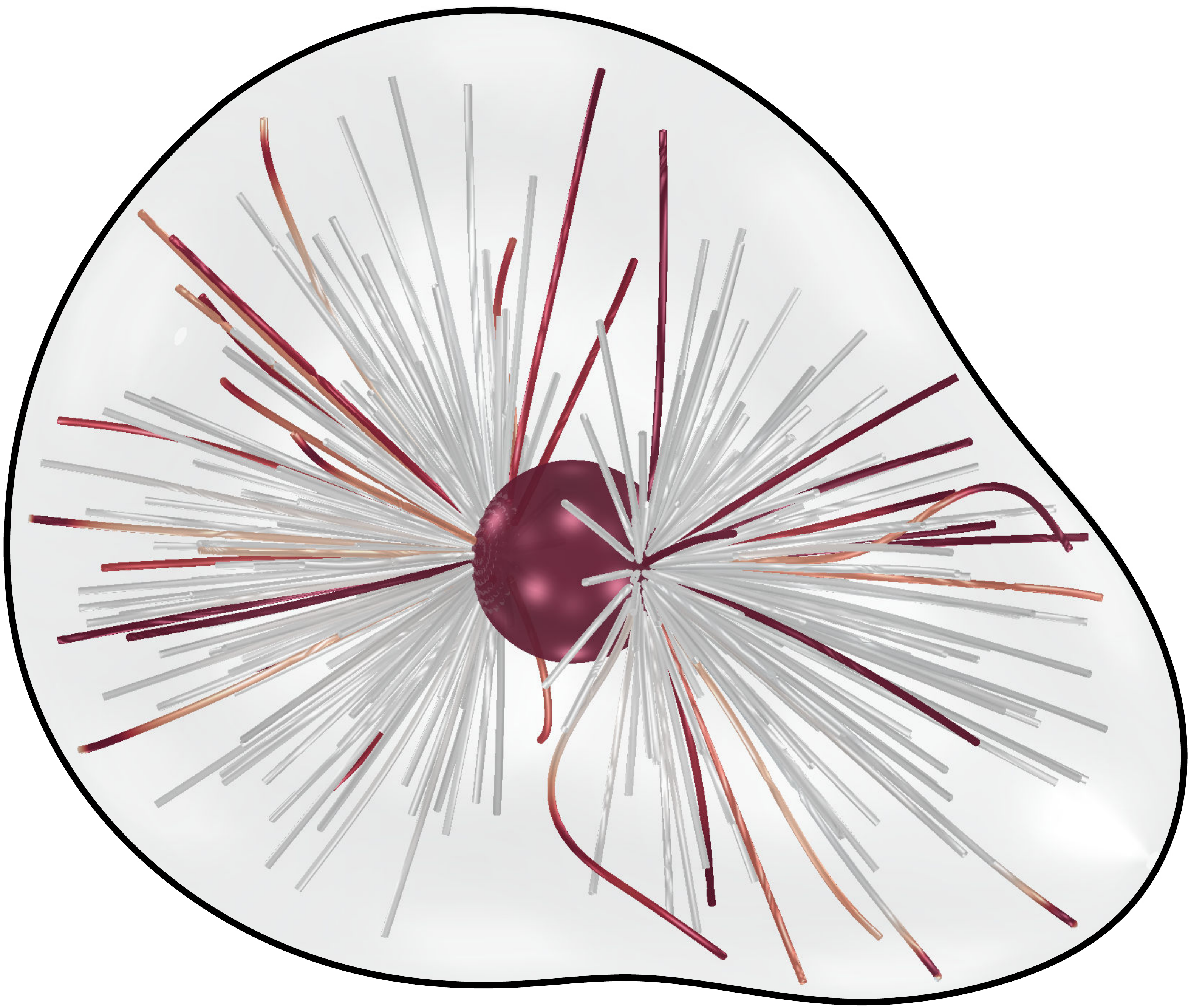}}
  \subfigure[\label{sfg:S2_pos_push}]{
    \setlength\figurewidth{1.0in}%
    \setlength\figureheight{.95in}%
    \hspace*{-7pt}%
    \includepgf{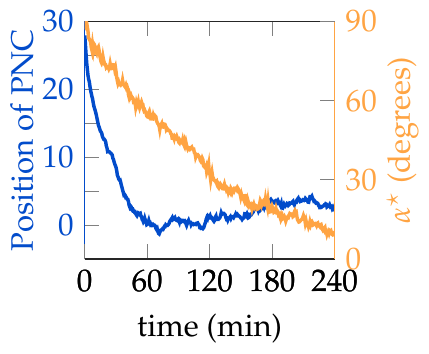}%
    \hspace*{-5pt}%
  }\\
  \subfigure[$t=  2\usep\minute$\label{sfg:S2_pull_1}]{\includegraphics[height=1.28in]{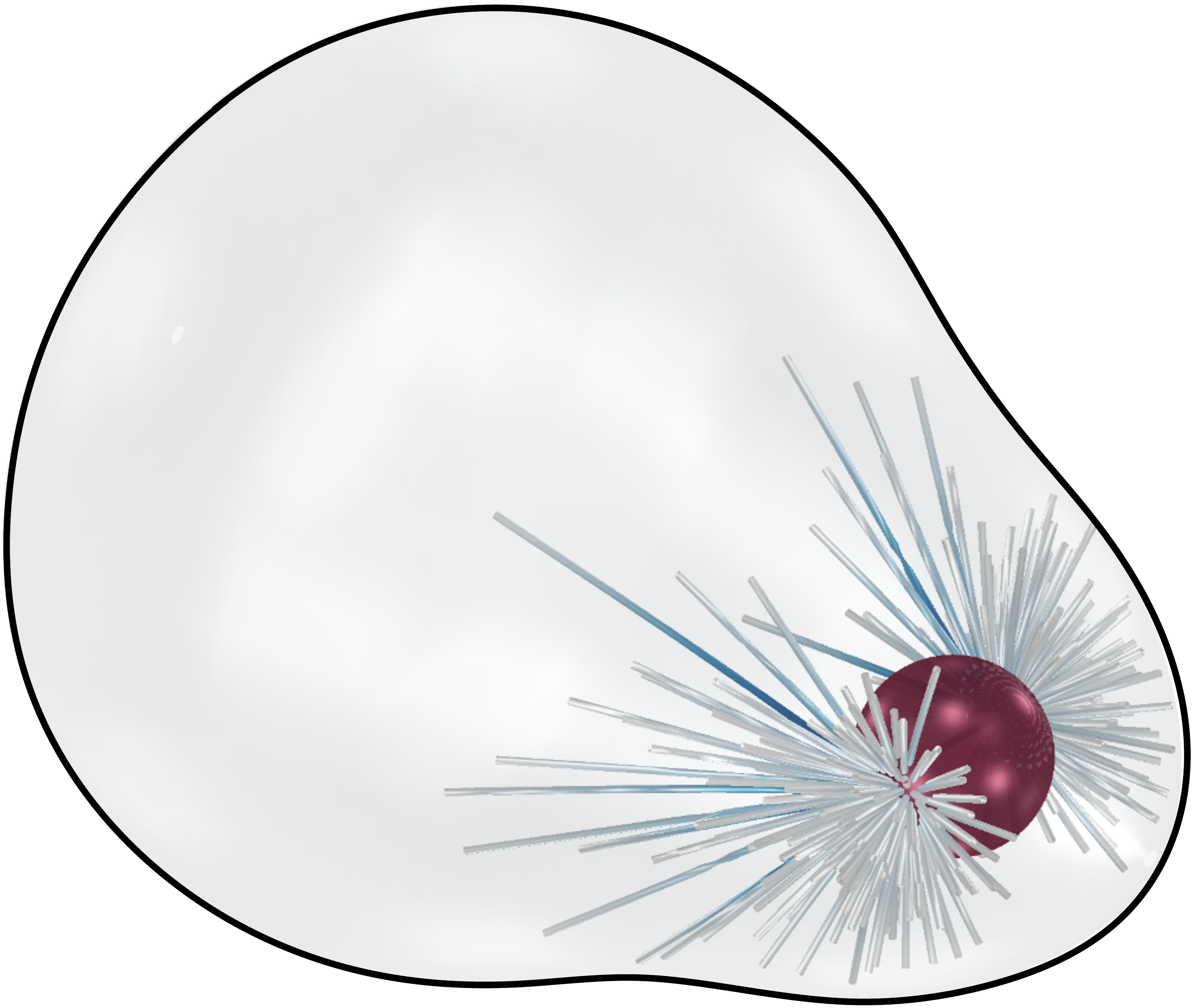}}
  \subfigure[$t= 25\usep\minute$\label{sfg:S2_pull_2}]{\includegraphics[height=1.28in]{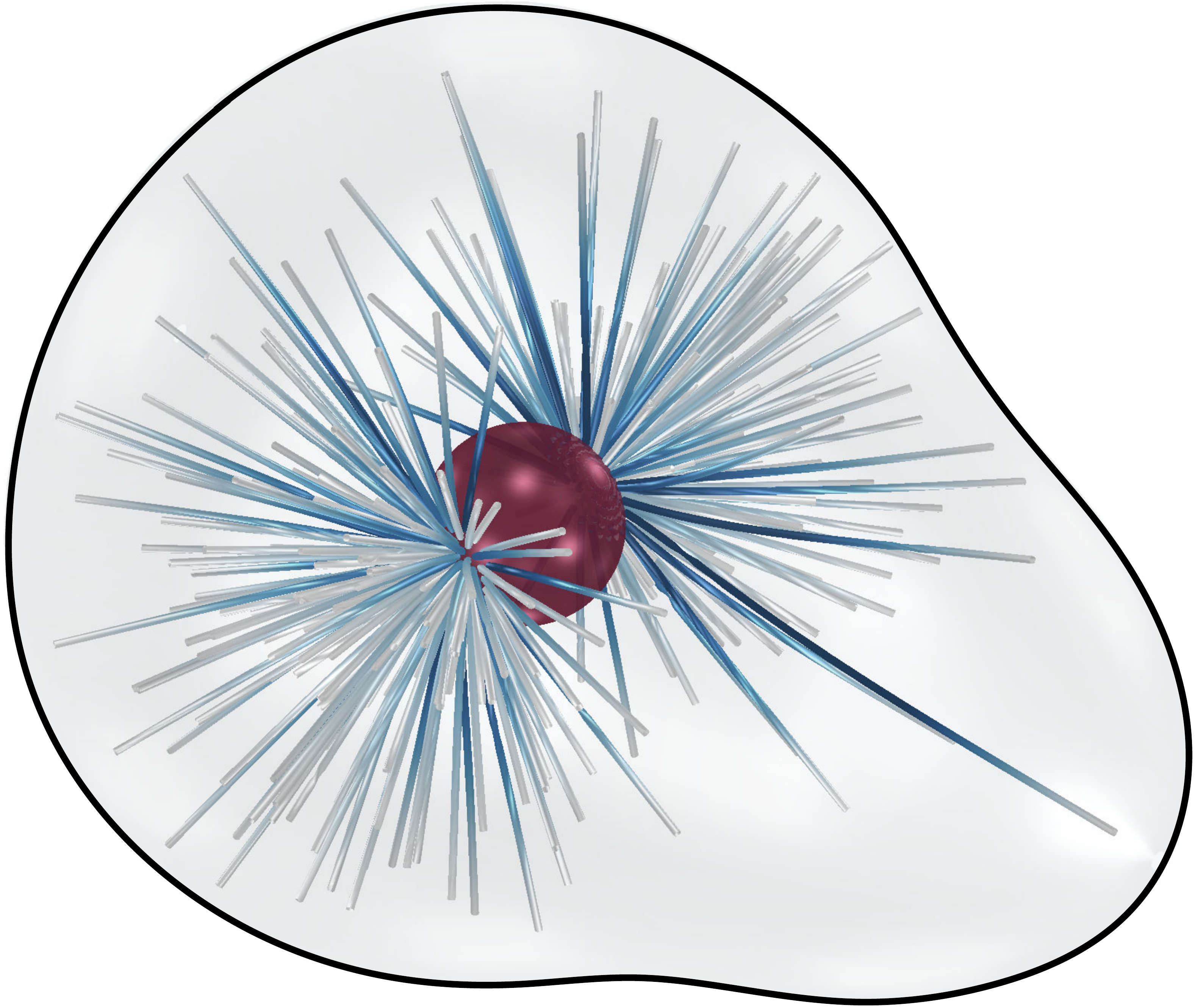}}
  \subfigure[$t=240\usep\minute$\label{sfg:S2_pull_3}]{\includegraphics[height=1.28in]{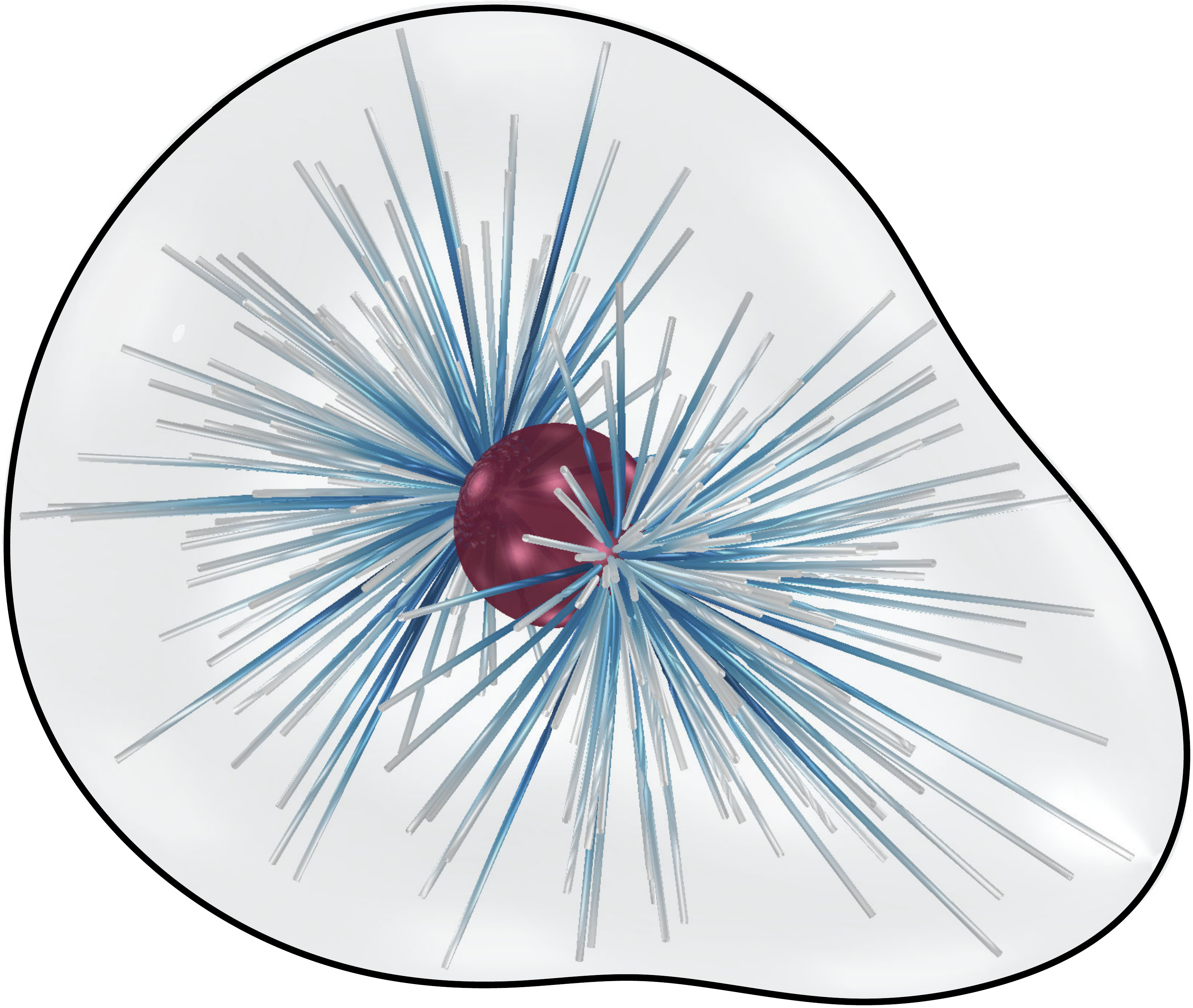}}
  \subfigure[\label{sfg:S2_pos_pull}]{
    \setlength\figurewidth{1.0in}%
    \setlength\figureheight{.95in}%
    \hspace*{-7pt}%
    \includepgf{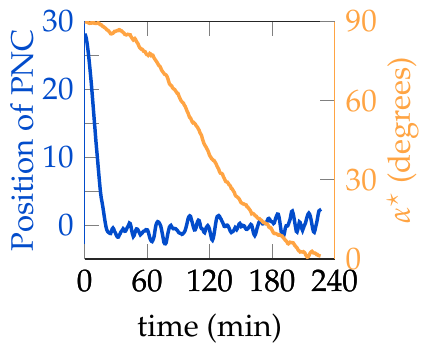}%
    \hspace*{-5pt}%
  }
  \mcaption{fig:S2}{Snapshots of centering process within $S_2$}{
  The descriptions and the parameters used are identical to the one provided
  in the caption of \pr{fig:S1}.}
\end{figure}

\begin{figure}[!tb]
  \centering
  \small
  \subfigure[$t=  2\usep\minute$\label{sfg:S3_push_1}]{\includegraphics[height=1.28in]{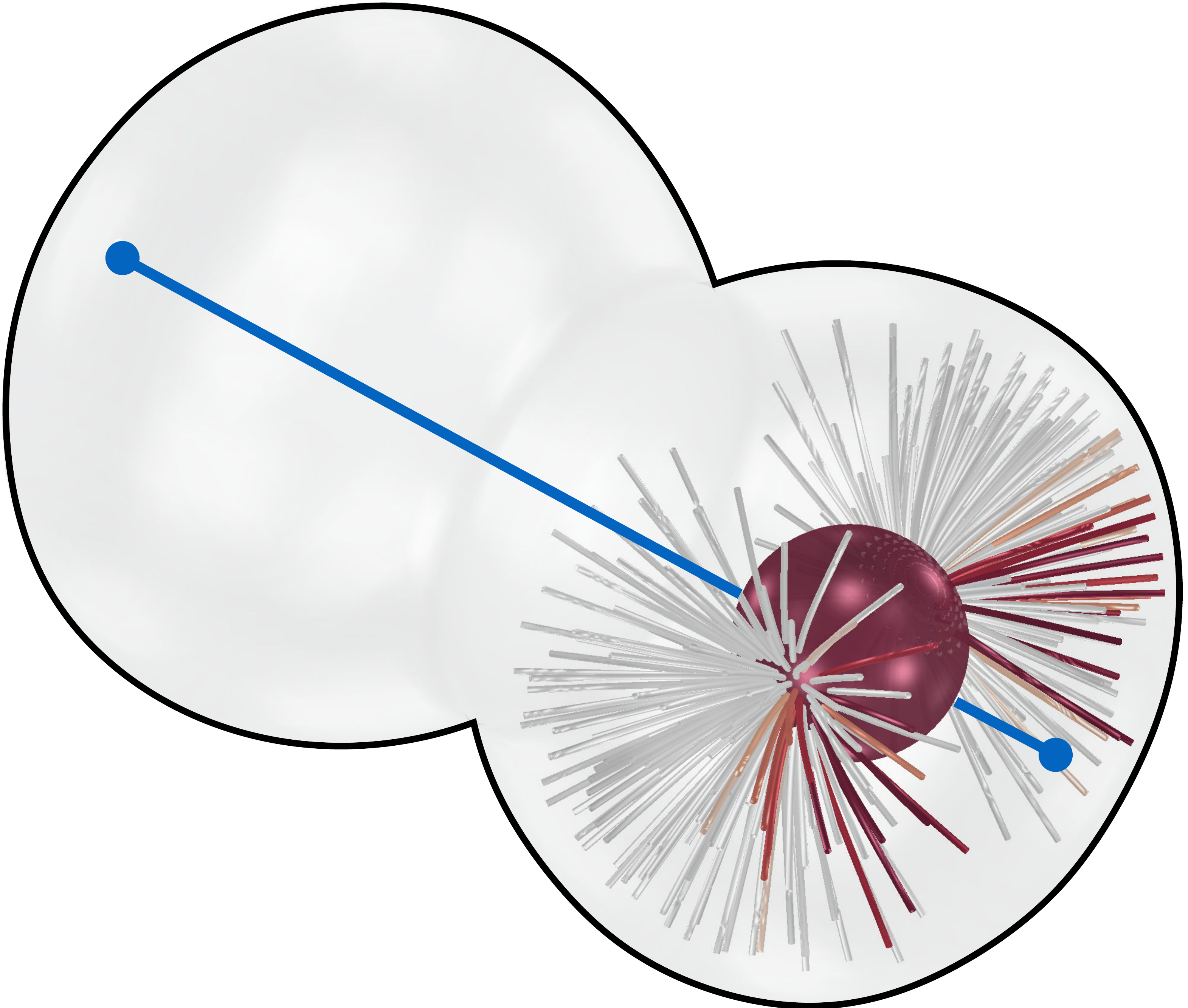}}
  \subfigure[$t= 90\usep\minute$\label{sfg:S3_push_2}]{\includegraphics[height=1.28in]{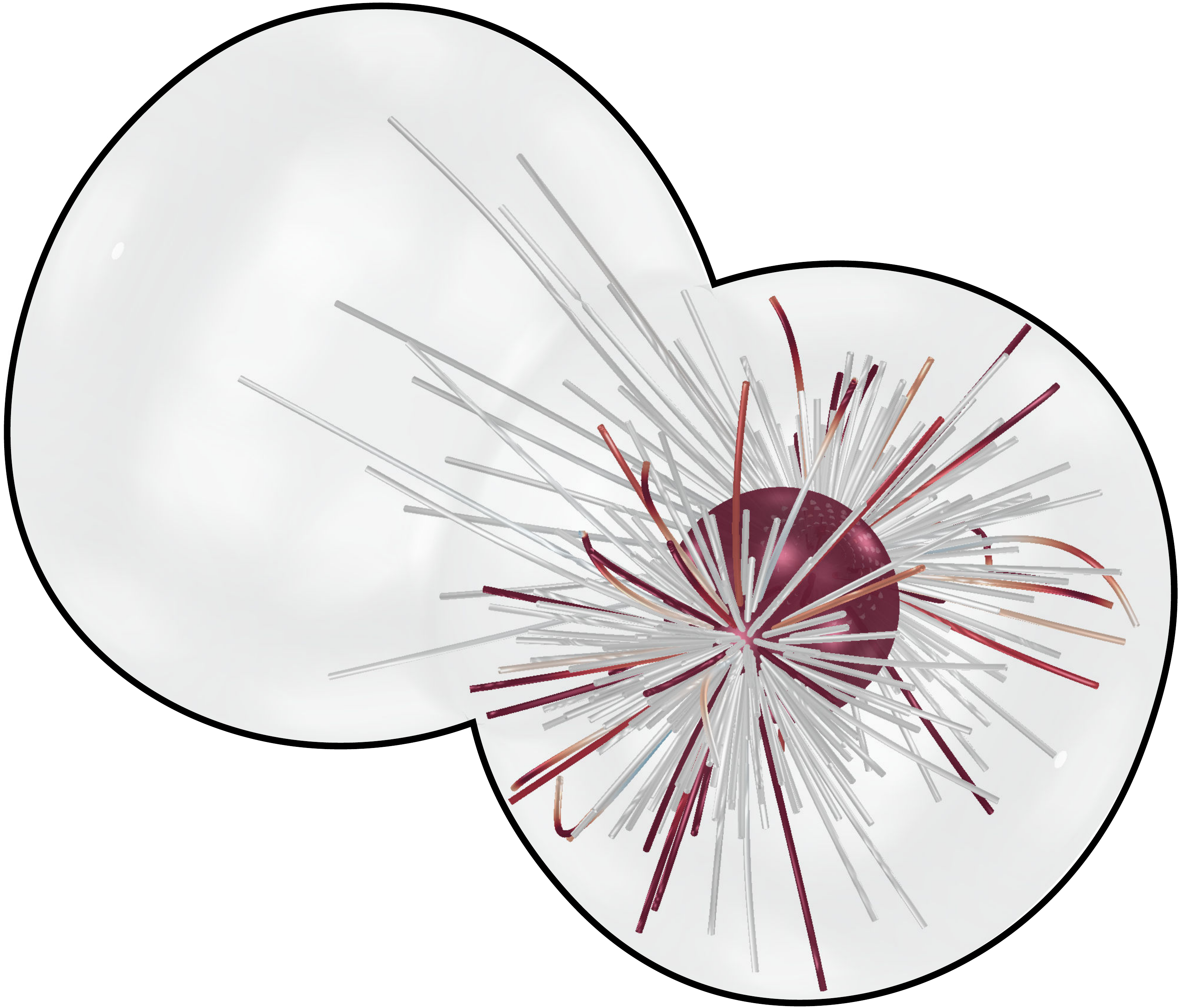}}
  \subfigure[$t=180\usep\minute$\label{sfg:S3_push_3}]{\includegraphics[height=1.28in]{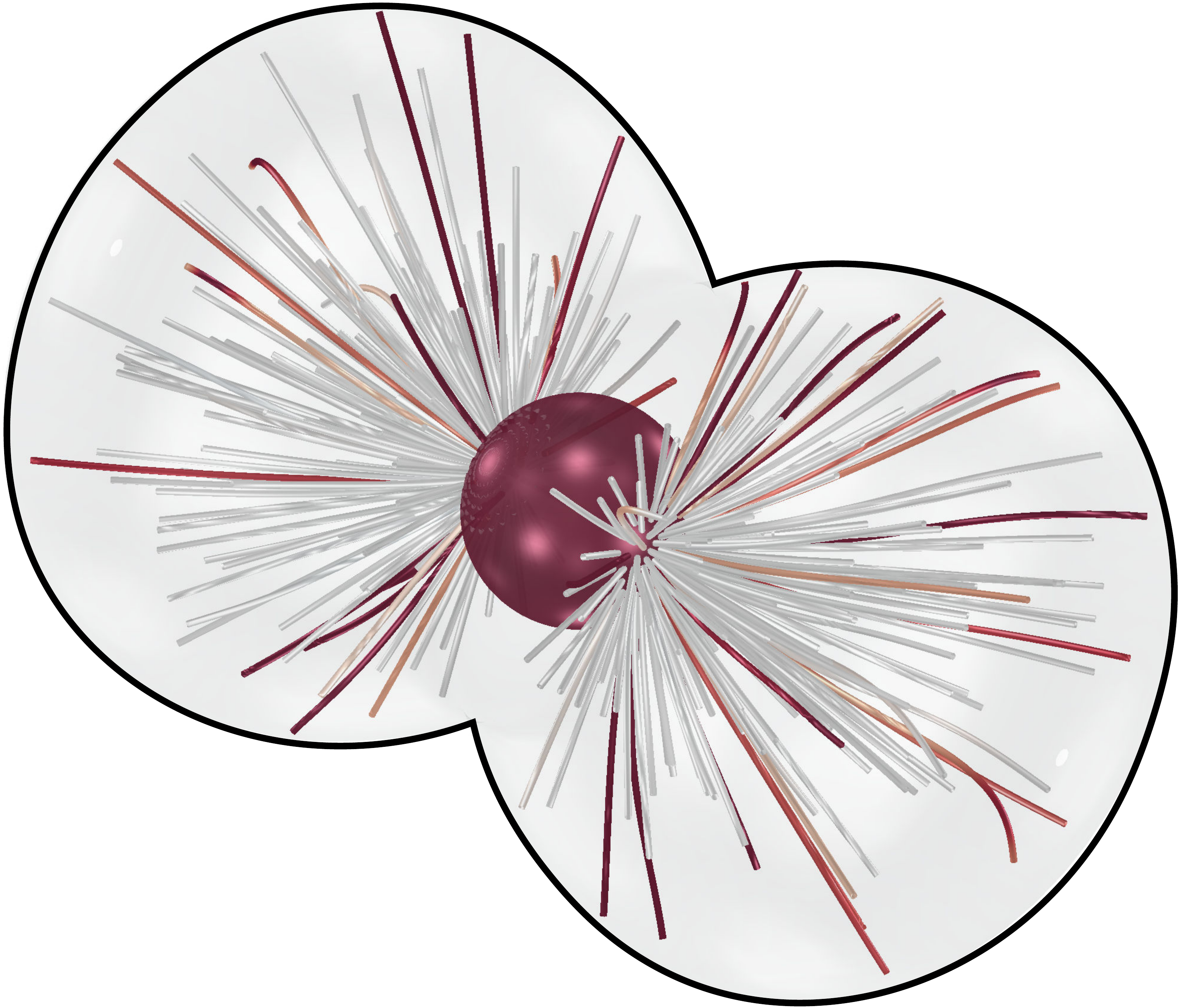}}
  \subfigure[\label{sfg:S3_pos_push}]{
    \setlength\figurewidth{1.0in}%
    \setlength\figureheight{.95in}%
    \hspace*{-7pt}%
    \includepgf{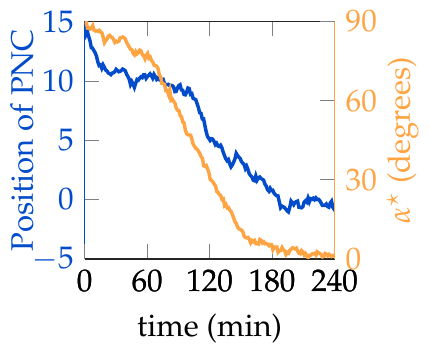}%
    \hspace*{-5pt}%
  }\\
  \subfigure[$t=  2\usep\minute$\label{sfg:S3_pull_1}]{\includegraphics[height=1.28in]{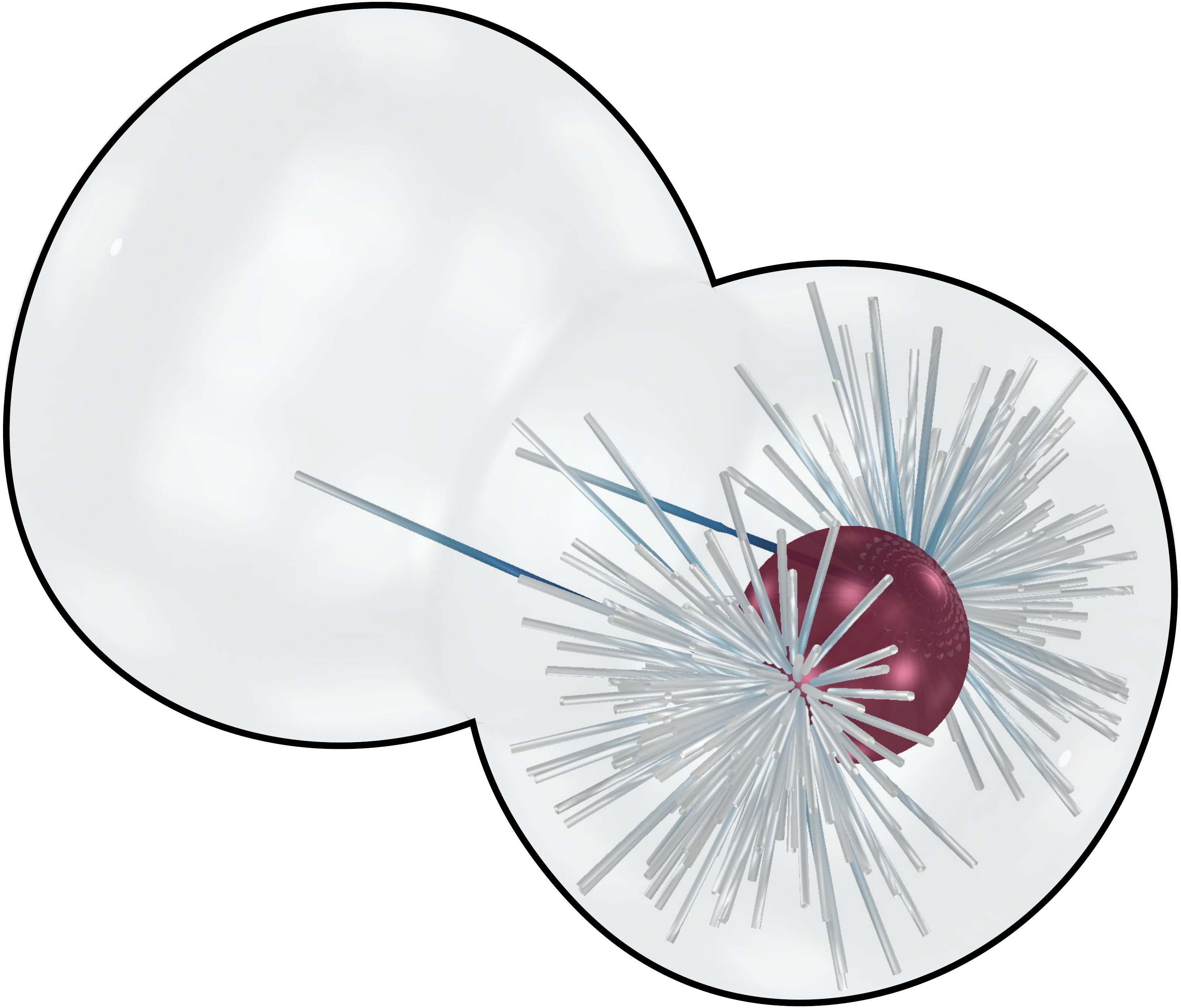}}
  \subfigure[$t= 90\usep\minute$\label{sfg:S3_pull_2}]{\includegraphics[height=1.28in]{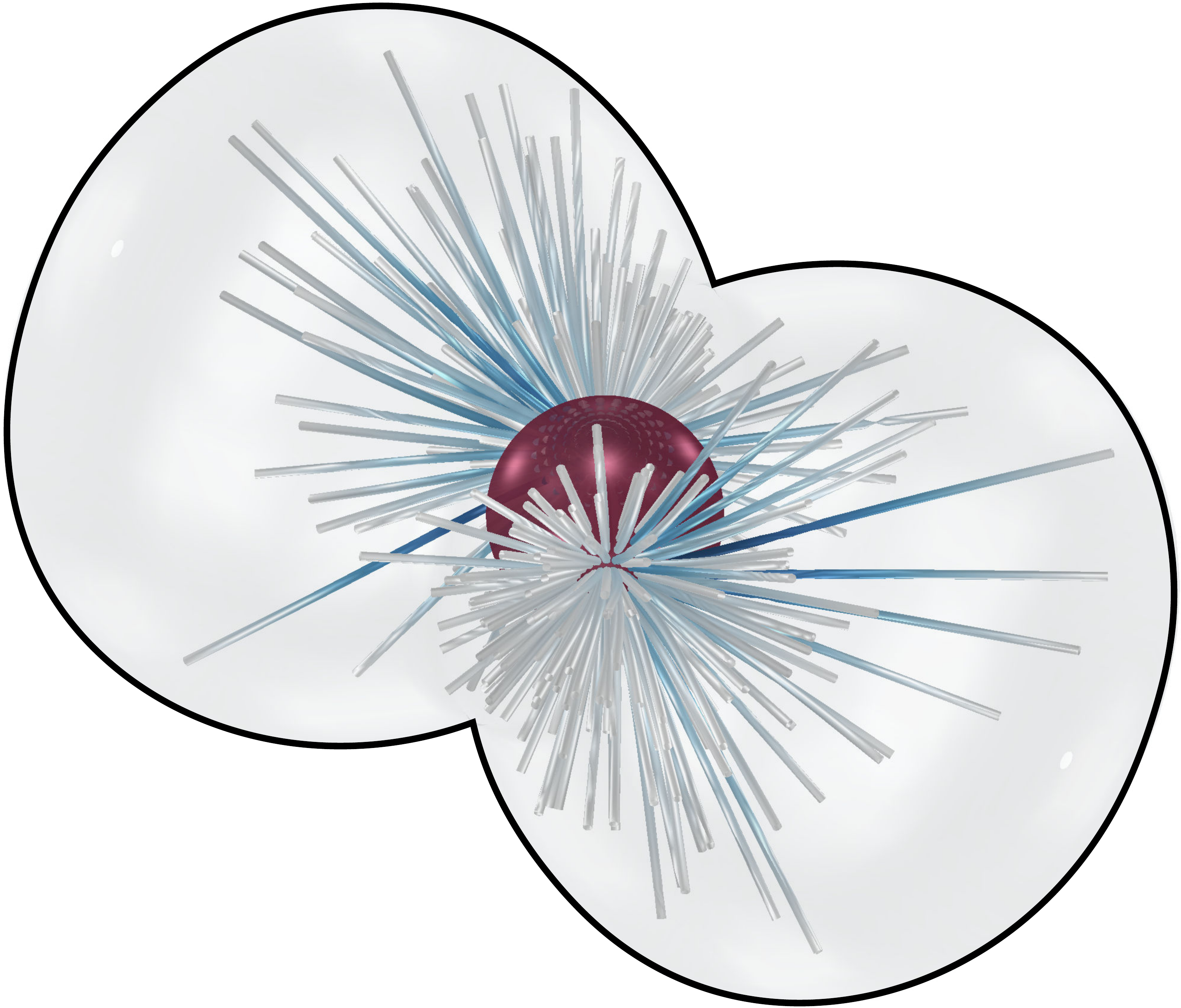}}
  \subfigure[$t=180\usep\minute$\label{sfg:S3_pull_3}]{\includegraphics[height=1.28in]{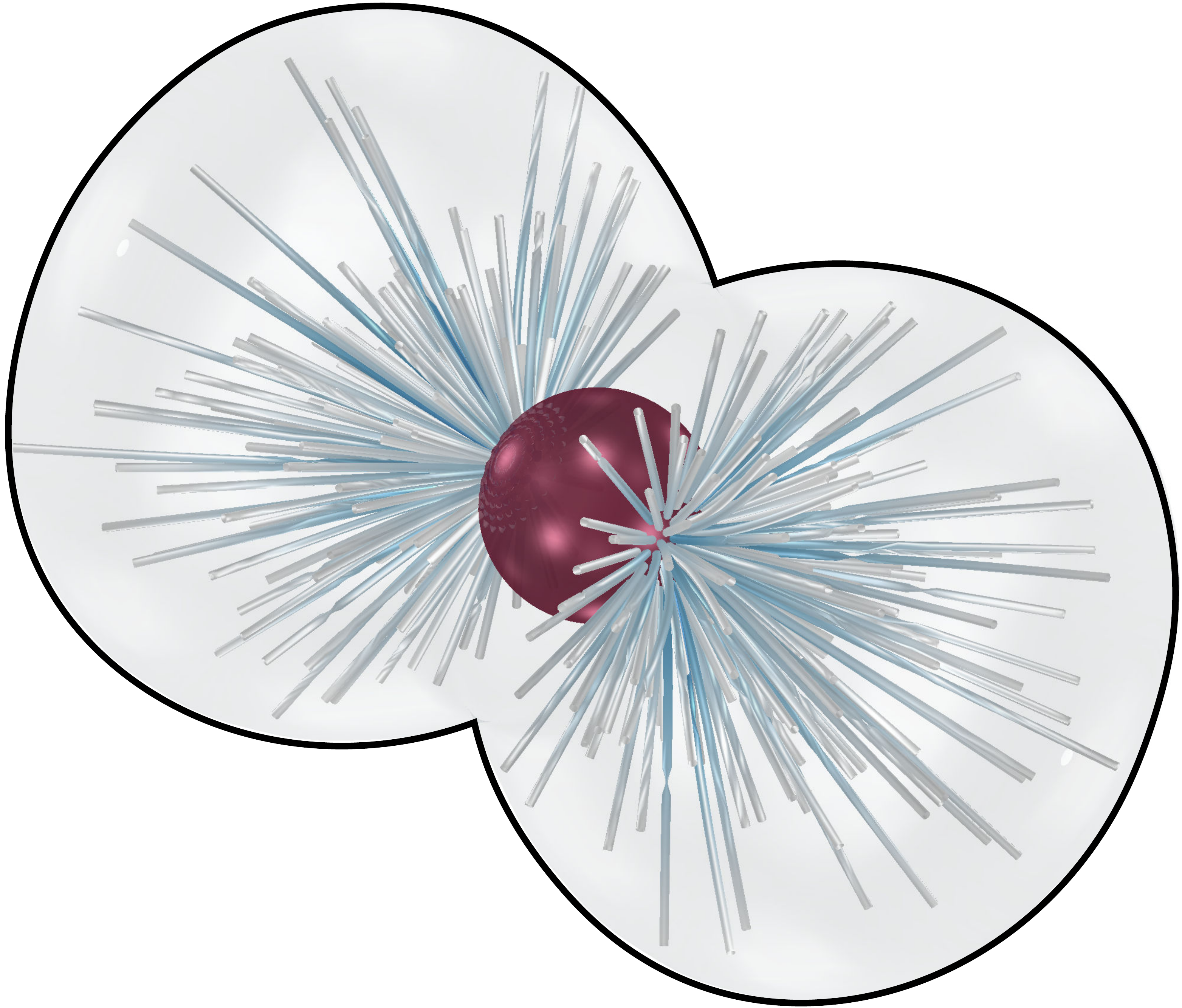}}
  \subfigure[\label{sfg:S3_pos_pull}]{
    \hspace{.02in}
    \setlength\figurewidth{1.0in}%
    \setlength\figureheight{.95in}%
    \hspace*{-7pt}%
    \includepgf{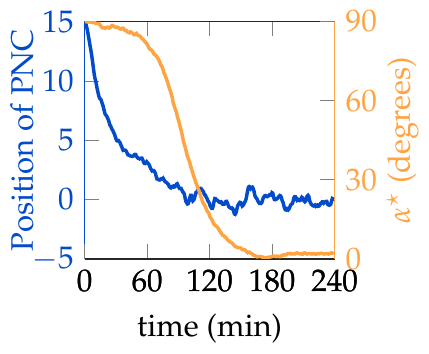}%
    \hspace*{-5pt}%
  }
  \mcaption{fig:S3}{Snapshots of centering within $S_3$}{
    For details refer to the caption of \pr{fig:S1}.}
\end{figure}

The most significant dynamical consequence of changing the cell
geometry is the large variation in centering and rotation times
between geometries and force-transduction models. For example, in
geometry $S_{1}$ the rotation times ($960\usep\minute$, and
$150\usep\minute$) are roughly 5 and 10 times larger than the
centering time ($200\usep\minute$ and $15\usep\minute$) in the
cortical pushing and cytoplasmic pulling simulations, respectively.

\emph{In vivo} observations of the first cell division in
normally-shaped \celegans\ embryos, however, show that rotation and
positioning occur on the same time-scales \citep{Kimura2007}. Further,
while cytoplasmic pulling and cortical pushing mechanisms do not
predict PNC rotation in spherically shaped eggshells, experiments on
such geometries show that the PNC can still properly center and align
with AP-axis when the embryonic cell is perturbed to be nearly
spherical \cite{HW1987,Tsou2002}. These observations suggest that
other active mechanisms, such as cortical pulling where the necessary
asymmetry is not induced by the geometry but by the anisotropic
distribution of cortical dyneins on the periphery, can also be
involved in pronuclear migration.

Unlike the $S_1$ case, the centering and rotation in the cortical
pushing model in the $S_{2}$ geometry occurs simultaneously
($240\usep\minute$) while in cytoplasmic pulling model the rotation
time ($240\usep\minute$) is again $10$ times larger than the
positioning time ($25\usep\minute$).  Finally in the $S_3$ geometry,
centering and rotation happen on the same time-scale
($240\usep\minute$) for both models. The faster rotation of the PNC in
case $S_{3}$ may be due to the fact that this is the most anisotropic
shape of the geometries.

A few more differences are evident. For example while the steady-state
positions of the PNC in the cortical pushing and cytoplasmic pulling
models are similar to each other in the $S_{1}$ and $S_3$ geometries,
the final positions for the two models differ by
$5\usep\microm$. Another observation is that in the $S_3$ cell shape,
using the cortical pushing model, the PNC becomes \emph{trapped} for
approximately one-half of the migration time in a dynamically
\emph{semi-stable} position $10\usep\microm$ from the center (shown in
\pr{sfg:S3_push_2}). Eventually the PNC escapes this dynamical cage as
it orients towards the $\zhat$-axis. To explain this behavior we first
note that the $S_3$ geometry is roughly two overlapped spheres. The
MTs that pass from the posterior sphere to the anterior are long and
their associated polymerization forces are weak. Thus, the net
cortical pushing force on the PNC, after it reaches to the center of
mass of the posterior sphere, is expected to be small. Hence these
positions are associated with slower dynamics. In contrast to the
cortical pushing model, for the cytoplasmic pulling model the longest
MTs that polymerize from posterior to anterior side of the cell induce
the largest centering forces. As a result this slower dynamics is not
observed in that model (see \pr{sfg:S3_pos_pull}).

So far, we have only considered the effect of perturbing the cell
shape on positioning of the PNC, while biophysical parameters were
kept fixed. Although we have chosen the parameters of the model based
on experimental measurements, the experimental uncertainties
associated with these parameters, including cytoplasmic viscosity, MT
polymerization rates, and motor protein densities, are quite
significant. Also, many of these parameters\emdash/for example the
density of the motors\emdash/are internally regulated by the cell in
different stages of the cell division. Thus, it is important to study
how the dynamics change over feasible ranges of the biophysical
parameters. To this end, we simulate PNC migration in $S_{1}$ after
changing the MT turnover time to $\tau_\lbl{cat}=8\usep\sec$ in
\pr{eq:tau_cat} for the cortical pushing model, and vary the motor
protein density as $n_\lbl{dyn}=0.01, 0.02, \text{and }
0.04\usep(\microm)^{-1}$ in the cytoplasmic pulling model.

\Pr{fig:tauvar} shows snapshots of the conformation of the MTs after
the PNC has reached its steady state position in the $S_{1}$ cell
geometry for $\tau_\lbl{cat}=4\usep\sec$ and $8\usep\sec$. It is clear
that the scale of the MT deformations, as well as the number of fibers
interacting with the periphery, are significantly larger for
simulations with larger turnover times.

\begin{figure}[!bt]
  \centering
  \subfigure[$\tau_\lbl{cat}=4\usep\sec$\label{sfg:S1_tauvar_1}]{\includegraphics[height=1.6in]{S1_push3}}
  \hspace{18pt}
  \subfigure[$\tau_\lbl{cat}=8\usep\sec$\label{sfg:S1_tauvar_2}]{\includegraphics[height=1.6in]{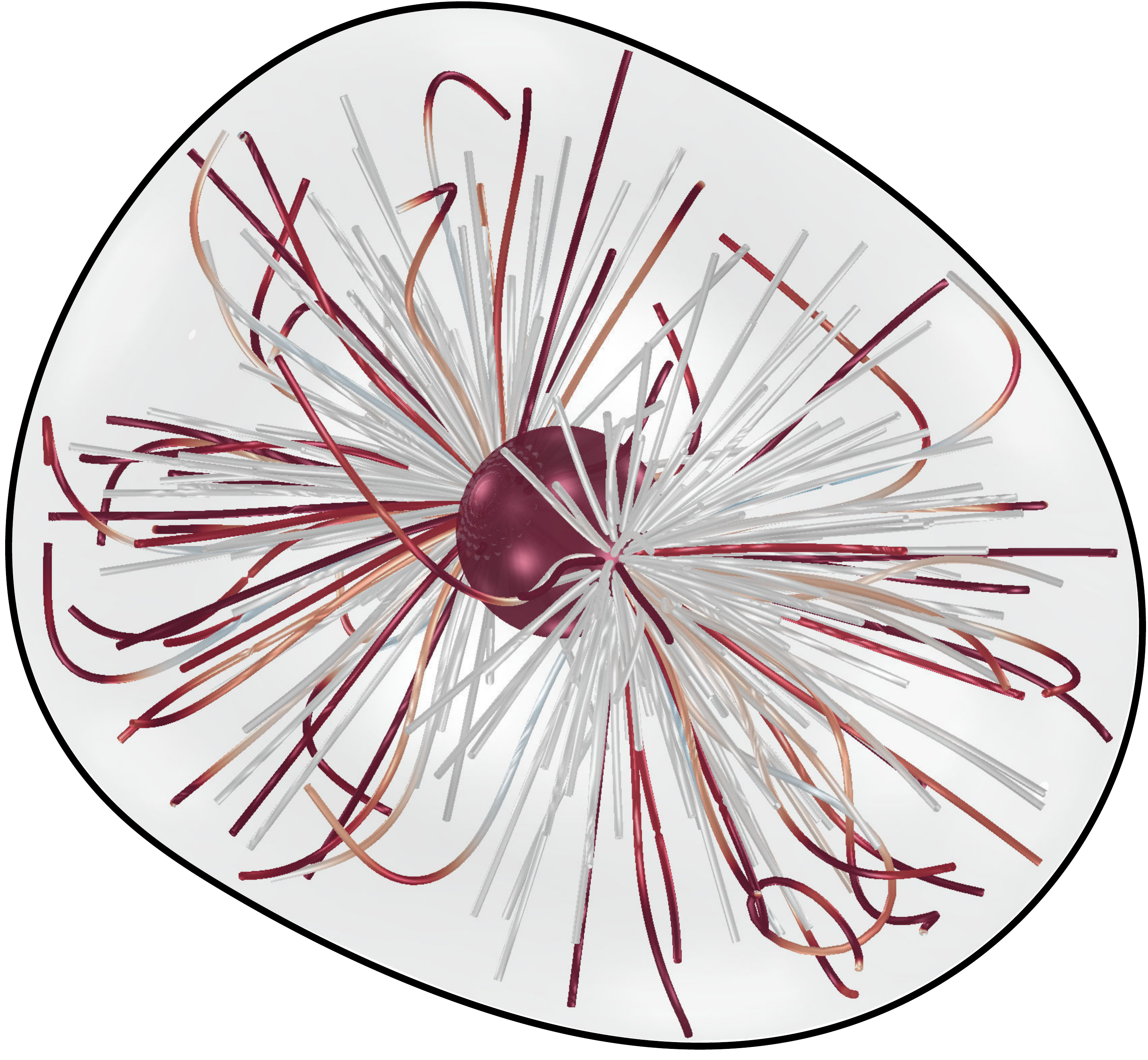}}\\
  \subfigure[Positioning and alignment\label{sfg:posangle_tauvar}]{
    \setlength\figurewidth{3.4in}%
    \setlength\figureheight{1.2in}%
    \ifUseTikz
    \tikzset{external/export=false}%for some reason externalization messes up the legend line colors
    \fi
    \includepgf{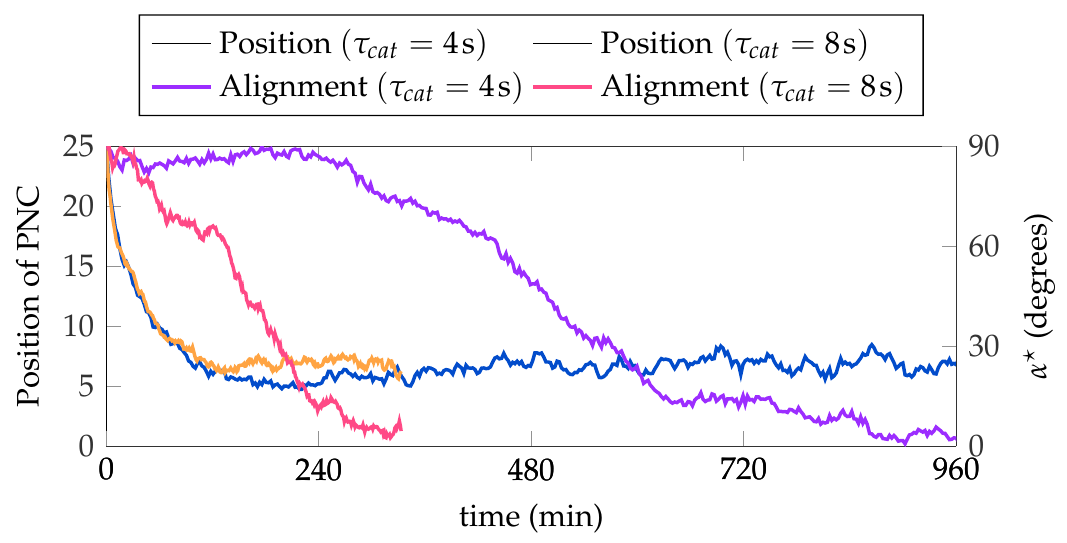}%
  }
  \mcaption{fig:tauvar}{Centering for different turnover
    time}{Snapshot of the conformation of the MTs after reaching the
    steady state position and alignment, using turnover times
    \subref{sfg:S1_tauvar_1} $\tau_\lbl{cat}=4\usep\sec$ (same as
    \pr{fig:S1}) and \subref{sfg:S1_tauvar_2}
    $\tau_\lbl{cat}=8\usep\sec$ respectively;
    \subref{sfg:posangle_tauvar} the position and alignment angle of
    the PNC with respect to time for $\tau_\lbl{cat}=4\usep\sec$ and
    $8\usep\sec$; see \pr{eq:tau_cat}.  }
\end{figure}

The centering and rotation dynamics of the PNC for these two turnover
times are compared in \pr{sfg:posangle_tauvar}. The main observation
is that while increasing $\tau_\lbl{cat}$ does not strongly affect the
migration dynamics of the PNC, it does strongly affect its rotation;
Rotation occurs approximately three times faster for
$\tau_\lbl{cat}=8\usep\sec$ than for $\tau_\lbl{cat}=4\usep\sec$. If
the only effect of increasing the turnover time was to increase the
average polymerization force from the periphery, we would have
expected the translational and rotational motion to be affected
similarly. However, this change plainly has a large effect, for
example, on the conformation of the centrosomal MTs, and hence on the
transmission of the polymerization force to the PNC.

Finally, we study the effect of motor protein density, $n_\lbl{dyn}$,
on PNC positioning dynamics for the cytoplasmic pulling model in the
$S_1$ geometry. For this case, \pr{fig:S_varpull} shows the evolution
of PNC position and alignment angle for the different values of
$n_\lbl{dyn}$. The simplest analysis suggests that the speed of PNC
migration is proportional to $\vector{F}^E$, which would in turn
suggest that the time to both translate and rotate to proper position
halves with each doubling of $n_\lbl{dyn}$. \Pr{sfg:pulling_osc}
clearly shows that the dynamics does not follow this scaling. For of
the lowest motor density, $n_\lbl{dyn}=0.01$, the PNC monotonically
migrates towards its steady-state position and thereafter shows small
fluctuations around that position. When the active force is increased
by $2$ and $4$ times, the time it takes the PNC to reach its
time-averaged steady-state position is reduced by a factor of 3 and 9,
respectively. This suggests that the speed of translation has a
stronger-than-linear scaling with the active force.

Another interesting effect of increasing motor density is that the PNC
starts to oscillate around its mechanical equilibrium position with a
well-defined amplitude and frequency.  The amplitude of oscillation is
greatly increased from $n_\lbl{dyn}=0.02$ to $n_\lbl{dyn}=0.04$, while
the frequency of these oscillations remain roughly unchanged.  The
source of this rocking motion has similarly been related to enhanced
activity of molecular motors in that particular stage with the
difference that in this case the motors are bound to the cell cortex
\citep{Pecreaux2006a}.  We have formulated a simple mathematical model
based on the relative magnitudes of growth and shrinking velocities
compared with the velocity of PNC migration to explain these
oscillations (to be published elsewhere).
\begin{figure}[!tb]
  \centering
  \subfigure[\label{sfg:pulling_osc}]{
    \hspace{-14pt}
    \setlength\figurewidth{2.7in}%
    \setlength\figureheight{1.12in}%
    \includepgf{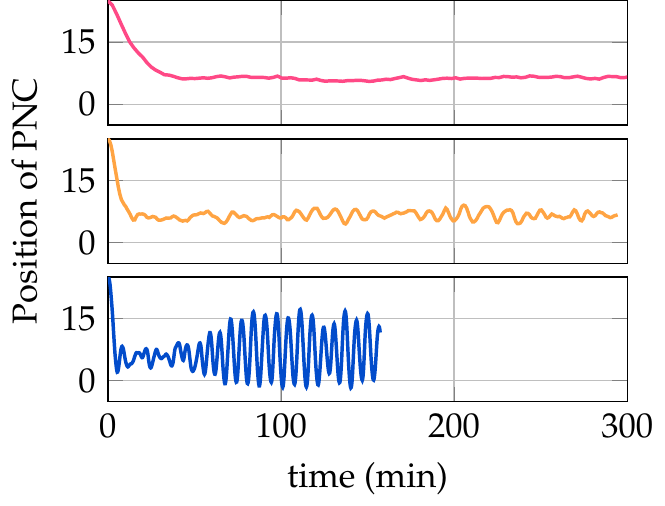}%
  }
  \subfigure[\label{sfg:pulling_osc_angle}]{
    \hspace{-14pt}
    \setlength\figurewidth{2.4in}%
    \setlength\figureheight{1.6in}%
    \includepgf{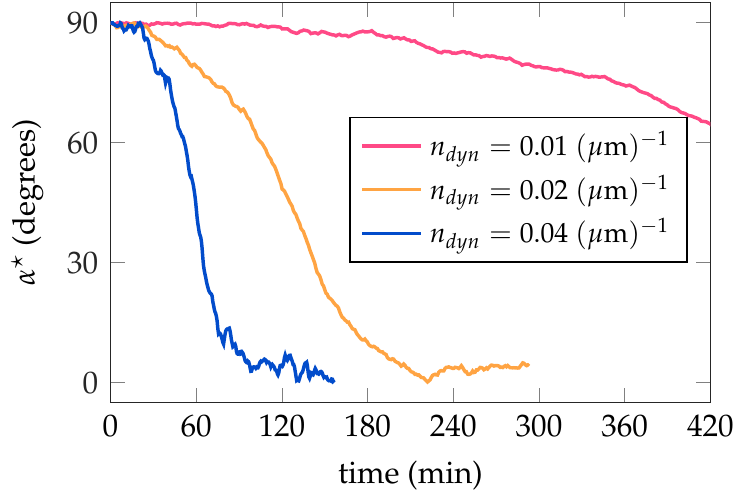}%
  }\\
  % \subfigure[\label{sfg:sche_osc}]{
  % \includegraphics[height=1.350in,draft=false]{sche_osc}
  % \def\figscale{.45}
  % \includepgf{schematic-osc}
  % }
  %
  \mcaption{fig:S_varpull}{Position and angle of the PNC for different
    dynein densities}{\subref{sfg:pulling_osc} The $\zhat$ component
    of the position of the PNC for three different densities of dynein
    motors, $n_{dyn}=0.01, \, 0.02, \text{and }
    0.04\usep(\microm)^{-1}$ \subref{sfg:pulling_osc_angle} the
    alignment angle of the MTOC axis with respect to the $\zhat$ axis
    with respect to time for different density of the motors,
    $n_\lbl{dyn}=0.01, 0.02, \text{and } 0.04\usep(\microm)^{-1}$.
    %\subref{sfg:sche_osc} a schematic demonstration of the interplay
    %between migration rate and the growth rate that can lead to
    %oscillations in the PNC position.
  }
\end{figure}

The angle between the axis of MTOCs and $\zhat$ axis at different
values of $n_\lbl{dyn}$ are plotted in \pr{sfg:pulling_osc_angle}.  As
is seen, the alignment time is greatly reduced when $n_\lbl{dyn}$ is
increased from $0.01$ to $0.02$ and continues to reduce by roughly a
factor of 2 for $n_\lbl{dyn}=0.04$.  This clearly shows that the
rotation time is also not a linear function of the total active force.

In summary we used our numerical platform to study a critical question
in the cell division: the effect of confinement geometry and
biophysical parameters upon the dynamics of pronuclear migration.  We
used three different geometries, and two models for the force
transduction on MTs required for the positioning of the PNC, i.e., a
cortical pushing and the cytoplasmic pulling models.
%% We investigated the effect of geometry by considering three different
%% geometries, and two different force generation models (\emph{cortical
%%   pushing} and \emph{cytoplasmic pulling}) for the centering of the
%% PNC.
Our results show that the time-scales of positioning and rotation of
the PNC can change by several factors depending upon the choice of
model. We showed also that varying the biophysical parameters in our
force transduction models can change the dynamics tremendously.
Finally the steady state orientations and positions of the PNC are in
general agreement with the experimental results of \cite{Minc2011},
for geometries that are qualitatively similar; direct comparison
between these two results is not possible since the geometries are not
identical. Also the variety of geometries considered in
\cite{Minc2011} is much large than what we have studied here. A more
systematic study of the cell size and shape anisotropy (as is done in
\cite{Minc2011}) would be an interesting extension of the current
study.

\subsection{Sedimentation of a cloud of flexible
fibers\label{ssc:sediment}} Due its importance in industrial settings
and fundamental studies, the sedimentation of particles in fluids has
been a key area of study in suspension mechanics over the past half
century.  The classical experimental and theoretical literature on
sedimentation of spherical particles in viscous fluids is reviewed by
\cite{Davis1985} and more recently by \cite{Guazzelli2010}.  There has
also been several experimental \citep{Herzhaft1996, Herzhaft1999,
  Metzger2005} and numerical \citep{Butler2002,SDS2005, Tornberg2006}
studies on sedimentation of suspensions of rigid fibers.  Despite its
long history and sustained interest, several aspects of sedimenting
suspensions remain poorly understood.  The origin of many of the
complex behaviors has been attributed to the coupling between the long
range many-body hydrodynamic interactions and the relative arrangement
of the particles, which eventually determines the collective behavior
of the suspension \citep{Ramaswamy2001}. For example, sedimenting
fiber suspensions show formation of inhomogeneous clusters of fibers
and subsequent enhancement of the sedimentation rate in both numerical
and experimental studies \citep{Metzger2005, SDS2005, Tornberg2006}.

An interesting problem in this context is the sedimentation of a
cloud\footnote{the volume of the fluid wherein the particles are
  dispersed} of particles in a viscous fluid. \cite{Adachi1977}
experimentally studied the evolution of a spherically-shaped cloud of
a spherical particles.  They showed that the cloud undergoes complex
deformation as it sediments in the fluid. In a repeating cycle, the
initial spherical shape of the cloud deforms to a torus-like shape.
The torus structure eventually breaks into smaller clouds which
themselves evolve into tori.  An extensive overview of this problem
and the relevant literature are given in \citep{Machu2001,
  Metzger2007}.  The majority of research involves studying clouds
formed by spherical particles and the effect of shape anisotropy has
been studied to a much lesser extent \citep{Metzger2005, Park2010}.
To the best of our knowledge, sedimentation of flexible fibers has not
been studied experimentally and simulation studies are limited to
sedimentation of a single semi-flexible fiber \citep{LMSS2013}.  We
also note that, recently \cite{Manikantan2016} studied the
sedimentation of weakly flexible fiber suspension using particle
simulations. In this study the fibers are treated as rigid rods to the
leading order and the effect of flexibility only enters through the
rate of rotation of the fibers using a continuum theory proposed by
the same authors \citep{Manikantan2014}.

conclusive study of different aspects of sedimentation of a cloud of
flexible fibers is beyond of the scope of this work.  Instead, our aim
here is to demonstrate that the present numerical framework can be
used to investigate some of these aspects by showcasing the dynamics
in a limited range of parameter space.  We will show that even this
limited set of results poses a number of interesting physical
questions to be further pursued.

\subsubsection{Background}
For simplicity and for consistency with the previous studies, we
nondimensionalize velocity and length on the initial velocity and
radius of the cloud, respectively.  Taking the cloud to be a spherical
droplet with effective excess weight of $N_F F$ and radius $R_0$, the
initial velocity of a cloud can be estimated as $V_{c}={N_F
  F}/{(5\pi \mu R_{0})}$ where $N_F$ is again the number of fibers, $F$
is the net force acting on each fiber due to their density difference
with the fluid.

In dilute suspensions, the flow induced by the sedimenting objects in
a Stokesian fluid can be approximated (to the first degree in volume
fraction) by representing each object as a point-force within the
fluid where the force is the total gravitational force acting on the
fiber. Experimental studies on clouds of spherical particles have
mostly been limited to the dilute regime ($\phi \le 0.20$). As a
result simulations treating spheres as point-particles have been
successful in reproducing many of the observed experimental behavior
in these systems \citep{Metzger2007}.  In contrast, for a cloud of
fibers the fiber lengths can be comparable to the radius of the cloud
itself and the average distance between the fibers may very well be
much smaller than their average length, resulting in large effective
volume fractions.  Thus, reducing the interactions of each fiber with
the fluid to a point-force becomes exceedingly inaccurate, and the
geometry of individual fibers and their HIs need to be explicitly
included. Such differences in geometry and parameter regime, as well
as the flexibility of the fibers is likely to result in differences in
the evolution of the shape and velocity of a cloud of sedimenting
fibers, with respect to clouds of spherical particles.

\subsubsection{Simulation setup}
We assume that all fibers have the same length and weight.  The
boundary conditions for the two ends of fibers are \emph{free}, and
given by relations in \pr{eqn:bc-free}.  The simulations are
initialized by randomly distributing, within an unbounded fluid, the
center-of-masses of $N_F=1024$ fibers in a spherical volume of radius
$R_{0}=1$. Fiber orientations are also taken as random.

For $t>0$ a downwards gravitational load, $\f^E=f_0 \zhat$, is applied
to all fibers. The simulations are performed for fibers of lengths
$L\in \{0.15 R_{0}, 0.25 R_0, 0.32 R_{0}, 0.40 R_{0}\}$ and aspect
ratio of $\epsilon=0.01$.

We define the effective volume fraction of the fibers as
\begin{equation}
  \phi^\star= \frac{N(L/2)^{3}}{(4/3) \pi R_{0}^{3}}.
\end{equation}
The lengths $L = 0.15 R_{0}, 0.25 R_0, 0.32 R_{0}, \text{and } 0.40
R_{0}$ correspond to $\phi^\star \approx 0.10, 0.50, 1.00, \text{and
}2.00$. The flexural rigidity of the fibers is chosen to be
$E=3.2\times 10^{-3} f_{0} R_{0}^{3}$, which results in very small
deformation of fibers in the shortest studied length and moderate
deformation for the longest.

\subsubsection{Results}
\Pr{fig:sed_1} shows snapshots of the evolution of the fiber cloud as
it sediments through the fluid.
\begin{figure}[!bt]
  \centering
  \raisebox{38pt}{\begin{minipage}[b]{0.4\textwidth}
    \centering
    \subfigure[\label{sfg:sed-init}]{
      \includegraphics[width=\textwidth]{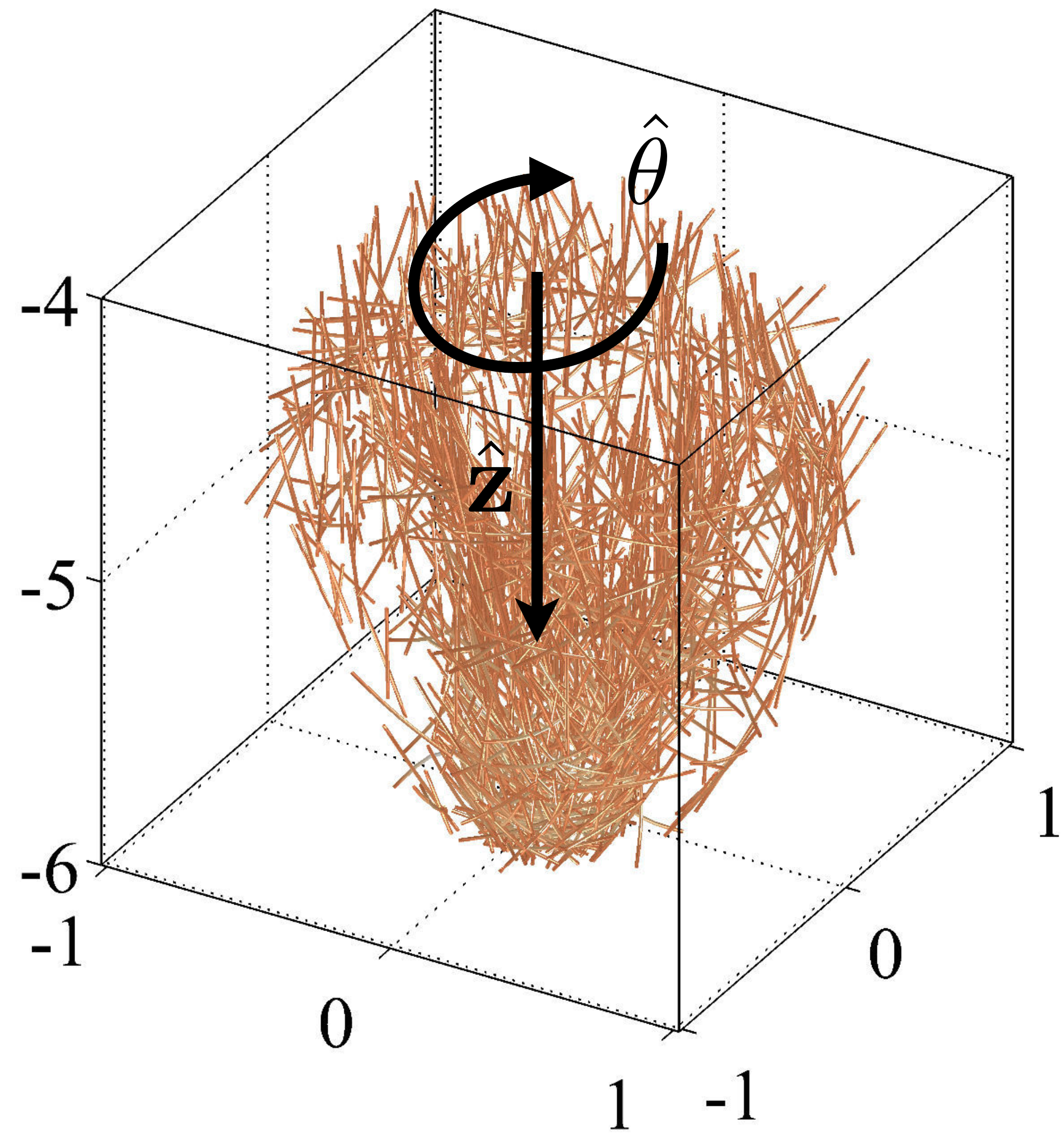}}\\
    \addtocounter{subfigure}{1}%adjusting the numbers to go horizontal
    \subfigure[\label{sfg:sed-pinch}]{
      \includegraphics[width=\textwidth]{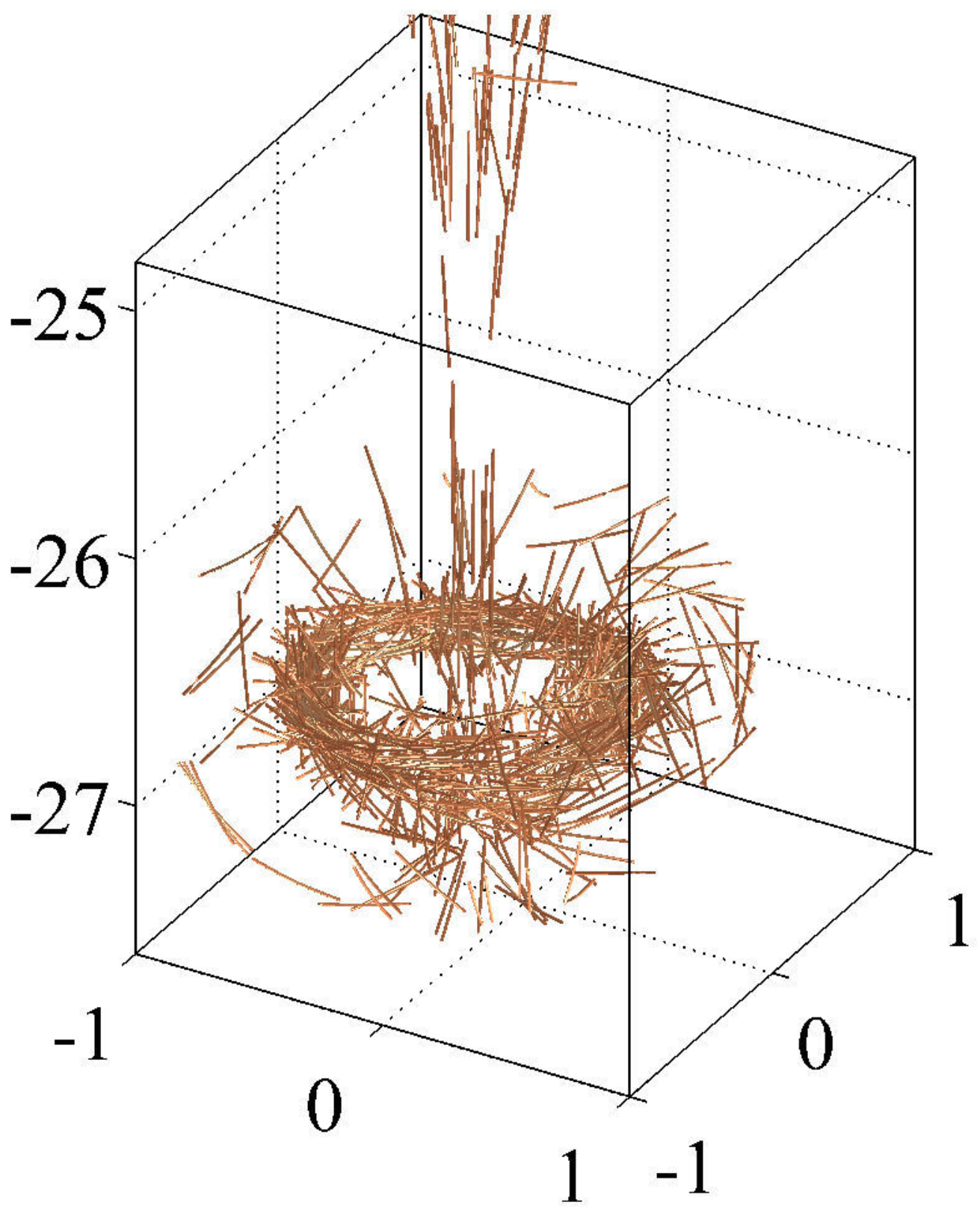}}
  \end{minipage}}
  \hspace{0.35in}
  \begin{minipage}[b]{0.4\textwidth}
    \centering
    \addtocounter{subfigure}{-2}%adjusting the numbers to go horizontal
    \subfigure[\label{sfg:sed-torus}]{
      \includegraphics[width=\textwidth]{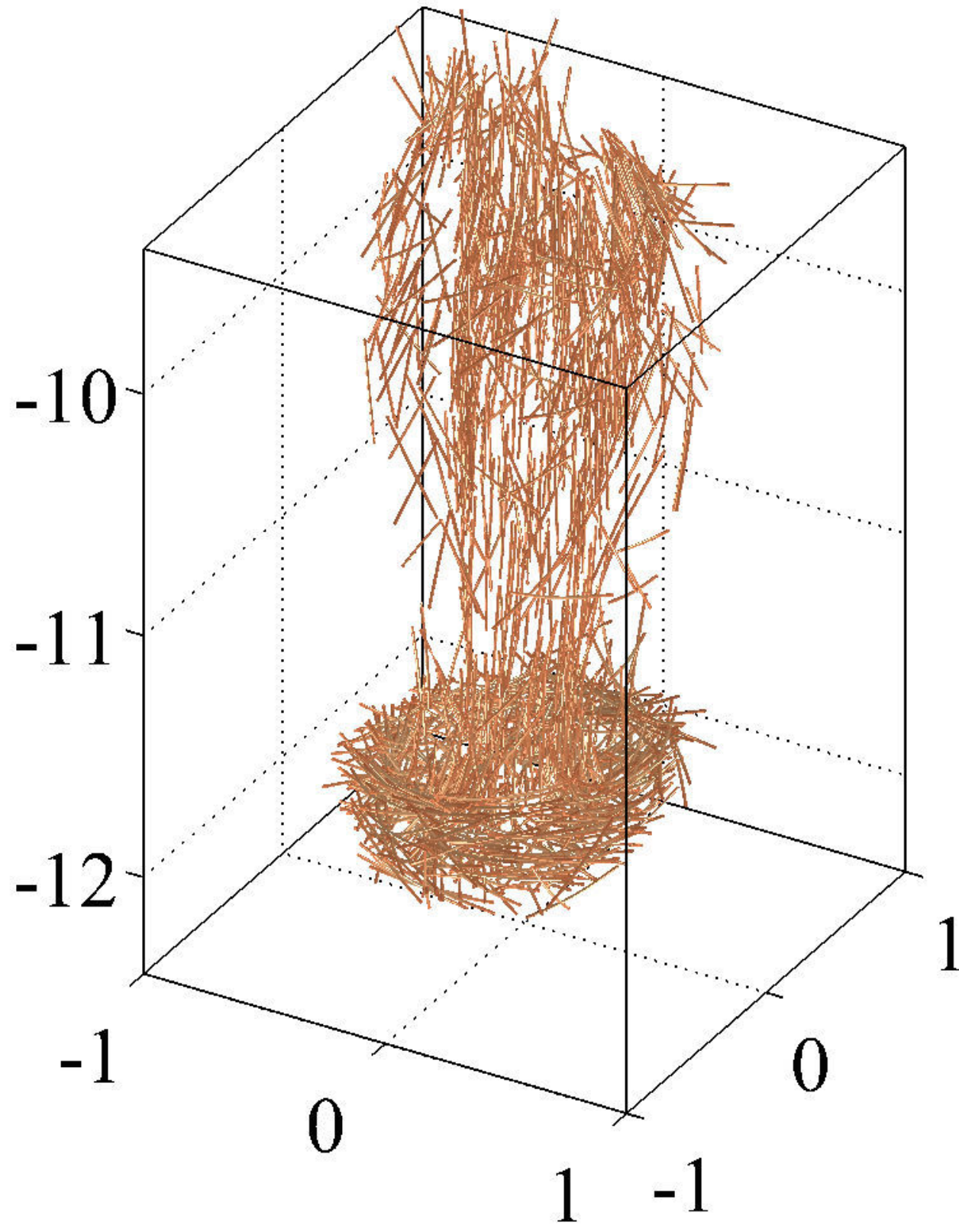}}\\
    \addtocounter{subfigure}{1}%adjusting the numbers to go horizontal
    \subfigure[\label{sfg:sed-bottom}]{
      \includegraphics[width=\textwidth]{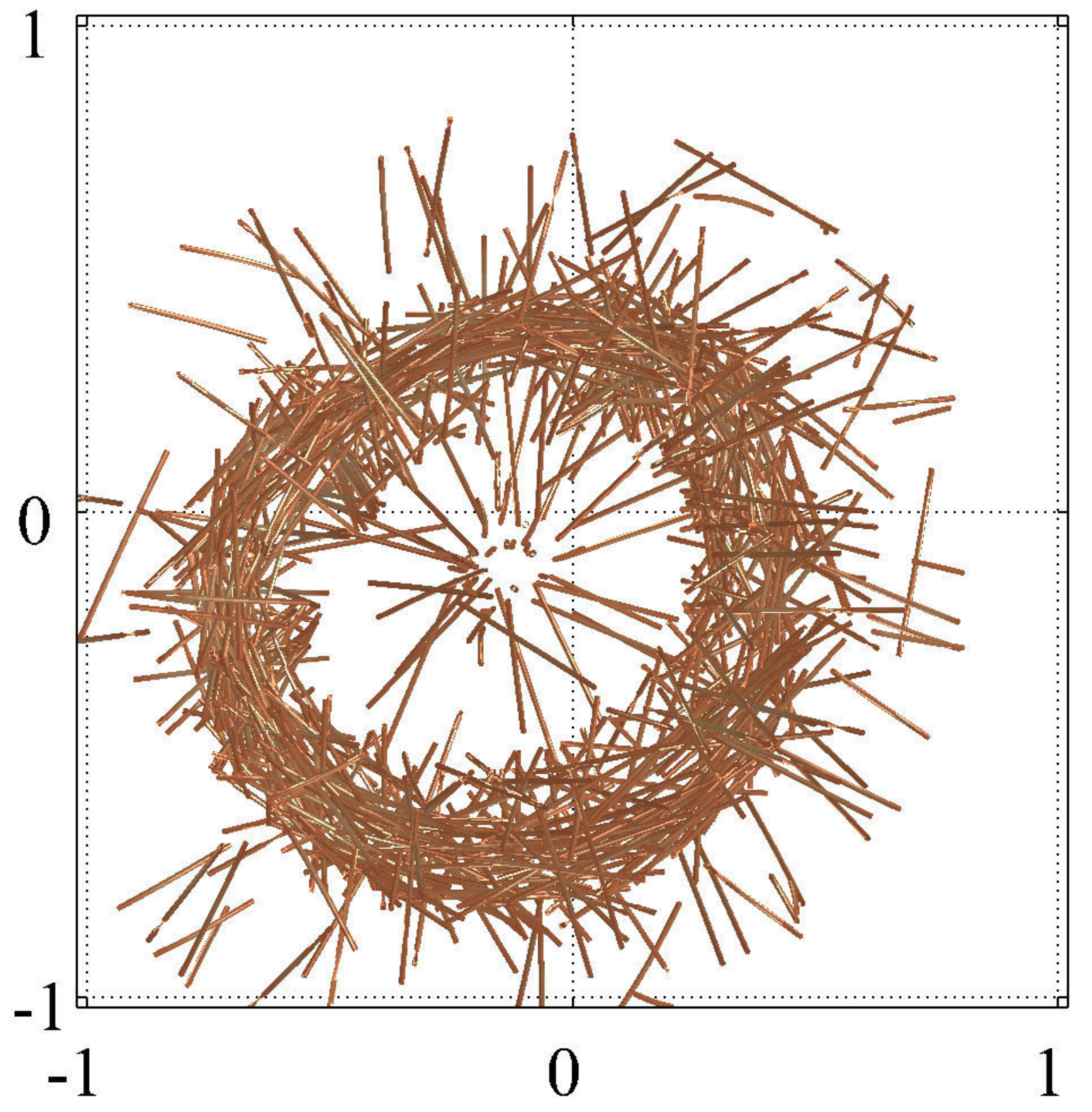}}
  \end{minipage}
  \mcaption{fig:sed_1}{key Stages in the sedimentation of fibers}{
    Simulation snapshots of the key stages in the shape evolution of
    the sedimenting cloud of $1024$ flexible fibers with
    $L/R_{0}=0.32$.  \subref{sfg:sed-init} The core of the cloud
    increases speed and the front of the cloud stretches while the
    rear contracts.  \subref{sfg:sed-torus} A torus-like structure is
    formed in the front part of the cloud, followed by a trail of
    fibers in the rear.  \subref{sfg:sed-pinch} The torus detaches
    from the trail.  \subref{sfg:sed-bottom} the projection of
    subfigure~\subref{sfg:sed-pinch} in $r{-}\theta$ plane to better
    depict the organization of the fibers inside the torus ring.  }
\end{figure}
Certain qualitative features are observed for all fiber lengths.  In
the initial stage of the sedimentation, the fibers align with the
$\zhat$ direction and an inward cusp is observed in the rear face of
the cloud, which is shown in \pr{sfg:sed-init}. Here the rear and
front are defined as the parts of the cloud volume with minimum and
maximum values of $z$, respectively.

The front face stretching continues until a torus of fibers is formed,
trailed by fibers aligned in the sedimentation direction; See
\pr{sfg:sed-torus}.  The torus diameter then grows in time while the
tail of downward aligned fibers gets longer and thinner and eventually
pinches off from the torus as is seen in \pr{sfg:sed-pinch}.  The
torus remains stable for up to $150 R_{0}$ in sedimentation distance.

These observations are in general agreement with experiments and
simulation results of \cite{Metzger2007} and \cite{Park2010}, who
studied clouds of spherical particles and rigid fibers.  These studies
also found that after sedimenting more than $600 R_{0}$ the torus
structure breaks up into smaller tori.  We do not see this in our
simulations, most likely due to the shorter simulation times. We find
that at sedimentation distances larger than $150 R_0$, the fibers in
the torus become exceedingly close and in near contact and the
simulations became numerically unstable.  Reducing the time step did
not resolve this issue.  A similar difficulty was reported by
\cite{Park2010} in their simulations of sedimentation of a cloud of
rigid fibers.  These authors removed this instability by not
accounting for HIs between the fibers closer than a cutoff distance.
As a result the fibers could pass one another without generating
numerical difficulties, though at the unknown cost of not properly
accounting for contact mechanics.  We did not pursue that route here,
and are currently working on including a model for steric interactions
between fibers. Another interesting observation is the azimuthal
alignment of the fibers within the torus.  This is seen in the bottom
view of the cloud in \pr{sfg:sed-bottom}.

\begin{figure}[!tb]
  \centering
  \setlength\figurewidth{3.5in}
  \setlength\figureheight{1.6in}
  %\fbox{%
  \hspace*{-6pt}\includepgf{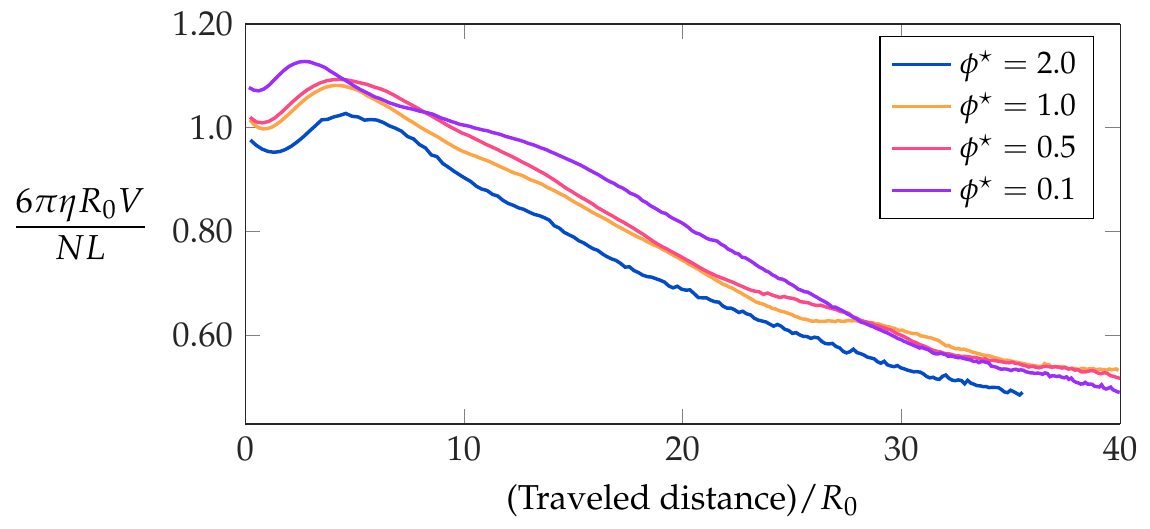}%
  %}
  \captionsetup{width=.8\linewidth}
  \mcaption{fig:sed_plot}{Cloud Sedimentation Velocity}{ The
    dimensionless velocity of the cloud versus the sedimentation
    distance at different volume fractions.}
\end{figure}

The dimensionless velocity of the cloud as a function of the
dimensionless sedimentation distance, for different volume fractions,
is given in \pr{fig:sed_plot}. The results for different volume
fractions all nearly collapse onto a single curve. At all volume
fractions, the velocity starts to increase after an initial decrease
and then reaches a maximum.  The maximum roughly corresponds to the
configuration shown in \pr{sfg:sed-init}.  Thereafter, the formation
of the tail and the leakage of fibers from the front of the cloud into
the tail results in a monotonic decrease of the velocity. This
decrease in sedimentation velocity levels off after the torus detaches
from the tail in $z/R_{0} > 35$, where $z$ is the net displacement of
the center of the mass of the cloud since the start of
sedimentation.

A key point is that the dimensionless velocity is primarily determined
by the cloud's traveled distance to the initial radius of the cloud, $z/R_0$, and is not explicitly a
function of time. We also observe that, similar to the sedimentation
velocity, the evolution of the shape of the cloud is primarily set by
the distance traveled (results not shown here), irrespective of the
fibers length and their effective volume fraction.

Next, to gain further insight at a more coarse-grained level, we
compute the local number density and the velocity of the fibers.
Considering that the structure is statistically symmetric around the
direction of gravitational force, $\zhat$, we use a cylindrical
coordinate and average the positions and velocities of the fibers in
$\hat{\vector{\theta}}$ direction.  We choose center of the coordinate
as $(x_c,y_c,z_c) = (0,0,(\max(z)+\min(z))/2)$.  Note that the exact
value of $z_c$ has no effect on the results presented
here. \Pr{fig:sed_den} illustrates the averaged quantities of interest
at the same times shown in
\pr{fig:sed_1}\subref{sfg:sed-init}--\subref{sfg:sed-pinch} in the
$r{-}z$ plane.  The number density is color-coded with white (light)
denoting zero and red (dark) denoting the largest number density. The
velocity vectors at any given point represent the velocity of the
fibers relative to the average sedimentation velocity of the cloud's
center of mass.
\begin{figure}[!htb]
  \centering
  \subfigure[\label{sfg:sed_den-a}]{\includegraphics[height=2.1in]{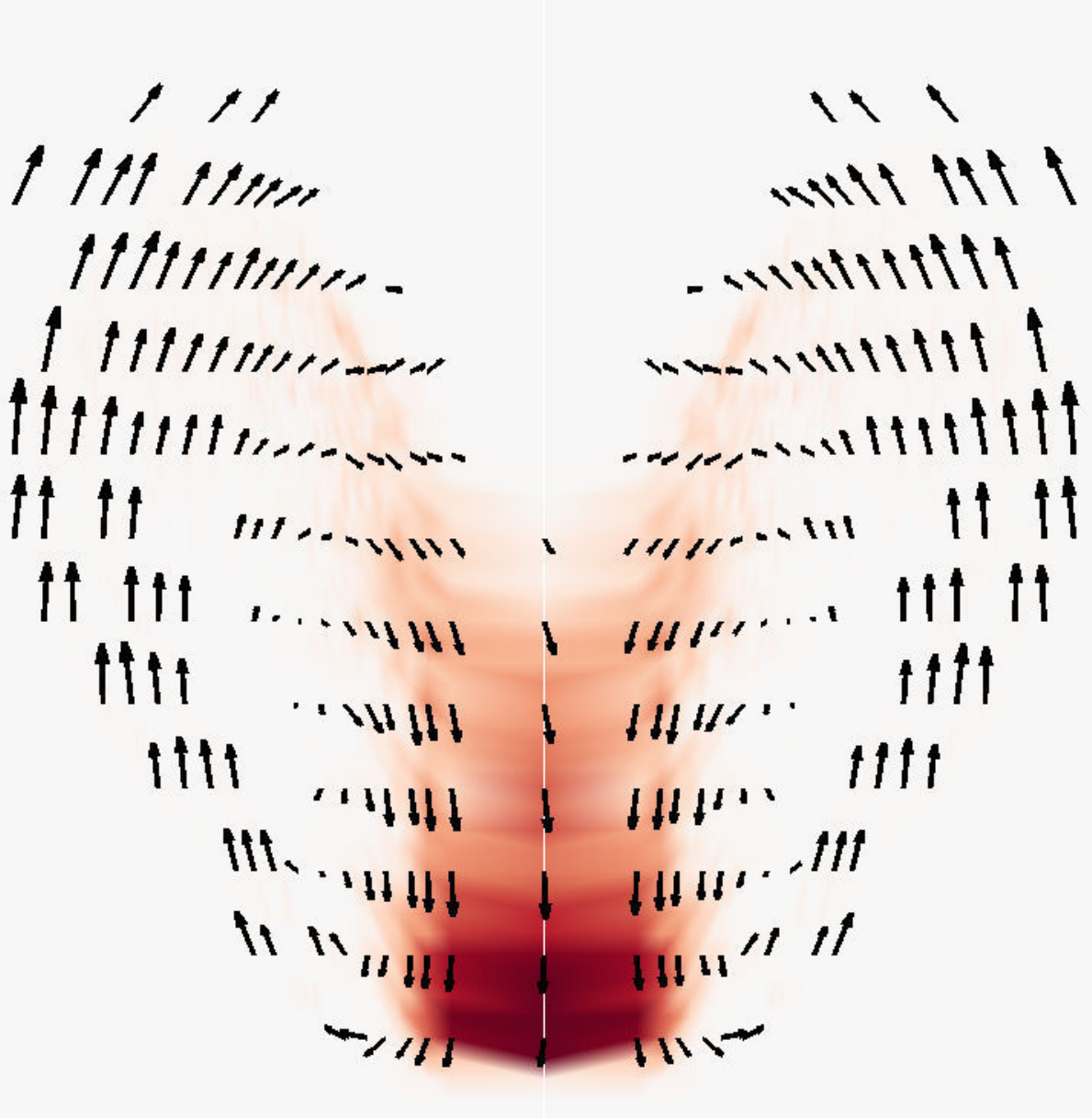}}
  \subfigure[\label{sfg:sed_den-b}]{\includegraphics[height=2.1in]{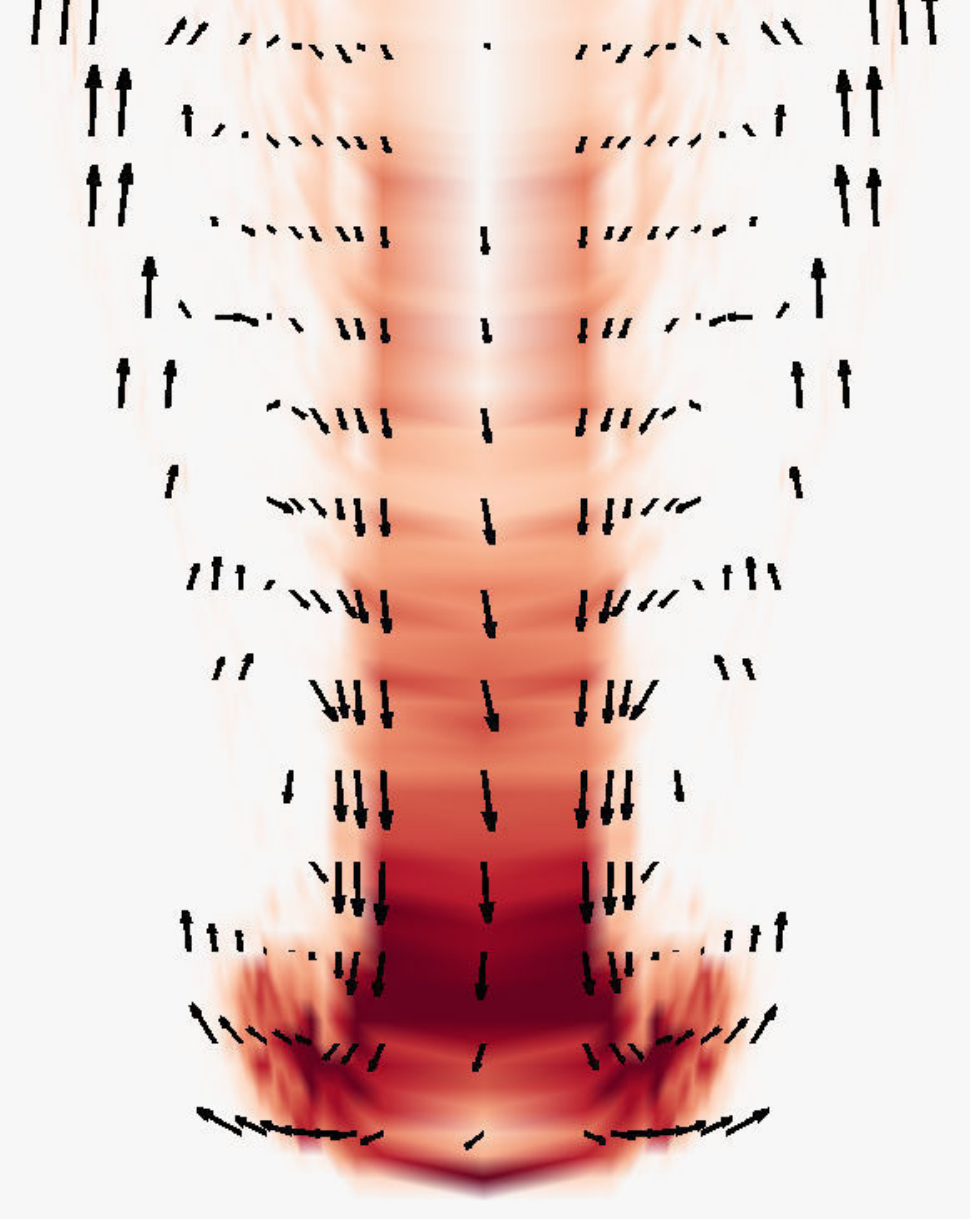}}
  \subfigure[\label{sfg:sed_den-c}]{\includegraphics[height=2.1in]{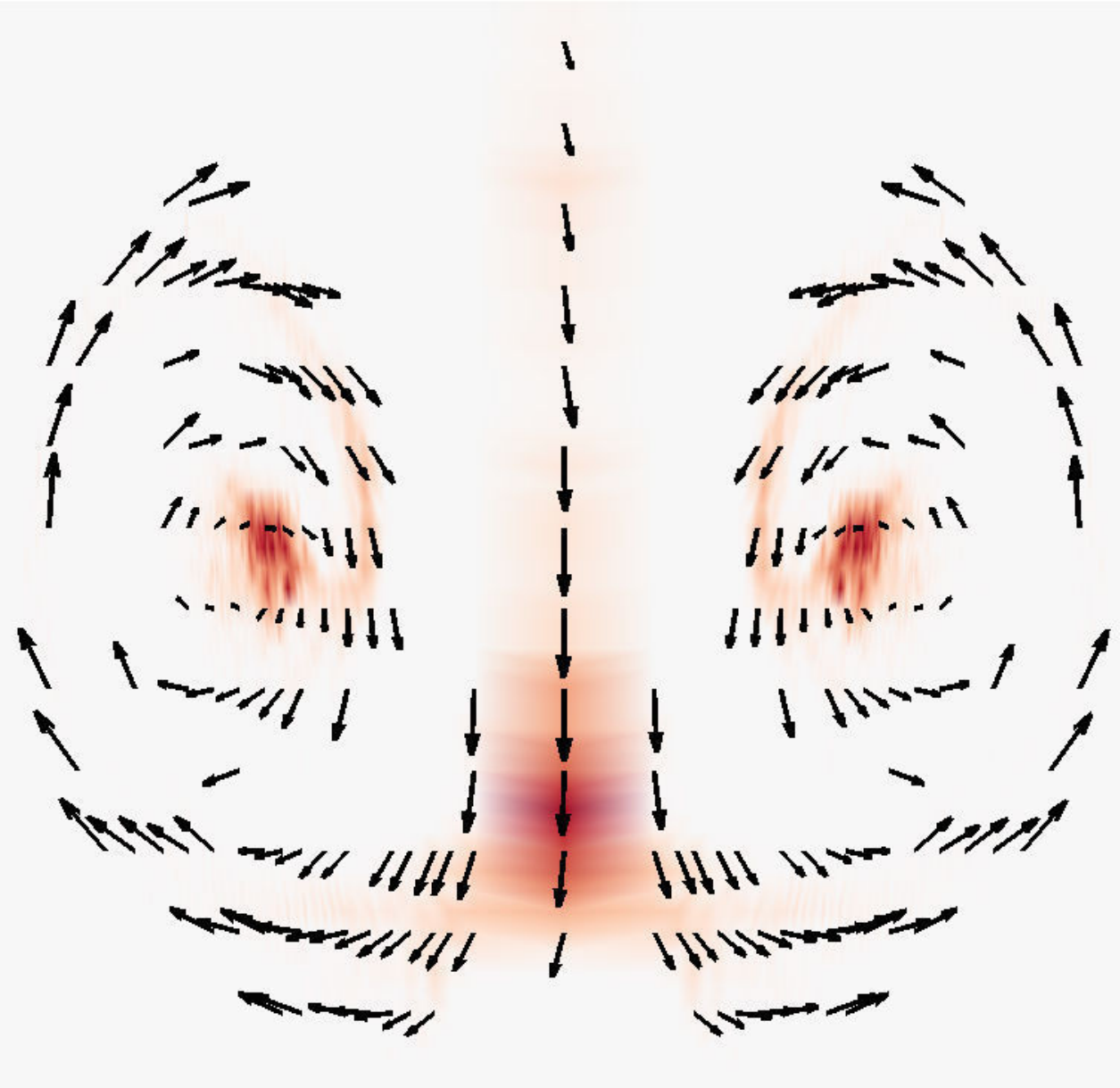}}
  \mcaption{fig:sed_den}{$\theta$-averaged number density and
    sedimentation velocity}{The number density of the fibers and their
    relative velocity with respect to the velocity of the center of
    the cloud, computed by sampling the configuration of the fibers
    from simulations.  The time and physical parameters of the
    simulations coincide with the ones illustrated in
    \pr{fig:sed_1}\subref{sfg:sed-init}--\subref{sfg:sed-pinch}.
    White (light) and red (dark) colors denote low and high volume
    fractions respectively.  }
\end{figure}

The initial form of the velocity field in \pr{fig:sed_den} resembles
the toroidal flows inside of a sedimenting droplet, which is in line
with the previous theoretical work that treats the cloud as a viscous
droplet \citep{Davis1985}. The flow structure and the concentration
variations become more complex in later stages. The velocity field in
\pr{sfg:sed_den-c} is reminiscent of the classical flow of a vortex
ring at high \emph{Reynolds} numbers.  In fact, \cite{Wen1996}
reported a similar vortex ring structure for sedimentation of
particles in inertially dominant flows, i.e. $\reynolds \gg 1$, where
the evolving shape of the cloud is determined by transition from
laminar to turbulent flow \citep{Wen1996, Bush2003}.

Along these lines, \cite{Tong1999} theoretically demonstrated that the
equations of motion describing the concentration and velocity
fluctuations in dilute suspensions of sedimenting non-Brownian
particles are analogous to the equations describing turbulent
convection at high Prandtl numbers. A more in-depth understanding of
this analogy involves studying the dynamics of the sedimenting clouds
in the framework of non-equilibrium statistical physics.  This
typically includes analyzing the correlation functions of number
density and velocity as well as their fluctuations in time and space.
For the relevant works discussing these aspects, the interested reader
is referred to \citep{Tong1999, Serge2001, Ramaswamy2001,
  Guazzelli2001, Nguyen2005}.  %
%can be composed of two primary velocity fields:
%\begin{inparaenum}[(i)]
%  \item the toroidal flow as in Fig. \pr{sfg:sed_den-a} and
%(ii) the velocity field directed towards the
% stagnation regions (sinks) that are located along the interior ring of the torus.
% This results in formation of a high concentration of aligned fibers along the ring as shown in
%figure \ref{fig:sed_1(d)}.
%\end{inparaenum}
%The

Finally, for a better understanding of the formation of the
azimuthally aligned fibers within the torus, we sample and monitor $q
= |\inner{\X_{s}}{\hat{\vector{\theta}}}|$ in space and time\footnote
{Note that in this definition, $\X_s$ varies along the length of the
  fibers}. The definition of $q$ gives zero for the fibers aligned in
$r{-}z$ plane and unity for fibers aligned in $\hat{\vector{\theta}}$
direction. The results are shown in \pr{fig:sed_nem}. Again, the times
coincide with those of \pr{fig:sed_1}.
\begin{figure}[!hbt]
  \centering
  \subfigure[]{\includegraphics[height=2.0in]{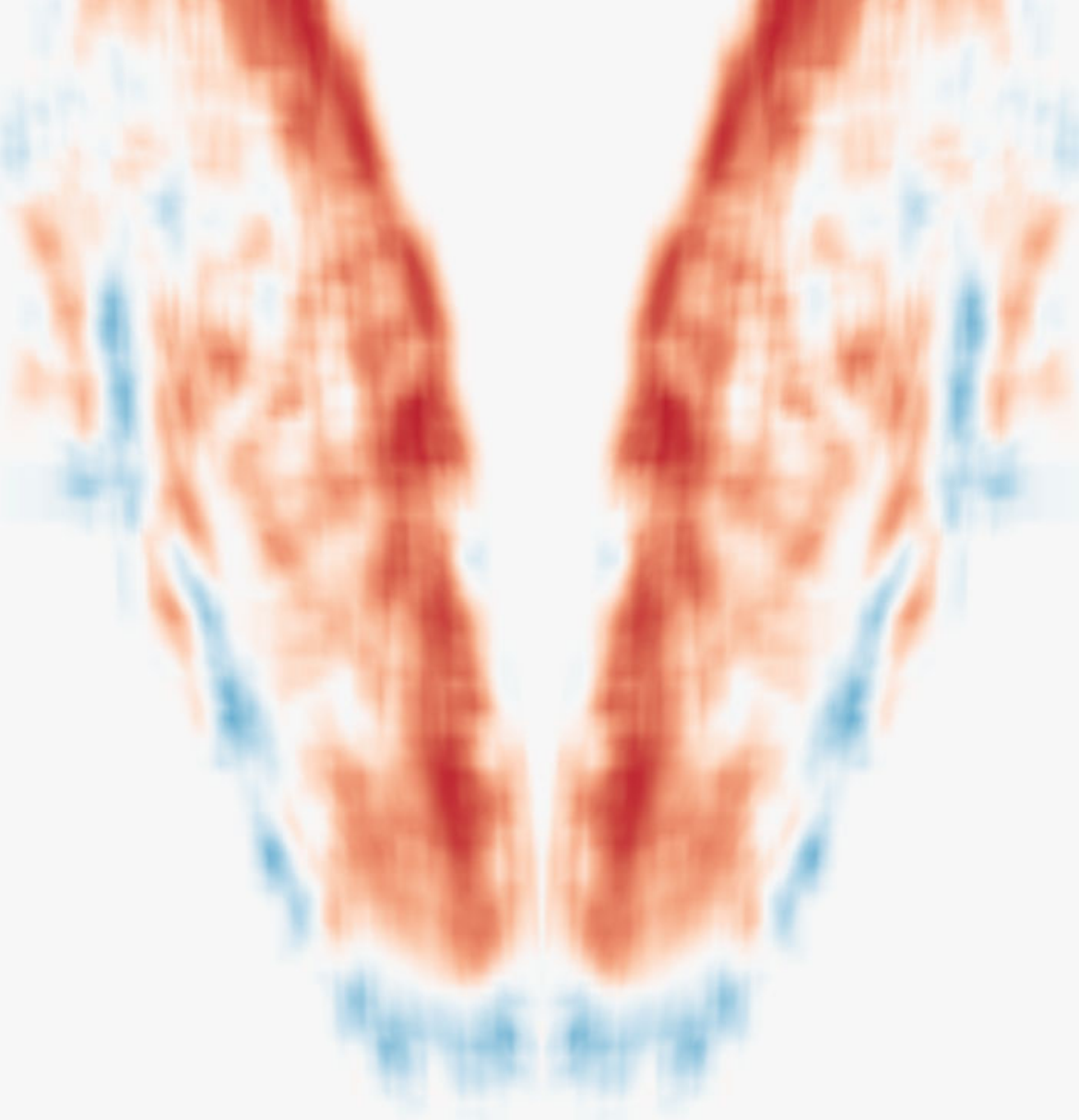}}
  \subfigure[]{\includegraphics[height=2.0in]{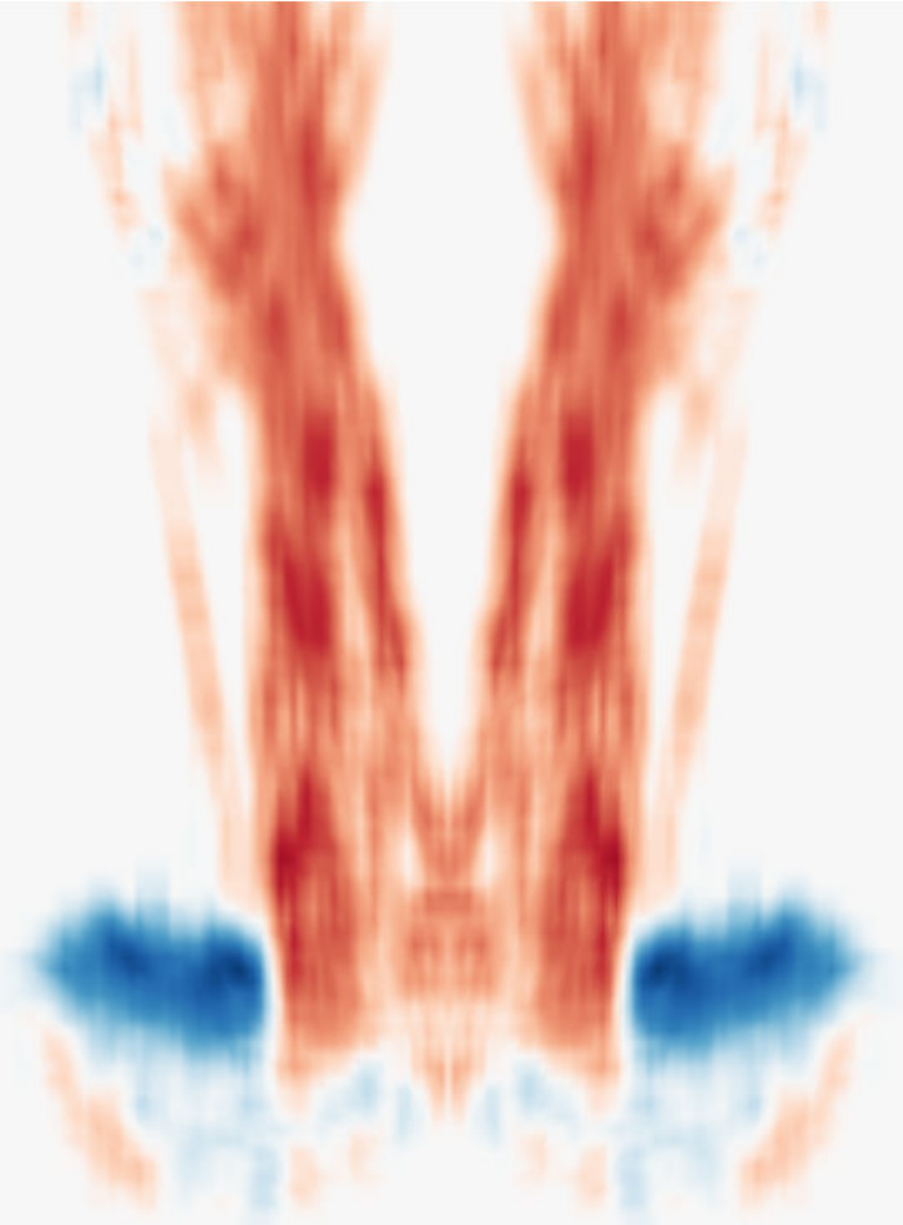}}
  \subfigure[]{\includegraphics[height=2.0in]{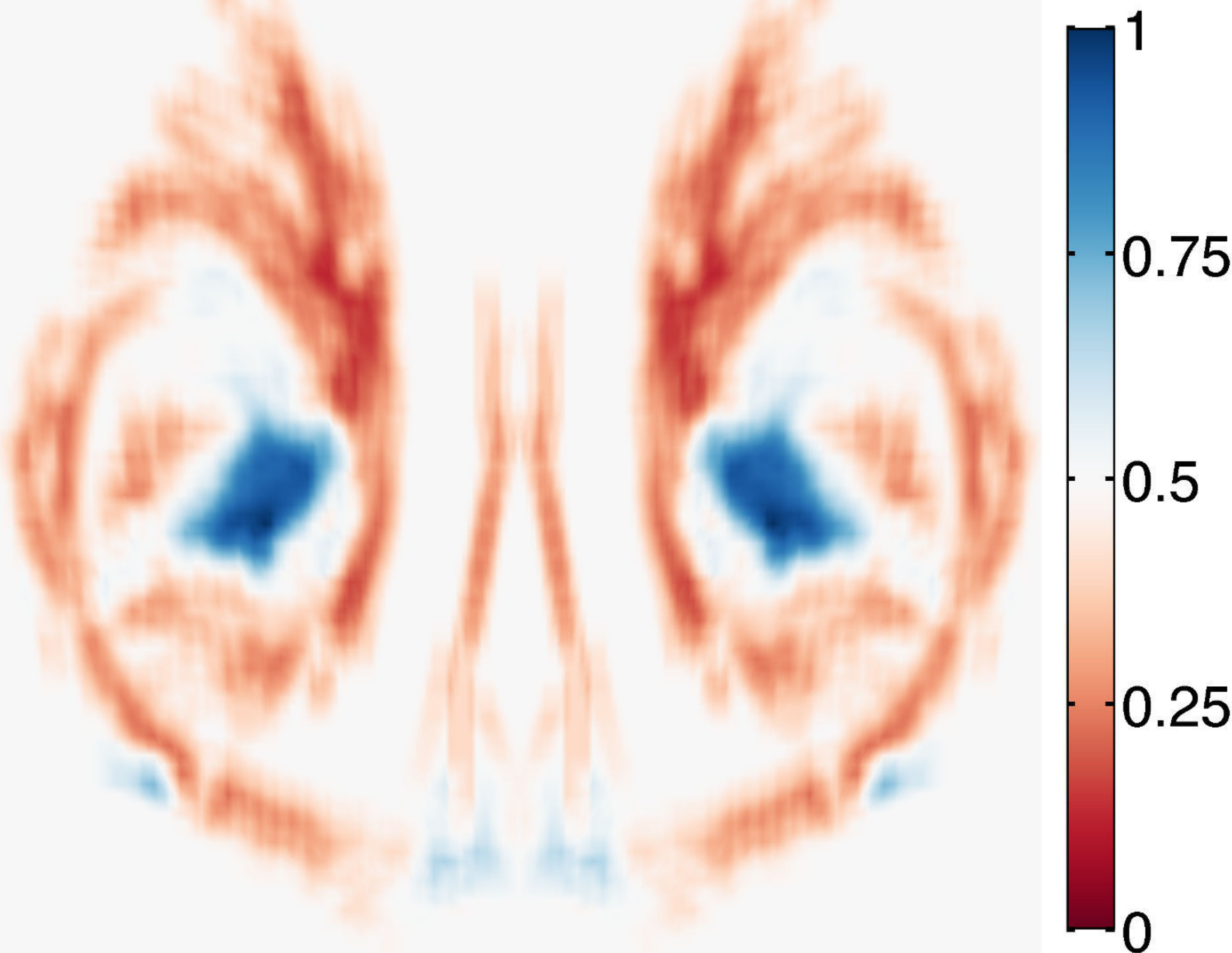}}
  \mcaption{fig:sed_nem}{Alignment of sedimenting fibers}{Variations
    of $q=|\inner{\X_{s}}{\hat{\vector{\theta}}}|$ in space at the
    time frames depicted in \pr{fig:sed_1} and \pr{fig:sed_den}.
    Recall that $\theta$ is the direction of the symmetry of the shape
    of the cloud and $\X_{s}$ is the tangent vector of the fiber at
    the particular sampling point.  The colors change from red
    denoting the alignment of the tangent (orientation) vector of the
    fiber in $r{-}z$ plane to blue denoting the alignment in
    $\hat{\vector{\theta}}$ direction.  }
\end{figure}

As can be seen, at early times the fibers are mostly aligned in
$r{-}z$ plane (red color), while a concentration of
$\hat{\vector{\theta}}$-aligned fibers are formed in the front region
of the cloud (blue color). At later times, these
$\hat{\vector{\theta}}$-aligned fibers are migrated to the ring-like
structure formed within the torus (see the supplementary material).

%We conclude this section with a few final remarks.  Using our
%numerical technique, we were able to confirm the previous
%experimental and numerical observations related to the sedimentation
%of a cloud of particles in viscous fluids.  Moreover, a preliminary
%analysis of the simulation data revealed many interesting aspects of
%the problem that are not yet explored.  For example, we showed that
%the generated flow in sedimenting fibers is analogous to the vortex
%ring in $\reynolds \gg 1$.  This pattern has also been reported for
%sedimenting clouds of particles in inertially dominated regimes.  We
%are currently working on developing the statistical mechanical
%framework that puts these seemingly different regimes in a unified
%perspective.

For the results presented here, although the longer fibers ($L=0.4
R_0$) undergo moderate bending, their deformations were not large
enough to qualitatively change the overall dynamics and shape of the
sedimenting cloud. Another interesting direction of research would be
to bridge between the observed dynamics of sedimenting polymeric
liquids (very flexible fibers) \citep{Sostarecz2003} to the case
studied here by systematically varying the flexibility of the fibers
and their volume fraction.

%% To circumvent this, in previous numerical study of this problem the
%% \emph{Stokeslet} singularity is regularized by a minimum cutoff
%% distance \citep{Park2010}. We did not pursue such modifications.

\section{Conclusions\label{sec:conclusion}}
We presented a platform for the dynamic simulation of semi-flexible
fibrous assemblies immersed in a Stokesian fluid with applications to
cellular mechanics and suspension mechanics.  The hydrodynamic
interactions (HIs) of the fibers with the fluid and other bodies are
described by the non-local slender body theory while a second-kind
boundary integral formulation is used for other rigid bodies and the
outer confinement.  The deformations of the fibers under external
flows and forces are modeled using the Euler--Bernoulli beam theory.
We introduced several modifications with respect to the previous
numerical treatment of semi-flexible fibers \citep{Tornberg2004},
including the spectral representation of the fibers, modeling growth
and shrinking dynamics of the fibers (important in biological
settings), using fast summation techniques (fast multipole method) for
computing HIs, parallelization of the entire computation scheme, and
implicit treatment of the HIs of fibers and particles in the
time-stepping scheme.

Our numerical tests demonstrate that the method is spectrally accurate
in space.  Also, due to the use of semi-implicit time-stepping, the
number of points per-fiber was found to have no effect on the stable
time-step.  Finally as a result of the efficient preconditioning of
the linear system, GMRES converges to the solution (with relative
error of $10^{-5}$) with less than $16$ iterations, for a relatively
concentrated assembly of fibers ( $N_f=2048$).

We validated our method against the available theoretical and
experimental observations in three representative problems.  Through
this process, we found several new directions of research within each
problem that can be pursued using this numerical technique.  For
example, our simulation technique can be used to investigate the
problems that arise in the flow of suspensions in complex geometries,
such as shear-induced migration \citep{Frank2003, Moraczewski2005,
  Gao2008} and dynamics of flexible filaments in microfluidic channels
\citep{Autrusson2011}.

%to study two problems in cell mechanics and
%cell mechanics and one problem to demonstrate the more general
%applicability of the method.
%\editt{In the first problem, we used an active microrheology setup to compute
%the mobility and viscoelastic behavior of a spherical particle with a
%shell formed by a large number of radially attached fibers, within a
%spherical confinement.
%In the second problem, we used our simulation to study the effect of
%cell geometry on the positioning of pronuclear complex (PNC) in the
%prophase stage of cell division.
%For this we used three geometries with different type and extent of
%asymmetry and two different active mechanisms for moving the PNC.
%In the third problem, we looked at the sedimentation of a cloud of
%flexible fibers under gravity, a classical suspension mechanics
%problem.
%In this problem the many-body hydrodynamic interactions result in
%evolution of the cloud to a torus-like structure at large enough
%sedimentation distances.
%In short, our simulation}
Our method also enables us to compute actively induced cytoplasmic
flows in cellular processes.  For example, in a concurrent work on the
dynamics of pronuclear migration in cell division, we demonstrate that
each active mechanism involved in moving the PNC has a distinct flow
signature.  These generic features of the flow can differentiate
between potential active mechanisms.  Another important application
arises in microrheology measurements, where the motion of a
microscopic probe (or probes) is used to infer the rheology of the
medium \citep{Mofrad2009, Wirtz2009, Squires2010}.  Within our method
it is straightforward to compute the motion of the probe that occurs
due to external forcing and the motion of other bodies and a variety
of internal active mechanisms.  Such fundamental studies can form the
basis of microscopically informed coarse-grained theories for
describing the mechanics of cytoskeleton.

Our numerical method can be improved in several ways.  While thermal
fluctuations are typically negligible for microtubule filaments in
cytoskeleton, they play a key role in the dynamics of actin and
intermediate filaments as well as the dynamics of polymer chains.
Thus, to extend the applicability of the platform, thermal
fluctuations need to be included and we are currently working on this
feature.
%The main numerical difficulty here is to efficiently compute the
%square root of the mobility matrix, which is needed for computing the
%displacement induced by thermal forces \commentt[AR]{rephrase?}.
%While this has already been implemented in methods such as
%particle-mesh Ewald \citep{Banchio, SDS2005} a parallelized fast
%multipole method (FMM) code for this purpose is still not
%available.
Steric interactions of fibers were also not imposed in our
simulations, and fibers were allowed to cross one another.  While this
was a very rare event in our simulations, modeling such interactions
is required if our framework is to be applied to dense assemblies such
as the mitotic spindle structure, dense actin networks, and polymer
melts above their entanglement concentration.

%\editt{An inherent difficulty with boundary integral techniques is evaluating
%the flow near the surfaces and singularities.
%While we did not consider any special treatment for these conditions,
%there are efficient and accurate methods to compute the flow near the
%surfaces \protect\citep{}.}{}
In this work, we used the first-order backward Euler method for the
time stepping with satisfactory results in terms of accuracy and
stability.  In modeling of pronuclear migration in \pr{ssc:position},
the fibers grow, shrink and alternate between these states, through a
Poisson process.  As a result, high order backward difference time
stepping methods are unnecessarily complex and hinder using adaptive
time stepping.  High order methods based on spectral deferred
correction (SDC) are successfully used for simulation of vesicles in
2D using boundary integral methods with adaptive time stepping
\cite{Quaife2014b}.  Extension of such method to fiber suspensions is
an interesting research direction.

Although several studies demonstrate the importance of HIs in defining
the rheology of synthetic systems such as fiber pulps, and the
collective dynamics of fibers in cilia and flagellar motion and
biofilms \citep{Shelley2015}, these interactions between the fibers in
cytoskeletal assemblies have been almost entirely ignored.
The typical justification for this was that the long-range HIs beyond
the nearest neighbor distance is screened and the dynamics of each
individual fiber is dictated by its local interactions with the
neighboring fibers \citep{BM2014}.
Our example of viscoelastic behavior of a sphere with the shell
composed of radial microtubules shows that this picture is not
correct in all length scales.
We observed that at times longer than the elastic relaxation time of
the fibers, the structure can be modeled as a sphere with a porous
shell using Brinkman equation.
Ignoring the HIs between the MTs would result in a drag coefficient
that would linearly increase with the number of fibers, in strong
contradictions with the predictions of Brinkman equation and our
simulation results.
The forced oscillation simulation on the same structure also revealed
that the characteristic relaxation time of the porous particle is $25$
times shorter than the relaxation time of individual fibers which
shows that many-body HIs substantially modify the elastic properties
of the material as well.
Thus in attempting to model the HIs through coarse-grained
relationships, the underlying assumptions and their correctness in the
relevant time- and length-scale must be scrutinized.
In this regard, direct simulation provides a powerful tool to inform
the reduced models.
%

%%- refs
\bibliographystyle{alpha}
\bibliography{refs}

\newcommand{\etalchar}[1]{$^{#1}$}
\begin{thebibliography}{FLMD{\etalchar{+}}05}

\bibitem[AGL{\etalchar{+}}11]{Autrusson2011}
Nicolas Autrusson, Laura Guglielmini, Sigolene Lecuyer, Roberto Rusconi, and
  Howard~A Stone.
\newblock The shape of an elastic filament in a two-dimensional corner flow.
\newblock {\em Physics of Fluids (1994-present)}, 23(6):063602, 2011.

\bibitem[AKY77]{Adachi1977}
K~Adachi, S.~Kiriyama, and N.~Yoshioka.
\newblock {The behavior of a swarm of particles moving in a viscous fluid}.
\newblock {\em Chemical Engineering Science}, 33:115--121, 1977.

\bibitem[ARW95]{Ascher1995a}
UM~Ascher, SJ~Steven J.~SJ Ruuth, and Brian T. R.~BTR Wetton.
\newblock {Implicit-Explicit Methods for Time-Dependent Partial Differential
  Equations}.
\newblock {\em SIAM Journal on Numerical Analysis}, 32(3):797--823, June 1995.

\bibitem[BK01]{Bruno2001}
Oscar~P. Bruno and Leonid~A. Kunyansky.
\newblock {A Fast, High-Order Algorithm for the Solution of Surface Scattering
  Problems: Basic Implementation, Tests, and Applications}.
\newblock {\em Journal of Computational Physics}, 169(1):80--110, May 2001.

\bibitem[BK06]{BK2006}
A.~R. Bausch and K.~Kroy.
\newblock {A bottom-up approach to cell mechanics}.
\newblock {\em Nature Physics}, 2(4):231--238, 2006.

\bibitem[BM14]{BM2014}
C.~P. Broedersz and F.~C. MacKintosh.
\newblock {Modeling semiflexible polymer networks}.
\newblock {\em Rev. of Modern Phys.}, 86(3):995--1036, 2014.

\bibitem[Boy01]{boyd2001}
John~P Boyd.
\newblock {\em {Chebyshev and Fourier spectral methods}}.
\newblock Courier Corporation, 2001.

\bibitem[Bri47]{Brinkman1947}
H.~C. Brinkman.
\newblock {A calculation of the viscous force exerted by a flowing fluid on a
  dense swarm of particles}.
\newblock {\em Applied Science Res}, A1:27--34, 1947.

\bibitem[BS01]{BS2001}
Leif~E. Becker and Michael~J. Shelley.
\newblock Instability of elastic filaments in shear flow yields
  first-normal-stress differences.
\newblock {\em Physical Review Letters}, 87(19):198301, 2001.

\bibitem[BS02]{Butler2002}
Jason~E Butler and Eric~SG Shaqfeh.
\newblock Dynamic simulations of the inhomogeneous sedimentation of rigid
  fibres.
\newblock {\em Journal of Fluid Mechanics}, 468:205--237, 2002.

\bibitem[BTB03]{Bush2003}
John W.~M. Bush, B.~A. Thurber, and F.~Blanchette.
\newblock {Particle clouds in homogeneous and stratified environments}.
\newblock {\em Journal of Fluid Mechanics}, 489:29--54, 2003.

\bibitem[BWV14]{Barnett2014}
Alex Barnett, Bowei Wu, and SK~Veerapaneni.
\newblock { Spectrally-accurate quadratures for evaluation of layer potentials
  close to the boundary for the 2D Stokes and Laplace equations}.
\newblock {\em arXiv preprint arXiv:1410.2187}, pages 1--21, 2014.

\bibitem[BZK09]{Bommes2009a}
David Bommes, Henrik Zimmer, and Leif Kobbelt.
\newblock {Mixed-Integer Quadrangulation}.
\newblock In {\em ACM Transactions On Graphics (TOG)}, page~77, 2009.

\bibitem[CFM05]{Cortez2005}
Ricardo Cortez, Lisa Fauci, and Alexei Medovikov.
\newblock {The method of regularized Stokeslets in three dimensions: Analysis,
  validation, and application to helical swimming}.
\newblock {\em Physics of Fluids}, 17(3):1--14, 2005.

\bibitem[CH04]{Cowan2004}
Carrie~R Cowan and Anthony~a Hyman.
\newblock {Asymmetric cell division in C. elegans: cortical polarity and
  spindle positioning.}
\newblock {\em Annual review of cell and developmental biology}, 20:427--453,
  2004.

\bibitem[Cor01]{Cortez2001}
Ricardo Cortez.
\newblock The method of regularized stokeslets.
\newblock {\em SIAM Journal on Scientific Computing}, 23(4):1204--1225, 2001.

\bibitem[CRZ15]{corona2015}
Eduardo Corona, Abtin Rahimian, and Denis Zorin.
\newblock A tensor-train accelerated solver for integral equations in complex
  geometries, 2015.

\bibitem[CWG11]{Gompper2011}
Raghunath Chelakkot, Roland~G Winkler, and Gerhard Gompper.
\newblock Semiflexible polymer conformation, distribution and migration in
  microcapillary flows.
\newblock {\em Journal of Physics: Condensed Matter}, 23(18):184117, 2011.

\bibitem[DA85]{Davis1985}
Robert~H Davis and Andreas Acrivos.
\newblock {Sedimentation of noncolloidal particles at low reynolds numbers}.
\newblock {\em Annual Review of Fluid Mechanics}, 2(17):91--118, 1985.

\bibitem[DKWD04]{Dahl2004}
Kris~Noel Dahl, Samuel~M Kahn, Katherine~L Wilson, and Dennis~E Discher.
\newblock {The nuclear envelope lamina network has elasticity and a
  compressibility limit suggestive of a molecular shock absorber.}
\newblock {\em Journal of cell science}, 117(20):4779--4786, 2004.

\bibitem[DM97]{Desai1997}
Arshad Desai and Timothy~J Mitchison.
\newblock Microtubule polymerization dynamics.
\newblock {\em Annual review of cell and developmental biology}, 13(1):83--117,
  1997.

\bibitem[DY97]{Dogterom1997}
M~Dogterom and B~Yurke.
\newblock {Measurement of the force-velocity relation for growing
  microtubules.}
\newblock {\em Science}, 278(5339):856--860, 1997.

\bibitem[FAWM03]{Frank2003}
Martin Frank, Douglas Anderson, Eric~R. Weeks, and Jeffrey~F. Morris.
\newblock {Particle migration in pressure-driven flow of a Brownian
  suspension}.
\newblock {\em J. Fluid Mech.}, 493:363--378, 2003.

\bibitem[FLMD{\etalchar{+}}05]{Flores2005}
Heather Flores, Edgar Lobaton, Stefan M{\'e}ndez-Diez, Svetlana Tlupova, and
  Ricardo Cortez.
\newblock A study of bacterial flagellar bundling.
\newblock {\em Bulletin of mathematical biology}, 67(1):137--168, 2005.

\bibitem[FM10]{Fletcher2010}
Daniel~a Fletcher and R~Dyche Mullins.
\newblock {Cell mechanics and the cytoskeleton.}
\newblock {\em Nature}, 463(7280):485--492, 2010.

\bibitem[FYJ15]{Fu2015}
Szu-Pei Fu, Y.-N. Young, and Shidong Jiang.
\newblock Efficient brownian dynamics simulation of dna molecules with
  hydrodynamic interactions in linear flows.
\newblock {\em Phys. Rev. E}, 91:063008, Jun 2015.

\bibitem[GG08]{Gao2008}
C.~Gao and J.~F. Gilchrist.
\newblock {Shear-induced particle migration in one-, two-, and
  three-dimensional flows}.
\newblock {\em Phys. Rev. E}, 77(2):025301, 2008.

\bibitem[GGSH01]{Grill2001}
S.~W. Grill, P.~G\"{o}nczy, E.~H. Stelzer, and A.~Hyman.
\newblock {Polarity controls forces governing asymmetric spindle positioning in
  the Caenorhabditis elegans embryo.}
\newblock {\em Nature}, 409(6820):630--633, 2001.

\bibitem[GH10]{Guazzelli2010}
E~Lisabeth Guazzelli and John Hinch.
\newblock {Fluctuations and instability in sedimentation}.
\newblock {\em Annual Review of Fluid Mechanics}, 2010.

\bibitem[GK10]{Gilden2010}
Julia Gilden and Matthew~F. Krummel.
\newblock {Control of cortical rigidity by the cytoskeleton: Emerging roles for
  septins}.
\newblock {\em Cytoskeleton}, 67(8):477--486, 2010.

\bibitem[G{\"o}t00]{gotz2000}
Thomas G{\"o}tz.
\newblock {\em Interactions of fibers and flow: asymptotics, theory and
  numerics}.
\newblock dissertation., 2000.

\bibitem[GPW98]{Goldstein1998a}
Raymond~E. Goldstein, Thomas~R. Powers, and Chris~H. Wiggins.
\newblock {The Viscous Nonlinear Dynamics of Twist and Writhe}.
\newblock {\em Physical Review Letters}, 80(23):9, 1998.

\bibitem[Gro12]{Groot2012}
Robert~D Groot.
\newblock {How to impose stick boundary conditions in coarse-grained
  hydrodynamics of Brownian colloids and semi-flexible fiber rheology.}
\newblock {\em The Journal of chemical physics}, 136(6):064901, 2012.

\bibitem[Gua01]{Guazzelli2001}
Elisabeth Guazzelli.
\newblock {Evolution of particle-velocity correlations in sedimentation}.
\newblock {\em Physics of Fluids}, 13(6):1537--1540, 2001.

\bibitem[HB83]{Happel1983}
J~Happel and H~Brenner.
\newblock {\em {Low Reynolds number hydrodynamics, 1965}}.
\newblock Martinus Nijhoff Publishers, 1983.

\bibitem[HF96]{Higdon1996}
J.~J.~L. Higdon and G.~D. Ford.
\newblock {Permeability of three-dimensional models of fibrous porous media}.
\newblock {\em Journal of Fluid Mechanics}, 308(1):341, 1996.

\bibitem[HG99]{Herzhaft1999}
Benjamin Herzhaft and {\'E}lisabeth Guazzelli.
\newblock Experimental study of the sedimentation of dilute and semi-dilute
  suspensions of fibres.
\newblock {\em Journal of Fluid Mechanics}, 384:133--158, 1999.

\bibitem[HGMS96]{Herzhaft1996}
Benjamin Herzhaft, {\'E}lisabeth Guazzelli, Michael~B Mackaplow, and Eric~SG
  Shaqfeh.
\newblock Experimental investigation of the sedimentation of a dilute fiber
  suspension.
\newblock {\em Physical review letters}, 77(2):290, 1996.

\bibitem[HO08]{Helsing2008}
Johan Helsing and Rikard Ojala.
\newblock {On the evaluation of layer potentials close to their sources☆}.
\newblock {\em Journal of Computational Physics}, 227(5):2899--2921, February
  2008.

\bibitem[How01]{Howard-Book}
Jonathon Howard.
\newblock {\em Mechanics of Motor Proteins and the Cytoskeleton}.
\newblock Sinauer Associates, Inc, first edition, 2001.

\bibitem[HW87]{HW1987}
Anthony~A Hyman and John~G White.
\newblock Determination of cell division axes in the early embryogenesis of
  caenorhabditis elegans.
\newblock {\em The J. Cell Biol.}, 105(5):2123--2135, 1987.

\bibitem[JDD03]{Janson2003}
Marcel~E. Janson, Mathilde~E. {De Dood}, and Marileen Dogterom.
\newblock {Dynamic instability of microtubules is regulated by force}.
\newblock {\em Journal of Cell Biology}, 161(6):1029--1034, 2003.

\bibitem[JF01]{Joung2001}
C~G Joung and X~J Fan.
\newblock {Direct simulation of flexible fibers}.
\newblock {\em J. Non-{N}ewtonian Fluid. Mech.}, 99:1--36, 2001.

\bibitem[Joh80]{Johnson1980}
Robert~E. Johnson.
\newblock {An improved slender-body theory for Stokes flow}.
\newblock {\em Journal of Fluid Mechanics}, 99(02):411, 1980.

\bibitem[KBGO13]{Klockner2012}
A.~Kl\"{o}ckner, A.~H. Barnett, L.~Greengard, and M.~O'Neil.
\newblock {Quadrature by Expansion: A New Method for the Evaluation of Layer
  Potentials}.
\newblock {\em Journal of Computational Physics}, 0(3):332--349, July 2013.

\bibitem[KK89]{Karrila1989}
Seppo~J. Karrila and Sangtae Kim.
\newblock {Integral Equations of the Second Kind for Stokes Flow: Direct
  Solution for Physical Variables and Removal of Inherent Accuracy
  Limitations}.
\newblock {\em Chemical Engineering Communications}, 82(1):123--161, 1989.

\bibitem[KK11]{Kimura2011}
Kenji Kimura and Akatsuki Kimura.
\newblock {Intracellular organelles mediate cytoplasmic pulling force for
  centrosome centration in the Caenorhabditis elegans early embryo.}
\newblock {\em Proceedings of the National Academy of Sciences of the United
  States of America}, 108(1):137--142, 2011.

\bibitem[KO05]{Kimura2005}
Akatsuki Kimura and Shuichi Onami.
\newblock {Computer simulations and image processing reveal length-dependent
  pulling force as the primary mechanism for C. elegans male pronuclear
  migration}.
\newblock {\em Developmental Cell}, 8(5):765--775, 2005.

\bibitem[KO07]{Kimura2007}
Akatsuki Kimura and Shuichi Onami.
\newblock {Local cortical pulling-force repression switches centrosomal
  centration and posterior displacement in C. elegans}.
\newblock {\em Journal of Cell Biology}, 179(7):1347--1354, 2007.

\bibitem[KR76]{Keller1976}
JB~Keller and SI~Rubinow.
\newblock {Slender-body theory for slow viscous flow}.
\newblock {\em Journal of Fluid Mechanics}, 75:705--714, 1976.

\bibitem[Kre99]{Kress1989}
Rainer Kress.
\newblock {\em {Linear integral equations}}, volume~82.
\newblock Springer-Verlag, New York,, 1999.

\bibitem[KS11]{Keaveny2011}
Eric~E. Keaveny and Michael~J. Shelley.
\newblock {Applying a second-kind boundary integral equation for surface
  tractions in Stokes flow}.
\newblock {\em Journal of Computational Physics}, 230(5):2141--2159, March
  2011.

\bibitem[LMSS13]{LMSS2013}
Lei Li, Harishankar Manikantan, David Saintillan, and Saverio. Spagnolie.
\newblock {The sedimentation of flexible filaments}.
\newblock {\em Journal of Fluid Mechanics}, 735:705--736, 2013.

\bibitem[LP04]{Lim2004}
Sookkyung Lim and Charles~S. Peskin.
\newblock {Simulations of the Whirling Instability by the Immersed Boundary
  Method}.
\newblock {\em SIAM Journal on Scientific Computing}, 25(6):2066--2083, 2004.

\bibitem[LS14]{LS2014}
A.~Linder and M.~Shelley.
\newblock Elastic fibers in flows.
\newblock {\em Fluid-Structure Interactions in Low Reynolds Number Flows,
  edited by C. Duprat and HA Stone (RSC Publishing, 2015)}, 2014.

\bibitem[MB14]{Malhotra2014}
Dhairya Malhotra and George Biros.
\newblock {A distributed memory fast multipole method for volume potentials}.
\newblock {\em users.ices.utexas.edu}, V(212), 2014.

\bibitem[MB15]{Malhotra2015}
Dhairya Malhotra and George Biros.
\newblock {PVFMM: A Parallel Kernel Independent FMM for Particle and Volume
  Potentials}.
\newblock {\em Communications in Computational Physics}, 18(03):808--830, 2015.

\bibitem[MBC11]{Minc2011}
Nicolas Minc, David Burgess, and Fred Chang.
\newblock {Influence of cell geometry on division-plane positioning}.
\newblock {\em Cell}, 144(3):414--426, 2011.

\bibitem[McN13]{McNally2013a}
Francis~J. McNally.
\newblock {Mechanisms of spindle positioning}.
\newblock {\em Journal of Cell Biology}, 200(2):131--140, 2013.

\bibitem[MGB05]{Metzger2005}
Bloen Metzger, Elisabeth Guazzelli, and Jason~E. Butler.
\newblock {Large-scale streamers in the sedimentation of a dilute fiber
  suspension}.
\newblock {\em Physical Review Letters}, 95(16):1--4, 2005.

\bibitem[MLSS14]{Manikantan2014}
Harishankar Manikantan, Lei Li, Saverio~E Spagnolie, and David Saintillan.
\newblock The instability of a sedimenting suspension of weakly flexible
  fibres.
\newblock {\em Journal of Fluid Mechanics}, 756:935--964, 2014.

\bibitem[MMNS01]{Machu2001}
Gunther Machu, Walter Meile, Ludwig~C. Nitsche, and Uwe Schaflinger.
\newblock {Coalescence, torus formation and breakup of sedimenting drops:
  experiments and computer simulations}.
\newblock {\em Journal of Fluid Mechanics}, 447(April 2000):299--336, 2001.

\bibitem[MNG07]{Metzger2007}
Bloen Metzger, Maxime Nicolas, and \'{E}lisabeth Guazzelli.
\newblock {Falling clouds of particles in viscous fluids}.
\newblock {\em Journal of Fluid Mechanics}, 580:283, 2007.

\bibitem[MNMG87]{Masliyah1987}
Jacob~H. Masliyah, Graham Neale, K.~Malysa, and Theodorus {G.M. Van De Ven}.
\newblock {Creeping flow over a composite sphere: Solid core with porous
  shell}.
\newblock {\em Chemical Engineering Science}, 42(2):245--253, 1987.

\bibitem[Mof09]{Mofrad2009}
Mohammad~R.K. Mofrad.
\newblock {Rheology of the Cytoskeleton}.
\newblock {\em Ann. Rev. of Fluid Mech.}, 41:433--453, 2009.

\bibitem[MS16]{Manikantan2016}
Harishankar Manikantan and David Saintillan.
\newblock Effect of flexibility on the growth of concentration fluctuations in
  a suspension of sedimenting fibers: Particle simulations.
\newblock {\em Physics of Fluids (1994-present)}, 28(1):013303, 2016.

\bibitem[MTS05]{Moraczewski2005}
T.~Moraczewski, H.~Y. Tang, and N.~C. Shapley.
\newblock {Flow of a concentrated suspension through an abrupt axisymmetric
  expansion measured by nuclear magnetic resonance imaging}.
\newblock {\em Journal of rheology}, 49(6):1409--1428, 2005.

\bibitem[NF07]{Nedelec2007}
Francois Nedelec and Dietrich Foethke.
\newblock {Collective Langevin dynamics of flexible cytoskeletal fibers}.
\newblock {\em New J. Physics}, 9:0--24, 2007.

\bibitem[NL05]{Nguyen2005}
Nhan-Quyen Nguyen and Anthony J.~C. Ladd.
\newblock {Sedimentation of hard-sphere suspensions at low Reynolds number}.
\newblock {\em Journal of Fluid Mechanics}, 525:73--104, 2005.

\bibitem[NM15]{NM2015}
E~Nazockdast and Jeffrey~F. Morris.
\newblock Active microrheology of colloidal suspensions: simulation and
  microstructural theory.
\newblock {\em Submitted to Journal of Rheology}, 2015.

\bibitem[NRNS15]{NRNS2015}
Ehssan Nazockdast, Abtin Rahimian, Daniel Needleman, and Michael Shelley.
\newblock Cytoplasmic flows as signatures for the mechanics of mitotic
  positioning.
\newblock {\em Preprint}, pages 1--26, 2015.

\bibitem[NS01]{Nedelec2001}
Fran{\c{c}}ois N{\'e}d{\'e}lec and Thomas Surrey.
\newblock Dynamics of microtubule aster formation by motor complexes.
\newblock {\em Comptes Rendus de l'Acad{\'e}mie des Sciences-Series
  IV-Physics-Astrophysics}, 2(6):841--847, 2001.

\bibitem[OLC13]{Olson2013}
Sarah~D. Olson, Sookkyung Lim, and Ricardo Cortez.
\newblock {Modeling the dynamics of an elastic rod with intrinsic curvature and
  twist using a regularized Stokes formulation}.
\newblock {\em Journal of Computational Physics}, 238:169--187, 2013.

\bibitem[OT14]{Ojala2014}
Rikard Ojala and Anna-Karin Tornberg.
\newblock {An accurate integral equation method for simulating multi-phase
  Stokes flow}.
\newblock pages 1--22, April 2014.

\bibitem[PM87]{Power1987}
Henry Power and Guillermo Miranda.
\newblock {Second Kind Integral Equation Formulation of Stokes' Flows Past a
  Particle of Arbitrary Shape}.
\newblock {\em SIAM Journal on Applied Mathematics}, 47(4):689, 1987.

\bibitem[PMGB10]{Park2010}
Joontaek Park, Bloen Metzger, \'{E}lisabeth Guazzelli, and Jason~E. Butler.
\newblock {A cloud of rigid fibres sedimenting in a viscous fluid}.
\newblock {\em Journal of Fluid Mechanics}, 648:351, 2010.

\bibitem[POO93]{Peskin1993}
C.~S. Peskin, G.~M. Odell, and G.~F. Oster.
\newblock Cellular motions and thermal fluctuations: the brownian ratchet.
\newblock {\em Biophysical journal}, 65(1):316, 1993.

\bibitem[Poz92]{pozrikidis1992}
C~Pozrikidis.
\newblock {\em {Boundary Integral and Singularity Methods for Linearized
  Viscous Flow}}.
\newblock Cambridge University Press, Cambridge, 1992.

\bibitem[Poz01]{Pozrikidis2001c}
C~Pozrikidis.
\newblock {Interfacial dynamics for Stokes flow}.
\newblock {\em Journal of Computational Physics}, 169:250301, 2001.

\bibitem[PR08]{Park2008}
Dae~Hwi Park and Lesilee~S. Rose.
\newblock {Dynamic localization of LIN-5 and GPR-1/2 to cortical force
  generation domains during spindle positioning}.
\newblock {\em Developmental Biology}, 315(1):42--54, 2008.

\bibitem[PRH{\etalchar{+}}06]{Pecreaux2006a}
Jacques Pecreaux, Jens-christian Ro, Anthony~A Hyman, Karsten Kruse, Frank Ju,
  Stephan~W Grill, and Jonathon Howard.
\newblock {Spindle sscillations during asymmetric cell division require a
  threshold number of active cortical force generators}.
\newblock {\em Current Biology}, pages 2111--2122, 2006.

\bibitem[QB14]{Quaife2014b}
Bryan Quaife and George Biros.
\newblock {High-order adaptive time stepping for vesicle suspensions with
  viscosity contrast}.
\newblock pages 1--11, 2014.

\bibitem[Ram01]{Ramaswamy2001}
Sriram Ramaswamy.
\newblock {Issues in the statistical mechanics of steady sedimentation}.
\newblock {\em Advances in Physics}, 50(3):297--341, 2001.

\bibitem[RGGG98]{Reinsch1998}
Sigrid Reinsch, Pierre G\"{o}nczy, P~Gonczy, and Pierre G\"{o}nczy.
\newblock {Mechanisms of nuclear positioning}.
\newblock {\em Journal of cell science}, 2295(6820):2283--2295, August 1998.

\bibitem[RK97]{Ross1997}
Russell~F Ross and Daniel~J Klingenberg.
\newblock {Dynamic simulation of flexible fibers composed of linked rigid
  bodies}.
\newblock {\em J. Chem. Phys.}, 106:2949--2960, 1997.

\bibitem[RVZB15]{Rahimian2015}
Abtin Rahimian, Shravan~K. Veerapaneni, Denis Zorin, and George Biros.
\newblock {Boundary integral method for the flow of vesicles with viscosity
  contrast in three dimensions}.
\newblock {\em Journal of Computational Physics}, 298:766--786, 2015.

\bibitem[Saa03]{saad2003}
Yousef Saad.
\newblock {\em {Iterative methods for sparse linear systems}}.
\newblock Siam, 2003.

\bibitem[SB03]{Sostarecz2003}
Michael~C. Sostarecz and Andrew Belmonte.
\newblock {Motion and shape of a viscoelastic drop falling through a viscous
  fluid}.
\newblock {\em Journal of Fluid Mechanics}, 497:235--252, 2003.

\bibitem[SD09]{Siller2009}
Karsten~H Siller and Chris~Q Doe.
\newblock {Spindle orientation during asymmetric cell division.}
\newblock {\em Nature cell biology}, 11(4):365--374, 2009.

\bibitem[SDS05]{SDS2005}
D~Saintillan, E~Darve, and E~S~G Shaqfeh.
\newblock {A smooth particle-mesh Ewald algorithm for Stokes suspension
  simulations: The sedimentation of fibers}.
\newblock {\em Physics of Fluids}, 17(3):33301, 2005.

\bibitem[SG98]{Stockie1998}
John~M Stockie and Sheldon~I Green.
\newblock Simulating the motion of flexible pulp fibres using the immersed
  boundary method.
\newblock {\em Journal of Computational Physics}, 147(1):147--165, 1998.

\bibitem[She16]{Shelley2015}
M.~Shelley.
\newblock The dynamics of microtubule/motor-protein assemblies in biology and
  physics.
\newblock {\em Ann. Rev. Fluid Mech.}, 48(1):null, 2016.

\bibitem[SLW01]{Serge2001}
P~N Segr\'{e}, F~Liu, and D~A Weitz.
\newblock {An effective gravitational temperature for sedimentation}.
\newblock {\em Nature}, 409, 2001.

\bibitem[SM10]{Squires2010}
Todd~M. Squires and Thomas~G. Mason.
\newblock {Fluid mechanics of microrheology}.
\newblock {\em Ann. Rev. of Fluid Mech.}, 42:413--438, January 2010.

\bibitem[Smi09]{Smith2010}
D.~J. Smith.
\newblock A boundary element regularized stokeslet method applied to cilia- and
  flagella-driven flow.
\newblock {\em Proceedings of the Royal Society of London A: Mathematical,
  Physical and Engineering Sciences}, 465(2112):3605--3626, 2009.

\bibitem[SMPS11]{Shinar2011}
Tamar Shinar, Miyeko Mana, Fabio Piano, and Michael~J Shelley.
\newblock {A model of cytoplasmically driven microtubule-based motion in the
  single-celled Caenorhabditis elegans embryo.}
\newblock {\em Proc. Natl. Acad. Sci. USA}, 108(26):10508--10513, 2011.

\bibitem[SU00]{Shelley2000}
Michael~J. Shelley and Tetsuji Ueda.
\newblock {The Stokesian hydrodynamics of flexing, stretching filaments}.
\newblock {\em Physica D: Nonlinear Phenomena}, 146(1-4):221--245, 2000.

\bibitem[TA99]{Tong1999}
P~Tong and B.~J. Ackerson.
\newblock {Analogies between colloidal sedimentation and turbulent convection
  at high Prandtl numbers}.
\newblock {\em Physical Review E}, 58(6):6931--6934, 1999.

\bibitem[TG06]{Tornberg2006}
Anna-Karin Tornberg and Katarina Gustavsson.
\newblock A numerical method for simulations of rigid fiber suspensions.
\newblock {\em Journal of Computational Physics}, 215(1):172--196, 2006.

\bibitem[THD{\etalchar{+}}02]{Tsou2002}
Meng-Fu~Bryan Tsou, Adam Hayashi, Leah~R DeBella, Garth McGrath, and Lesilee~S
  Rose.
\newblock {LET-99 determines spindle position and is asymmetrically enriched in
  response to PAR polarity cues in C. elegans embryos.}
\newblock {\em Development (Cambridge, England)}, 129(19):4469--4481, 2002.

\bibitem[Tre00]{trefethen2000}
Lloyd~N Trefethen.
\newblock {\em {Spectral methods in MATLAB}}, volume~10.
\newblock Siam, 2000.

\bibitem[TS04]{Tornberg2004}
Anna-Karin Tornberg and Michael~J. Shelley.
\newblock {Simulating the dynamics and interactions of flexible fibers in
  Stokes flows}.
\newblock {\em Journal of Computational Physics}, 196(1):8--40, May 2004.

\bibitem[vDTMD00]{VanDoorn2000}
G~S van Doorn, C~Tanase, B~M Mulder, and M~Dogterom.
\newblock {On the stall force for growing microtubules.}
\newblock {\em European biophysics journal : EBJ}, 29(1):2--6, 2000.

\bibitem[WA10]{Wu2010}
Jingshu Wu and Cyrus~K. Aidun.
\newblock {A method for direct simulation of flexible fiber suspensions using
  lattice Boltzmann equation with external boundary force}.
\newblock {\em International Journal of Multiphase Flow}, 36(3):202--209, 2010.

\bibitem[WHD01]{WHD2001}
Torsten Wittmann, Anthony Hyman, and Arshad Desai.
\newblock {The spindle : a dynamic assembly of microtubules and motors}.
\newblock {\em Nature Cell Biology}, 3:E28--E34, 2001.

\bibitem[Wir09]{Wirtz2009}
Denis Wirtz.
\newblock {Particle-tracking microrheology of living cells: principles and
  applications.}
\newblock {\em Ann. Rev. of biophysics}, 38:301--326, 2009.

\bibitem[WN96]{Wen1996}
F~Wen and A~Nacamuli.
\newblock {The effect of the Rayleigh number on a particle cloud}.
\newblock {\em Hydrodynamics}, pages 1275--1280, 1996.

\bibitem[WR75]{Weisenberg1975}
Richard~C Weisenberg and Arline~C Rosenfeld.
\newblock In vitro polymerization of microtubules into asters and spindles in
  homogenates of surf clam eggs.
\newblock {\em The Journal of cell biology}, 64(1):146--158, 1975.

\bibitem[WS15]{WS2015}
Jeffrey~K Wiens and John~M Stockie.
\newblock Simulating flexible fiber suspensions using a scalable immersed
  boundary algorithm.
\newblock {\em Computer Methods in Applied Mechanics and Engineering},
  290:1--18, 2015.

\bibitem[YBZ04]{ying2004kernel}
Lexing Ying, George Biros, and Denis Zorin.
\newblock A kernel-independent adaptive fast multipole algorithm in two and
  three dimensions.
\newblock {\em Journal of Computational Physics}, 196(2):591--626, 2004.

\bibitem[YBZ06]{Ying2006}
Lexing Ying, George Biros, and Denis Zorin.
\newblock {A high-order 3D boundary integral equation solver for elliptic PDEs
  in smooth domains}.
\newblock {\em Journal of Computational Physics}, 219(1):247--275, 2006.

\bibitem[YM95]{Yamamoto1995}
S~Yamamoto and T~Matsuoka.
\newblock Dynamic simulation of fiber suspensions in shear flow.
\newblock {\em J. Chem. Phys.}, 102:2254--2260, 1995.

\bibitem[YM96]{Yamamoto1996}
S~Yamamoto and T~Matsuoka.
\newblock Dynamic simulation of microstructure and rheology of fiber
  suspensions.
\newblock {\em Polym. Eng. Sci.}, 36(19):2396--2403, 1996.

\bibitem[ZD00]{Zinchenko2000}
AZ~Zinchenko and RH~Davis.
\newblock {An efficient algorithm for hydrodynamical interaction of many
  deformable drops}.
\newblock {\em Journal of Computational Physics}, 157(2):539--587, 2000.

\end{thebibliography}

\end{document}